\documentclass{elsarticle}
\usepackage{moreverb} 
\usepackage[colorlinks,bookmarksopen,bookmarksnumbered,citecolor=red,urlcolor=red,pdfauthor=author]{hyperref}
\usepackage{amsfonts,amsmath,amssymb}
\usepackage{amsthm}
\newtheorem{theorem}{Theorem}
\newtheorem{proposition}[theorem]{Proposition}
\newtheorem{corollary}[theorem]{Corollary}
\DeclareMathOperator{\sig}{sigmoid}

\usepackage[section]{placeins}
\usepackage{subcaption}
\usepackage{mathrsfs}
\usepackage{graphicx}
\usepackage{algorithm}
\usepackage{algpseudocode}
\usepackage{cleveref}
\usepackage{fancyhdr}
\usepackage[export]{adjustbox} 
\usepackage{enumerate} 
\usepackage{comment}
\usepackage[margin=2.5cm]{geometry}
\usepackage{bm}
\usepackage{booktabs}
\usepackage{multirow}
\usepackage{threeparttable}
\usepackage{siunitx}
\newcommand{\panellabel}[1]{\parbox[t]{0.32\textwidth}{\centering #1}}

\usepackage{pgfplots}
\pgfplotsset{compat=1.13}

\begin{document}

\begin{frontmatter}
    \renewcommand{\thefootnote}{\fnsymbol{footnotemark}}

    \fancypagestyle{plain}{%
        \fancyhf{} 
        \fancyhead[RO,RE]{\thepage} 
    }
    \title{Gappy probabilistic manifold decomposition for nonlinear field reconstruction}
    \author[lab1]{Qihan Feng}
   \author[lab2]{Jiaming Guo}
   \author[lab1]{Jiarun Meng}
   \author[lab1]{Dunhui Xiao\corref{cor1}}
   \cortext[cor1]{Corresponding author. Email: \texttt{xiaodunhui@tongji.edu.cn} (Dunhui Xiao)}
   \address[lab1]{School of Mathematical Sciences, Key Laboratory of Intelligent
      Computing and Applications (Ministry of Education), Tongji University, Shanghai
      200092, CHINA}
   \address[lab2]{Shanghai Research Institute for Intelligent Autonomous
      Systems, Tongji University, Shanghai 201210, CHINA}

    \begin{abstract}
        This paper proposes gappy probabilistic manifold decomposition (Gappy PMD), a nonlinear method for reconstructing high-dimensional fields from extremely sparse measurements. Gappy PMD reconstructs the field on the nonlinear manifold learned by probabilistic manifold decomposition (PMD). We further propose a differentiable point selection method
        for reduced-order model (ROM)-based field reconstruction (DPS).
        Using differentiable meshless interpolation within the ROM-based reconstruction framework, DPS makes the full-field reconstruction error differentiable with respect to the sampling locations and directly optimizes these locations. In addition, a theoretical error analysis for Gappy PMD is also given. It splits the squared reconstruction error into two orthogonal parts: one normal to the reconstruction manifold and the other induced by sparse sampling and observation noise. Under a stability condition on the sampling operator, this error vanishes with the PMD approximation error and the noise. The Gappy PMD is evaluated on three numerical test cases: flow past a cylinder, lid-driven cavity flow, and backward-facing step flow. For the same reduced dimension and sampling points, Gappy PMD attains mean relative $L^2$ errors one to two orders of magnitude below Gappy POD. Optimizing the sampling points with DPSR further improves reconstruction accuracy and robustness.

    \end{abstract}

\end{frontmatter}

\section{Introduction}
\vspace{-2pt}

Approximating a field variable over a spatial domain from sparse measurements is a fundamental problem in computational science and engineering \cite{loiseau2018sparse,kefal2021smoothed,kefal2024particle}. In practice, measurements are available at only a small fraction of the spatial domain and may be corrupted by noise \cite{willcox2006unsteady,introini2023stabilization}. Consequently, direct full-field reconstruction is underdetermined and sensitive to measurement errors \cite{manohar2018data,binev2018greedy,maday2015generalized,peherstorfer2020stability}. Spatial interpolation techniques such as radial basis function interpolation \cite{barba2010global} and Gaussian process regression \cite{gunes2006gappy} are commonly employed to this end, but reconstruction accuracy deteriorates sharply as the measurements become extremely sparse \cite{fukami2021machine,santos2023senseiver}. 

This limitation motivates the incorporation of prior information about the underlying field structure into the reconstruction. Such prior information can be extracted from a collection of field snapshots available offline. Deep learning methods along this line have advanced rapidly, with representative approaches including convolutional neural networks \cite{fukami2019super,fukami2021machine}, Senseiver \cite{santos2023senseiver}, and RecFNO \cite{zhao2024recfno}. These methods learn mappings from sparse measurements to complete fields and have demonstrated strong reconstruction performance across diverse settings. However, these methods typically require extensive training datasets, and the learned mapping is encoded implicitly in the network parameters with no explicit connection to the underlying physical structure \cite{fukami2021spatiotemporal, Erichson2020, Callaham2019, Dubois2022}. In fact, the snapshot data of many physical phenomena concentrate near a low-dimensional manifold \cite{holmes2012turbulence, Berkooz1993, fu2023nonlinear, guo2025probabilistic}. A compact representation can therefore be constructed offline. The online reconstruction then reduces to determining a small number of coordinates from the sparse measurements, and the resulting problem is far better constrained than recovering the full field directly. Gappy proper orthogonal decomposition (Gappy POD) \cite{everson1995karhunen} builds on this idea, exploiting the low-rank structure captured by POD to reconstruct missing data from sparse measurements. In Gappy POD, the field is represented by a small number of dominant POD modes, and the full state is recovered by estimating their modal coefficients from the available sparse measurements. It has since been widely studied and applied to field reconstruction problems \cite{gunes2006gappy,raben2012adaptive,peherstorfer2020stability,nekkanti2023gappy}.

Such methods rely on linear subspace approximation and achieve satisfactory reconstruction accuracy whenever the field data are well represented by a small number of POD modes. However, for strongly nonlinear or convection-dominated phenomena, the Kolmogorov $n$-width decays slowly, so linear methods must retain a large number of modes to reach acceptable accuracy \cite{hesthaven2026nonlinear, greif2019decay}. This severely degrades computational efficiency. Although adaptive approximation spaces and localized bases can alleviate this limitation \cite{peherstorfer2020model,amsallem2012nonlinear}, they do not eliminate the fundamental geometric restriction of linear subspace approximation \cite{Lee2020,Kim2022}.

To break through this bottleneck, nonlinear dimensionality reduction has been introduced into field reconstruction from sparse measurements. Martin-Linares et al.\ \cite{martinlinares2023gappy} proposed the gappy diffusion maps (Gappy DMAPs) method, which employs diffusion maps \cite{coifman2006diffusion} to extract the low-dimensional manifold structure of the data. It represents the intrinsic geometry with far fewer coordinates than POD modes while maintaining good reconstruction accuracy even with few or poorly located sampling points. Kim et al.\ \cite{kim2024gappy} proposed the gappy autoencoder (Gappy AE) method, which replaces the linear subspace representation used in Gappy POD with a nonlinear manifold learned by an autoencoder. Gappy AE has shown improved reconstruction accuracy over Gappy POD for problems with slowly decaying Kolmogorov $n$-width.

The performance of a field reconstruction method depends not only on the reconstruction algorithm itself but also critically on how the sampling points are selected \cite{manohar2018data}. For Gappy POD, the point selection procedure of the discrete empirical interpolation method (DEIM) \cite{chaturantabut2010nonlinear} can be used to determine sampling locations greedily. QDEIM \cite{drmac2016new} replaces the greedy procedure with a pivoted QR factorization and provides stronger theoretical guarantees. For reconstruction based on deep learning, Marcato et al.\ \cite{marcato2024differentiable} proposed a differentiable sensor placement strategy on the Senseiver model via bilinear interpolation. Liu et al.\ \cite{liu2025dspo} proposed DSPO, which optimizes sensor locations for field reconstruction with deep learning models using gradients of the reconstruction error. 
However, point selection for field reconstruction based on nonlinear low-dimensional representations has received much less attention.

Recently, Guo and Xiao \cite{guo2025probabilistic} proposed probabilistic manifold decomposition (PMD), a nonlinear model reduction method. PMD combines linear reduction via singular value decomposition (SVD) with a probabilistic embedding of the nonlinear residual manifold. This enables PMD to capture nonlinear structures that are difficult to represent efficiently in a linear subspace. Building on PMD, we propose Gappy PMD for field reconstruction from sparse measurements. The PMD representation is constructed offline from field snapshots. Online reconstruction then estimates the linear coefficients and manifold coordinates by solving a nonlinear least-squares problem using the available measurements. 


Based on this framework, we further propose a differentiable point selection method for reduced-order model (ROM)-based field reconstruction (DPS). Related methods such as DEIM and QDEIM can be adopted for point selection but rely on proxy criteria rather than the reconstruction error itself \cite{chaturantabut2010nonlinear,drmac2016new,manohar2018data}. DPS optimizes sampling locations directly with respect to the reconstruction error. To this end, we relax the discrete indices into continuous coordinates and evaluate the measurements via smooth interpolation. The reconstruction error then varies smoothly with the sampling locations and can be minimized directly. The approach is developed here for Gappy PMD and carries over to the classical Gappy POD.

We analyze the reconstruction error of Gappy PMD and establish an error identity. It decomposes the error into a component orthogonal to the local tangent space of the reconstruction manifold and a term governed by measurement noise and sampling stability on this tangent space. We further prove that the reconstruction error vanishes as the PMD approximation error and the measurement noise tend to zero, provided that the sampling operator remains stable on the reconstruction set. In the absence of the nonlinear lift mapping, the analysis reduces to the classical Gappy POD error estimate, thereby covering both linear and nonlinear reconstruction within a unified framework.

The remainder of this paper is organized as follows. Section~\ref{sec:gappy-pmd-DPS} presents the PMD representation, the Gappy PMD reconstruction method, and the differentiable point selection strategy. Section~\ref{sec:theoretical-analysis} presents a theoretical analysis of the reconstruction error of Gappy PMD, establishing error estimates and analyzing its convergence. Section~\ref{sec:numerical-results} validates the proposed methods through three numerical experiments: the two-dimensional flow past a cylinder, the lid-driven cavity flow, and the three-dimensional backward-facing step flow. Section~\ref{Conclusions} draws conclusions and outlines directions for future research.

\section{Gappy PMD and differentiable point selection}\label{sec:gappy-pmd-DPS}
The PMD representation is reviewed and adapted to the field reconstruction problem considered here. Gappy PMD is then formulated to recover the linear coefficients and manifold coordinates from sparse measurements. A differentiable point selection strategy is therefore introduced to optimize the measurement locations.

\subsection{Probabilistic manifold decomposition}\label{sec:pmd}
PMD represents a field by its projection onto a truncated SVD basis, with the projection residual modeled on a nonlinear manifold \cite{guo2025probabilistic}.

\subsubsection{Linear reduction}\label{sec:pmd-linear}
Let the training snapshot matrix be $X = [x_1, \dots, x_m] \in \mathbb{R}^{n \times m}$, where $n$ denotes the number of spatial degrees of freedom and $m$ denotes the number of snapshots. The training mean is
\begin{equation}
\bar{x} = \frac{1}{m}\sum_{j=1}^{m} x_j .
\label{eq:pmd-mean}
\end{equation}
Each snapshot is centered by subtracting the mean in Eq.~\eqref{eq:pmd-mean}. Denoting the centered snapshot matrix by $\widetilde{X} = [x_1 - \bar{x},\, \dots,\, x_m - \bar{x}]$, its economy singular value decomposition is
\begin{equation}
\widetilde{X} = Q\,\Sigma\,V^\top,
\label{eq:pmd-svd}
\end{equation}
where $Q \in \mathbb{R}^{n \times m}$ and $V \in \mathbb{R}^{m \times m}$ collect the left and right singular vectors as columns, and $\Sigma = \operatorname{diag}(\sigma_1, \dots, \sigma_m)$ contains singular values ordered as $\sigma_1 \geq \sigma_2 \geq \cdots \geq \sigma_m \geq 0$. The fraction of energy captured by the leading $r$ modes is
\begin{equation}
E(r) = \frac{\displaystyle\sum_{i=1}^{r} \sigma_i^{2}}{\displaystyle\sum_{i=1}^{m} \sigma_i^{2}},
\label{eq:pmd-energy}
\end{equation}
and $r$ is chosen so that $E(r) \geq 1-\delta$ for a prescribed tolerance $\delta$. The first $r$ columns of $Q$ form an orthonormal basis $\Phi \in \mathbb{R}^{n \times r}$. The projection coefficients and the residual of each snapshot are
\begin{equation}
a_j = \Phi^\top(x_j - \bar{x}), \qquad
u_j = x_j - \bar{x} - \Phi\,a_j .
\label{eq:pmd-linear-residual}
\end{equation}
By the Eckart--Young--Mirsky theorem \cite{eckart1936approximation,mirsky1960symmetric}, the $r$-dimensional subspace spanned by $\Phi$ minimizes the total reconstruction error $\sum_j \lVert u_j \rVert_2^{2}$ over all linear subspaces of the same dimension. Equivalently, the truncated SVD in Eq.~\eqref{eq:pmd-svd} provides the optimal rank-$r$ approximation in the Frobenius norm.

The approximation accuracy of any linear subspace is nevertheless limited by the Kolmogorov $n$-width of the solution manifold. For convection-dominated and other strongly nonlinear problems, this width can decay slowly with $r$, so that increasing the number of linear modes gives diminishing reductions in the residual. PMD addresses this limitation by modeling the residuals in Eq.~\eqref{eq:pmd-linear-residual} on a nonlinear manifold rather than by extending the linear basis.

\subsubsection{Residual manifold coordinates}\label{sec:pmd-residual}
After the linear reduction of section~\ref{sec:pmd-linear}, the residuals $\{u_j\}_{j=1}^{m}$ retain nonlinear features outside the span of $\Phi$. PMD represents these residuals through a low-dimensional nonlinear manifold $\mathcal{M}$ embedded in $\mathbb{R}^n$. The geometry of $\mathcal{M}$ is learned by constructing a Markov process on the residual samples. The spectral decomposition of its transition matrix then provides the coordinates of the probabilistic manifold \cite{guo2025probabilistic}.

The transition probabilities are built from pairwise similarities that measure proximity along $\mathcal{M}$. A Gaussian similarity is defined as
\begin{equation}
W_{ij} = \exp\!\left(-\frac{d_{\mathcal{M}}^{\,2}(u_i,\,u_j)}{\varepsilon}\right), \qquad i,j = 1, \dots, m,
\label{eq:pmd-similarity}
\end{equation}
where $\varepsilon > 0$ is a bandwidth parameter and $d_{\mathcal{M}}(u_i,u_j)$ denotes the geodesic distance between $u_i$ and $u_j$ on $\mathcal{M}$. Since the true geodesic distance is not available in closed form, it is approximated by the Floyd--Warshall algorithm on a $k_{\mathrm{nn}}$-nearest-neighbor graph constructed from the residual samples \cite{guo2025probabilistic}:
\begin{equation}
d_{\mathcal{M}}^{(s)}(u_i,u_j) = \min\!\left(d_{\mathcal{M}}^{(s-1)}(u_i,u_j),\; d_{\mathcal{M}}^{(s-1)}(u_i,u_s) + d_{\mathcal{M}}^{(s-1)}(u_s,u_j)\right),
\label{eq:pmd-geodesic-update}
\end{equation}
initialized with $\lVert u_i-u_j\rVert_2$ for neighboring pairs and $\infty$ otherwise. Using geodesic rather than Euclidean distance allows Eq.~\eqref{eq:pmd-similarity} to reflect separation along a curved residual manifold.

When residual samples are not uniformly distributed over $\mathcal{M}$, the row sums of $W$ are affected by sampling density. To reduce this bias, a density normalization is applied \cite{coifman2006diffusion}. Let $\rho_i = \sum_j W_{ij}$ estimate the local density at $u_i$. The corrected similarity is
\begin{equation}
\widetilde{W}_{ij} = \frac{W_{ij}}{\rho_i\,\rho_j} .
\label{eq:pmd-density-normalization}
\end{equation}
The corrected similarity matrix is then normalized by rows to define the Markov transition matrix
\begin{equation}
P(i,j) = \frac{\widetilde{W}_{ij}}{\displaystyle\sum_{k=1}^{m}\widetilde{W}_{ik}}, \qquad i,j = 1, \dots, m .
\label{eq:pmd-transition}
\end{equation}
Each entry $P(i,j)$ is the probability that a random walker at residual $u_i$ transitions to $u_j$ in one step. Nearby residuals on $\mathcal{M}$ have high transition probabilities, whereas distant residuals have negligible probability. The $t$-step transition matrix $P^t$ therefore captures larger-scale geometric connectivity by accumulating all paths of length $t$.

The manifold coordinates are extracted from the spectral structure of the transition matrix. The eigenvalue problem is
\begin{equation}
P\,\varphi_k = \lambda_k\,\varphi_k, \qquad
1 = \lambda_0 \geq \lambda_1 \geq \cdots \geq \lambda_{m-1} \geq 0 .
\label{eq:pmd-eigenproblem}
\end{equation}
The leading eigenvalue $\lambda_0=1$ corresponds to the stationary distribution and carries no geometric coordinate. Since $P$ and $P^t$ share eigenvectors, the $t$-step transition satisfies $P^t\varphi_k = \lambda_k^t\varphi_k$. The first $r_1$ nontrivial eigenvectors, weighted by the corresponding eigenvalue powers, define the residual manifold coordinates
\begin{equation}
\alpha_j = \left(\lambda_1^{\,t}\varphi_1(j),\; \dots,\; \lambda_{r_1}^{\,t}\varphi_{r_1}(j)\right) \in \mathbb{R}^{r_1} .
\label{eq:pmd-manifold-coordinates}
\end{equation}
Here $\varphi_k(j)$ denotes the $j$-th component of the $k$-th eigenvector. The factor $\lambda_k^t$ suppresses rapidly mixing modes associated with small eigenvalues, so the embedding retains the dominant geometric degrees of freedom. The truncation dimension $r_1$ is selected at a spectral gap where $\lambda_{r_1}$ is well separated from $\lambda_{r_1+1}$. Together with the linear coefficients in Eq.~\eqref{eq:pmd-linear-residual}, each training snapshot is parameterized by the low-dimensional PMD coordinates $(a_j,\alpha_j) \in \mathbb{R}^{r+r_1}$.

\subsubsection{Lift mapping and PMD representation}\label{sec:pmd-lift}
The construction above maps each training snapshot to the low-dimensional PMD coordinates $(a_j,\alpha_j)$. To recover the full state, PMD requires a lift mapping from the manifold coordinates back to the original state space. In the present formulation, the linear reconstruction $\bar{x}+\Phi a$ is already given by the SVD basis, so the lift mapping targets only the nonlinear residuals by constructing a function $f$ such that $f(\xi_j) \approx u_j$.

In the original PMD formulation, the lift mapping takes only the manifold coordinates $\alpha_j$ as input \cite{guo2025probabilistic}. However, the nonlinear PMD embedding is not guaranteed to be injective, so distinct residuals may be mapped to nearby or overlapping regions of the coordinate space \cite{otto2022inadequacy}. To reduce this ambiguity, we augment the input with the linear coefficients $a_j$. Since $\alpha_j$ is constructed from the residuals alone, $a_j$ carries complementary information about the linear component of the snapshot. Because $\alpha_j$ and $a_j$ generally live on different scales, each component is standardized to zero mean and unit variance:
\begin{equation}
\widetilde{\alpha}_j = \frac{\alpha_j - \bar{\alpha}}{\sigma_{\alpha}}, \qquad
\widetilde{a}_j = \frac{a_j - \bar{a}}{\sigma_a},
\label{eq:pmd-standardization}
\end{equation}
where $\bar{\cdot}$ and $\sigma_{\cdot}$ denote componentwise training means and standard deviations. The standardized components are concatenated to form the joint feature
\begin{equation}
\xi_j = [\widetilde{\alpha}_j,\; \widetilde{a}_j] \in \mathbb{R}^{d}, \qquad d = r_1 + r .
\label{eq:pmd-joint-feature}
\end{equation}
Standardization in Eq.~\eqref{eq:pmd-standardization} ensures that all coordinates contribute on comparable scales to the kernel evaluation.

The lift mapping is constructed by kernel ridge regression (KRR) with Gaussian radial basis functions. Given the training pairs $\{(\xi_j,u_j)\}_{j=1}^{m}$, the mapping is
\begin{equation}
f(\xi) = \sum_{j=1}^{m} c_j \exp\!\left(-\frac{\lVert \xi-\xi_j\rVert_2^2}{\varepsilon_f}\right),
\label{eq:pmd-lift-map}
\end{equation}
where $c_j \in \mathbb{R}^{n}$ are coefficient vectors and $\varepsilon_f>0$ is the kernel bandwidth. The bandwidth is a fixed multiple $c_{\varepsilon}>0$ of the median of pairwise squared distances among the training features,
\begin{equation}
\varepsilon_f = c_{\varepsilon}\,\operatorname{median}\!\left\{\lVert \xi_i-\xi_j\rVert_2^2 : 1 \leq i < j \leq m\right\},
\label{eq:pmd-kernel-bandwidth}
\end{equation}
which adapts the kernel scale to the data. Let $K \in \mathbb{R}^{m \times m}$ denote the Gram matrix with entries $K_{ij}=\exp(-\lVert \xi_i-\xi_j\rVert_2^2/\varepsilon_f)$, and let $U=[u_1,\dots,u_m]\in\mathbb{R}^{n\times m}$ collect the training residuals. The coefficient matrix $C=[c_1,\dots,c_m]^\top\in\mathbb{R}^{m\times n}$ is obtained by solving
\begin{equation}
(K+\lambda_f I)\,C = U^\top,
\label{eq:pmd-krr-system}
\end{equation}
where $\lambda_f>0$ is a regularization parameter. The complete PMD representation of a snapshot is
\begin{equation}
\widehat{x} = \bar{x} + \Phi\,a + f(\xi), \qquad \xi = \xi(\alpha,a) .
\label{eq:pmd-representation}
\end{equation}
The linear term in Eq.~\eqref{eq:pmd-representation} supplies the optimal rank-$r$ approximation of Eq.~\eqref{eq:pmd-svd}. The nonlinear term recovers what no $r$-dimensional linear subspace can reach. Accuracy therefore improves by adding manifold coordinates rather than by enlarging $r$.

\subsection{Gappy PMD}\label{sec:gappy-pmd}
Gappy PMD recovers the full state $x\in\mathbb{R}^{n}$ from $q \ll n$ sparse measurements by estimating the coordinates of the PMD representation in Eq.~\eqref{eq:pmd-representation}. Let $\mathcal{S}=\{i_1,\dots,i_q\}\subset\{1,\dots,n\}$ denote the sampled spatial indices, and let $S\in\mathbb{R}^{q\times n}$ be the corresponding sampling matrix. The measurements are
\begin{equation}
y = S\,x + \eta \in \mathbb{R}^{q},
\label{eq:gappy-observation}
\end{equation}
where $\eta$ accounts for measurement noise. The goal is to determine $(\alpha,a)$ from the measurements in Eq.~\eqref{eq:gappy-observation} and then recover the full state using Eq.~\eqref{eq:pmd-representation}.

In the full-field setting of section~\ref{sec:pmd-lift}, the lift mapping reconstructs the residual from $(\alpha,a)$ alone. With only $q$ measurements, however, the measurements $y$ and a candidate coefficient vector $a$ determine the nonlinear residual at the sampling points:
\begin{equation}
b(a;y) = y - S\bar{x} - S\Phi\,a \in \mathbb{R}^{q} .
\label{eq:gappy-observed-residual}
\end{equation}
Including this sampled residual as additional input allows the lift mapping to extrapolate from $q$ known local residual values to the unsampled mesh nodes. The joint feature in Eq.~\eqref{eq:pmd-joint-feature} is therefore extended to
\begin{equation}
\xi = [\widetilde{\alpha},\; \widetilde{a},\; \widetilde{b}] \in \mathbb{R}^{d_{\xi}}, \qquad d_{\xi}=r_1+r+q,
\label{eq:gappy-augmented-feature}
\end{equation}
where each component is standardized as in Eq.~\eqref{eq:pmd-standardization}. The lift mapping is retrained with this augmented input, using $b_j=S u_j$ computed from the training residuals. Because $b$ depends on the sampling configuration $\mathcal{S}$, the lift mapping must be retrained whenever $\mathcal{S}$ changes. The unknowns remain $(\alpha,a)$ only, since $b$ is determined by $a$ and $y$ through Eq.~\eqref{eq:gappy-observed-residual}. The complete offline construction is summarized in Algorithm~\ref{alg:offline}.

\begin{algorithm}[t]
\caption{Offline PMD training.}
\label{alg:offline}
\begin{algorithmic}[1]
\Procedure{GappyPMDTrain}{$X\in\mathbb{R}^{n\times m},\;
S\in\mathbb{R}^{q\times n},\; r,\; r_1,\; k_{\mathrm{nn}},\;
\varepsilon,\; t,\; c_\varepsilon,\; \lambda_f$}

\State Compute the snapshot mean
       $\displaystyle
       \bar{x}=\frac{1}{m}\sum_{j=1}^{m}x_j$

\State Compute the leading $r$ left singular vectors $\Phi$ of
       $[x_1-\bar{x},\,\dots,\,x_m-\bar{x}]$

\State For $j=1,\dots,m$, define
       $a_j=\Phi^\top(x_j-\bar{x})$ and
       $u_j=x_j-\bar{x}-\Phi a_j$

\State Build a $k_{\mathrm{nn}}$-nearest-neighbor graph on
       $\{u_1,\dots,u_m\}$ and approximate the geodesic distances
       $\{d_{\mathcal{M}}(u_i,u_j)\}$ by the Floyd--Warshall recursion
       \hfill \eqref{eq:pmd-geodesic-update}

\State Define
       $W_{ij}=
       \exp\!\bigl(-d_{\mathcal{M}}^2(u_i,u_j)/\varepsilon\bigr)$,
       $i,j=1,\dots,m$
       \hfill \eqref{eq:pmd-similarity}

\State Compute
       $\rho_i=\sum_j W_{ij}$ and define
       $\widetilde{W}_{ij}=W_{ij}/(\rho_i\rho_j)$ and
       $P(i,j)=\widetilde{W}_{ij}/\sum_l \widetilde{W}_{il}$
       \hfill
       \eqref{eq:pmd-density-normalization}--\eqref{eq:pmd-transition}

\State Compute the leading $r_1$ nontrivial eigenpairs
       $\{\lambda_k,\varphi_k\}_{k=1}^{r_1}$ of $P$
       \hfill \eqref{eq:pmd-eigenproblem}

\State For $j=1,\dots,m$, define
       $\alpha_j=
       \bigl(
       \lambda_1^t\varphi_1(j),\,\dots,\,
       \lambda_{r_1}^t\varphi_{r_1}(j)
       \bigr)$
       \hfill \eqref{eq:pmd-manifold-coordinates}

\State For $j=1,\dots,m$, compute
       $b_j=S u_j$

\State Standardize each component:
       $\widetilde{\alpha}_j=
       (\alpha_j-\bar{\alpha})/\sigma_\alpha$,
       $\widetilde{a}_j=
       (a_j-\bar{a})/\sigma_a$, and
       $\widetilde{b}_j=
       (b_j-\bar{b})/\sigma_b$
       \hfill \eqref{eq:pmd-standardization}

\State For $j=1,\dots,m$, define the augmented feature
       $\xi_j=
       [\widetilde{\alpha}_j;\,
        \widetilde{a}_j;\,
        \widetilde{b}_j]
       \in\mathbb{R}^{d_{\xi}}$
       \hfill \eqref{eq:gappy-augmented-feature}

\State Set the kernel bandwidth
       $\displaystyle
       \varepsilon_f=
       c_\varepsilon\,
       \operatorname{median}
       \left\{
       \lVert\xi_i-\xi_j\rVert_2^2:
       1\le i<j\le m
       \right\}$
       \hfill \eqref{eq:pmd-kernel-bandwidth}

\State Form the kernel matrix
       $K_{ij}=
       \exp\!\left(
       -\frac{\lVert\xi_i-\xi_j\rVert_2^2}{\varepsilon_f}
       \right)$,
       $i,j=1,\dots,m$

\State Solve
       $(K+\lambda_f I)C=U^\top$
       for
       $C=[c_1,\dots,c_m]^\top
       \in\mathbb{R}^{m\times n}$
       \hfill \eqref{eq:pmd-krr-system}

\State \Return
       $\bar{x},\,\Phi,\,
       \{a_j,\alpha_j\}_{j=1}^{m},\,
       \{\xi_j,c_j\}_{j=1}^{m},\,
       \varepsilon_f,\,
       (\bar{\alpha},\sigma_\alpha,\,
        \bar{a},\sigma_a,\,
        \bar{b},\sigma_b)$

\EndProcedure
\end{algorithmic}
\end{algorithm}

The reconstruction problem is to find coordinates whose PMD prediction is consistent with the measurements. The unknowns are collected into a normalized vector $\theta=[\widetilde{\alpha},\,\widetilde{a}]\in\mathbb{R}^{d}$, where $d=r_1+r$. The physical coordinates $\alpha(\theta)$ and $a(\theta)$ are recovered by inverting Eq.~\eqref{eq:pmd-standardization}. At any candidate $\theta$, Eq.~\eqref{eq:gappy-observed-residual} gives the sampled residual $b$, and Eq.~\eqref{eq:gappy-augmented-feature} gives the augmented feature $\xi(\theta)$. The measurement mismatch is
\begin{equation}
R(\theta) = b(a(\theta);y) - S\,f(\xi(\theta)) \in \mathbb{R}^{q} .
\label{eq:gappy-residual}
\end{equation}
The coordinates are estimated by solving the nonlinear least-squares problem
\begin{equation}
\theta^{\star} = \arg\min_{\theta\in\Theta}\,\lVert R(\theta)\rVert_2^2,
\label{eq:gappy-nls}
\end{equation}
where $\Theta=[\theta_L,\theta_U]$ is a box constraint obtained by expanding the range of training coordinates by a fixed margin. This constraint confines the search to the region where the lift mapping has been trained. In all experiments, $q=d+2$, so the residual system in Eq.~\eqref{eq:gappy-residual} is moderately overdetermined.

Because $R$ depends nonlinearly on $\theta$ through the kernel expansion in Eq.~\eqref{eq:pmd-lift-map}, an iterative solver is required. The linear coefficients are initialized by a regularized least-squares fit with regularization parameter $\lambda_a>0$,
\begin{equation}
a_0 = (\Phi_s^\top\Phi_s + \lambda_a I)^{-1}\Phi_s^\top(y-S\bar{x}), \qquad \Phi_s=S\Phi,
\label{eq:gappy-linear-init}
\end{equation}
and the manifold coordinates are initialized by a soft nearest-neighbor average over the training set $\{\alpha_j\}$. Let $\mathcal{N}_{k_\alpha}(\widetilde{a}_0)$ denote the $k_\alpha$ nearest training samples to $\widetilde{a}_0$ in the space of standardized linear coefficients. The initial manifold coordinates are
\begin{equation}
\alpha_0 = \sum_{j\in\mathcal{N}_{k_\alpha}(\widetilde{a}_0)} w_j\,\alpha_j, \qquad
w_j = \frac{\exp\!\left(-\lVert \widetilde{a}_0-\widetilde{a}_j\rVert_2^2/\tau\right)}{\displaystyle\sum_{l\in\mathcal{N}_{k_\alpha}(\widetilde{a}_0)}\exp\!\left(-\lVert \widetilde{a}_0-\widetilde{a}_l\rVert_2^2/\tau\right)},
\label{eq:gappy-alpha-init}
\end{equation}
where $\tau>0$ is a temperature parameter. Starting from $\theta_0=[\widetilde{\alpha}_0,\,\widetilde{a}_0]$, Eq.~\eqref{eq:gappy-nls} is solved by a trust-region reflective method \cite{Branch1999} with an analytic Jacobian. Differentiating Eq.~\eqref{eq:gappy-residual} with respect to $\theta$ gives
\begin{equation}
J(\theta) = \frac{\partial R}{\partial\theta}
= -S\Phi\,\frac{\partial a}{\partial\theta}
  - S\,\frac{\partial f}{\partial\xi}\,\frac{\partial\xi}{\partial\theta},
\label{eq:gappy-jacobian}
\end{equation}
where both terms can be evaluated in closed form because the linear prediction and the kernel expansion are smooth functions of $\theta$. The full state is recovered from the solution as
\begin{equation}
\widehat{x} = \bar{x} + \Phi\,a(\theta^{\star}) + f(\xi(\theta^{\star})) .
\label{eq:gappy-reconstruction}
\end{equation}
The online reconstruction is summarized in Algorithm~\ref{alg:online}.

\begin{algorithm}[t]
\caption{Online Gappy PMD reconstruction.}
\label{alg:online}
\begin{algorithmic}[1]
\Procedure{GappyPMDReconstruct}{$y,\,S,\,\Theta,\,\lambda_a,\,\tau,\,k_\alpha$}

\State Set $\Phi_s=S\Phi$

\State Compute the linear initial guess
$\displaystyle
 a_0=(\Phi_s^\top\Phi_s+\lambda_a I)^{-1}
 \Phi_s^\top(y-S\bar{x})$
\hfill \eqref{eq:gappy-linear-init}

\State Standardize $a_0$ to $\widetilde{a}_0$ and determine
$\mathcal{N}_{k_\alpha}(\widetilde{a}_0)$

\State Compute the manifold initial guess
$\displaystyle
 \alpha_0=\sum_{j\in\mathcal{N}_{k_\alpha}(\widetilde{a}_0)}
 w_j\alpha_j$,
where
$\displaystyle
 w_j=
 \frac{
 \exp\!\left(
 -\lVert\widetilde{a}_0-\widetilde{a}_j\rVert_2^2/\tau
 \right)}
 {
 \sum_{l\in\mathcal{N}_{k_\alpha}(\widetilde{a}_0)}
 \exp\!\left(
 -\lVert\widetilde{a}_0-\widetilde{a}_l\rVert_2^2/\tau
 \right)}
$
\hfill \eqref{eq:gappy-alpha-init}

\State Form the initial optimization variable
$\theta_0=[\widetilde{\alpha}_0,\widetilde{a}_0]\in\Theta$
from $(\alpha_0,a_0)$

\State For $\theta\in\Theta$, recover
$\alpha(\theta)$ and $a(\theta)$ and compute
$\displaystyle
 b(a(\theta);y)
 =
 y-S\bar{x}-S\Phi a(\theta)$
\hfill \eqref{eq:gappy-observed-residual}

\State Form the augmented feature
$\displaystyle
 \xi(\theta)=
 [\widetilde{\alpha},\,
  \widetilde{a},\,
  \widetilde{b}(\theta)]$
\hfill \eqref{eq:gappy-augmented-feature}

\State Define the measurement mismatch
$\displaystyle
 R(\theta)
 =
 b(a(\theta);y)-S f(\xi(\theta))$
\hfill \eqref{eq:gappy-residual}

\State Compute the analytic Jacobian
$\displaystyle
 J(\theta)
 =
 -S\Phi\frac{\partial a}{\partial\theta}
 -S\frac{\partial f}{\partial\xi}
   \frac{\partial\xi}{\partial\theta}$
\hfill \eqref{eq:gappy-jacobian}

\State Starting from $\theta_0$, solve
$\displaystyle
 \theta^\star
 =
 \arg\min_{\theta\in\Theta}
 \lVert R(\theta)\rVert_2^2$
by the trust-region reflective method
\hfill \eqref{eq:gappy-nls}

\State Compute the full-state reconstruction
$\displaystyle
 \widehat{x}
 =
 \bar{x}
 +\Phi a(\theta^\star)
 +f(\xi(\theta^\star))$
\hfill \eqref{eq:gappy-reconstruction}

\State \Return $\widehat{x}$

\EndProcedure
\end{algorithmic}
\end{algorithm}
       
\subsection{Differentiable point selection for ROM-based field reconstruction}\label{sec:DPS}

The accuracy of the Gappy PMD reconstruction in section~\ref{sec:gappy-pmd} is sensitive to the sampling set $\mathcal{S}$.  Classical point selection methods such as DEIM and QDEIM choose discrete mesh indices by algebraic criteria related to the reduced basis. These criteria improve numerical stability, but do not directly minimize the reconstruction error. However, directly optimizing sampling locations from the reconstruction error is difficult on a discrete mesh. DPS addresses this issue by relaxing the selected mesh indices to continuous sampling coordinates and then differentiating the resulting reconstruction error with respect to those coordinates.

Let $\{g_1,\ldots,g_n\}\subset\mathbb{R}^{d_g}$ be the physical mesh nodes, and let $\{\bar g_1,\ldots,\bar g_n\}\subset[0,1]^{d_g}$ denote the corresponding coordinates after componentwise normalization. In DPS, the \(q\) sampling points are parameterized by normalized coordinates \(\boldsymbol{\zeta}=(\zeta_1,\ldots,\zeta_q)\). To avoid constrained optimization over the box $[0,1]^{d_g}$, each point is parameterized by an unconstrained vector $\ell_i\in\mathbb{R}^{d_g}$ through
\begin{equation}
\zeta_i=\sig(\ell_i):=\bigl(1+e^{-\ell_i}\bigr)^{-1}\in(0,1)^{d_g},
\qquad i=1,\ldots,q .
\label{eq:DPS-continuous-coordinates}
\end{equation}
Thus the optimization is carried out over $\ell$, while $\zeta=\sig(\ell)$ remains inside the computational domain.

When $\zeta_i$ is allowed to vary continuously, it may fall between mesh nodes, so the discrete sampling matrix cannot be used directly. We therefore replace the discrete sampling operation by a meshless interpolation operator $\mathcal{I}_{\zeta}$. For a mesh function $v\in\mathbb{R}^n$, the value sampled at $\zeta_i$ is approximated from a local stencil $\mathcal{N}_i\subset\{1,\ldots,n\}$ of nearby mesh nodes:
\begin{equation}
\bigl[\mathcal{I}_{\zeta}v\bigr]_i
=
\sum_{j\in\mathcal{N}_i}
\omega_{ij}(\zeta_i)\,v_j,
\qquad i=1,\ldots,q .
\label{eq:DPS-interpolation}
\end{equation}
In this work the weights are computed from a radial basis function-generated finite difference (RBF-FD) interpolation system \cite{fornberg2015primer,shankar2015radial}. For the $i$-th sampling point, let $\phi$ be a radial basis function and let $\{p_1,\ldots,p_{N_p}\}$ be a polynomial basis of total degree at most $p_{\max}$. The weights $\omega^{(i)}=\{\omega_{ij}\}_{j\in\mathcal{N}_i}$ are obtained from
\begin{equation}
\begin{bmatrix}
A^{(i)} & B^{(i)}\\
(B^{(i)})^\top & 0
\end{bmatrix}
\begin{bmatrix}
\omega^{(i)}\\
\nu^{(i)}
\end{bmatrix}
=
\begin{bmatrix}
r^{(i)}\\
s^{(i)}
\end{bmatrix},
\label{eq:DPS-rbffd-system}
\end{equation}
where $A^{(i)}_{jj'}=\phi(\lVert \bar g_j-\bar g_{j'}\rVert_2)$, $B^{(i)}_{jl}=p_l(\bar g_j-\zeta_i)$, $r^{(i)}_j=\phi(\lVert \zeta_i-\bar g_j\rVert_2)$, $s^{(i)}_l=p_l(0)$, and $\nu^{(i)}$ contains the Lagrange multipliers that enforce exact reproduction of polynomials up to degree \(p_{\max}\). For each candidate $\zeta$, Eq.~\eqref{eq:DPS-rbffd-system} defines the weights used in Eq.~\eqref{eq:DPS-interpolation}. When differentiating with respect to the sampling coordinates, each local stencil $\mathcal{N}_i$ is held fixed. Within this fixed stencil, the interpolation weights vary smoothly with $\zeta$. Therefore, $\mathcal{I}_{\zeta}v$ is differentiable with respect to the sampling coordinates. In the experiments we use the cubic radial basis function $\phi(r)=r^3$, polynomial reproduction up to $p_{\max}=2$, and stencil size $k_s=40$.

Using $\mathcal{I}_{\zeta}$, we express the Gappy PMD reconstruction in terms of the continuous sampling coordinates. Each occurrence of the discrete sampling matrix~$S$ is replaced
by~$\mathcal{I}_{\zeta}$. Section~\ref{sec:gappy-pmd} constructs the augmented feature by using the sampled residuals $b_j=Su_j$ for the training snapshots and $b(a;y)=y-S\bar{x}-S\Phi a$ for a new measurement $y$. After replacing $S$ by $\mathcal{I}_{\zeta}$, the training residuals are first sampled as $b_j(\zeta)=\mathcal{I}_{\zeta}u_j$ and standardized as $\widetilde{b}_j(\zeta)=\bigl(b_j(\zeta)-\bar b(\zeta)\bigr)/\sigma_b(\zeta)$. These standardized residuals are then appended to the reduced coordinates to define the training features and Gram matrix,
\begin{equation}
\begin{alignedat}{2}
\xi_j(\zeta)
&=\bigl[\widetilde{\alpha}_j,\,
\widetilde{a}_j,\,
\widetilde{b}_j(\zeta)\bigr],\qquad&
K_{ij}(\zeta)
&=\exp\!\bigl(-\lVert \xi_i(\zeta)-\xi_j(\zeta)\rVert_2^2/\varepsilon_f\bigr).
\end{alignedat}
\label{eq:DPS-coordinate-krr}
\end{equation}
The corresponding KRR coefficients and lift mapping are then given by
\begin{equation}
\bigl(K(\zeta)+\lambda_f I\bigr)C(\zeta)=U^\top,\qquad
f_{\zeta}(\xi)=C(\zeta)^\top\kappa_{\zeta}(\xi).
\label{eq:DPS-coordinate-lift}
\end{equation}
Consequently, with the interpolation stencils fixed as described above, the assembled lift mapping $f_{\zeta}$ is differentiable with respect to~$\zeta$.

The online reconstruction map for a target snapshot is formulated as a differentiable function of the sampling coordinates~$\zeta$. This map is defined through the sparse measurements generated by $\mathcal{I}_{\zeta}$ and the PMD coordinates recovered from these measurements. For a target snapshot $x$, the sparse measurement is $y_\zeta=\mathcal{I}_\zeta x$. Given a candidate PMD coordinate $\theta=[\widetilde{\alpha},\widetilde{a}]\in\Theta$, the sampled residual and feature used to evaluate the lift are
\begin{equation}
b_{\zeta}(a(\theta);x)
=
\mathcal{I}_{\zeta}
\bigl(x-\bar{x}-\Phi a(\theta)\bigr),
\qquad
\xi_{\zeta}(\theta;x)
=
\bigl[
\widetilde{\alpha}(\theta),\,
\widetilde{a}(\theta),\,
\widetilde{b}_{\zeta}(a(\theta);x)
\bigr],
\label{eq:DPS-query-feature}
\end{equation}
which are Eqs.~\eqref{eq:gappy-observed-residual} and~\eqref{eq:gappy-augmented-feature} with $S$ replaced by $\mathcal{I}_{\zeta}$. The PMD coordinate $\theta_{\zeta}(x)$ is obtained by solving the online Gappy PMD nonlinear least-squares problem for the sparse measurement $y_\zeta$. This solve is represented by a fixed number $T$ of damped Gauss--Newton steps initialized as in Eqs.~\eqref{eq:gappy-linear-init} and~\eqref{eq:gappy-alpha-init}. If
\begin{equation}
R_{\zeta}(\theta;x)
=
b_{\zeta}(a(\theta);x)
-\mathcal{I}_{\zeta}
f_{\zeta}\!\left(\xi_{\zeta}(\theta;x)\right),
\label{eq:DPS-coordinate-residual}
\end{equation}
and $J_{\zeta}(\theta;x)=\partial R_{\zeta}(\theta;x)/\partial\theta$, one damped Gauss--Newton step solves
\begin{equation}
\left(
J_{\zeta}^{\top}J_{\zeta}
+\mu\,\operatorname{diag}(J_{\zeta}^{\top}J_{\zeta})
\right)
\Delta_{\zeta}^{(t)}
=
J_{\zeta}^{\top}R_{\zeta},
\qquad
\theta_{\zeta}^{(t+1)}(x)
=
\Pi_{\Theta}\!\left(
\theta_{\zeta}^{(t)}(x)
-\Delta_{\zeta}^{(t)}
\right),
\qquad t=0,\ldots,T-1 .
\label{eq:DPS-gn-update}
\end{equation}
where $R_{\zeta}$ and $J_{\zeta}$ are evaluated at $\theta=\theta_{\zeta}^{(t)}(x)$, and $\Pi_{\Theta}$ denotes projection onto $\Theta$. After $T$ steps, the recovered PMD coordinate and the corresponding full-state reconstruction are
\begin{equation}
\theta_{\zeta}(x)=\theta_{\zeta}^{(T)}(x),
\qquad
\mathcal{R}(\zeta;x)
=
\bar{x}
+\Phi\,a\bigl(\theta_{\zeta}(x)\bigr)
+f_{\zeta}\!\left(\xi_{\zeta}\bigl(\theta_{\zeta}(x);x\bigr)\right).
\label{eq:DPS-coordinate-map}
\end{equation}
This finite-step formulation makes the dependence of $\mathcal{R}(\zeta;x)$ on $\zeta$ explicit and allows it to be differentiated within the fixed-stencil computational graph.

The resulting differentiable reconstruction map provides a direct offline objective for choosing the sampling coordinates. Let $\mathcal{V}$ be a separate selection set that is not used in the offline PMD construction. DPS minimizes the relative reconstruction error
\begin{equation}
\mathcal{L}(\zeta)
=
\frac{1}{|\mathcal{V}|}
\sum_{x\in\mathcal{V}}
\frac{
\lVert \mathcal{R}(\zeta;x)-x\rVert_2
}{
\lVert x\rVert_2
}.
\label{eq:DPS-selection-loss}
\end{equation}
Gradients of $\mathcal{L}$ with respect to $\zeta$ are computed by automatic differentiation through the unrolled reconstruction map in Eq.~\eqref{eq:DPS-coordinate-map}. The objective is minimized with Adam using the coordinate parametrization in Eq.~\eqref{eq:DPS-continuous-coordinates}. After optimization, the continuous coordinates are snapped to distinct mesh nodes by solving a linear assignment problem, yielding the discrete sampling set $\mathcal{S}^{\star}$. The final PMD lift mapping is then rebuilt using the discrete sampling set $\mathcal{S}^{\star}$, and test snapshots are reconstructed by the standard Gappy PMD procedure in section~\ref{sec:gappy-pmd}.

DPS also applies to Gappy POD as a simpler special case. This corresponds to removing the nonlinear dimensionality reduction component from Gappy PMD. The reconstruction then reduces to $\widehat{x}=\bar{x}+\Phi a$, with the POD coefficients $a$ recovered online by linear least squares. The DPS procedure described above, including the continuous sampling relaxation and selection loss optimization, carries over directly.

\section{Theoretical analysis}\label{sec:theoretical-analysis}

A local error identity and bound are established for Gappy PMD. The stability
of QDEIM sampling is analyzed in relation to reconstruction convergence. The
classical reconstruction error estimate for Gappy POD is recovered as a linear
special case.

The Gappy PMD reconstruction of the full state associated with $\theta$ in
Eq.~\eqref{eq:pmd-representation} is denoted by $\mathcal F(\theta;y)$.
At any $\widetilde\theta\in\Theta$ at which $\mathcal F(\cdot;y)$ is
differentiable, the corresponding reconstruction and error are
\[
  \widetilde x:=\mathcal F(\widetilde\theta;y),\qquad
  \widetilde e:=x-\widetilde x.
\]
Let $\widetilde\Psi\in\mathbb{R}^{n\times\widetilde p}$ be any orthonormal
basis of
$\operatorname{range}\bigl(D_\theta\mathcal F(\widetilde\theta;y)\bigr)$,
where $\widetilde p$ is the dimension of this space.
The sampled basis and the normal component of the reconstruction error are
\[
  \widetilde\Psi_s:=S\widetilde\Psi,\qquad
  \widetilde e_\perp
  :=(I-\widetilde\Psi\widetilde\Psi^{\!\top})\widetilde e.
\]

\begin{theorem}[Local error identity and bound]
\label{thm:local-error}
Assume $1\leq\widetilde p\leq q$ and
$\operatorname{rank}(\widetilde\Psi_s)=\widetilde p$. With
\[
  \widetilde\Lambda
  :=\lVert\widetilde\Psi_s^\dagger\rVert_2
  =\frac{1}{\sigma_{\min}(\widetilde\Psi_s)},\qquad
  \widetilde r_{\mathrm{stat}}
  :=\widetilde\Psi^{\!\top}S^{\!\top}(y-S\widetilde x),
\]
the reconstruction error satisfies
\begin{align}
  \widetilde e
  &=\widetilde e_\perp
   -\widetilde\Psi\,\widetilde\Psi_s^\dagger
    (S\widetilde e_\perp+\eta)
   +\widetilde\Psi\,
    (\widetilde\Psi_s^{\!\top}\widetilde\Psi_s)^{-1}
    \widetilde r_{\mathrm{stat}},
  \label{eq:error-identity}\\[4pt]
  \lVert\widetilde e\rVert_2
  &\leq
  \widetilde\Lambda\bigl(
    \lVert\widetilde e_\perp\rVert_2
    +\lVert\eta\rVert_2
  \bigr)
  +\widetilde\Lambda^2\lVert\widetilde r_{\mathrm{stat}}\rVert_2.
  \label{eq:error-bound}
\end{align}
At an interior stationary point,
$\widetilde r_{\mathrm{stat}}=0$, so the corresponding terms vanish in
\eqref{eq:error-identity} and \eqref{eq:error-bound}. In particular,
$\lVert\widetilde e\rVert_2
\leq\widetilde\Lambda(\lVert\widetilde e_\perp\rVert_2
+\lVert\eta\rVert_2)$.
\end{theorem}

\begin{proof}
The orthogonal decomposition of $\widetilde e$ with respect to the local
tangent space is
\[
  \widetilde e
  =\widetilde\Psi\widetilde c+\widetilde e_\perp,
  \qquad
  \widetilde c:=\widetilde\Psi^{\!\top}\widetilde e.
\]
Using $y-S\widetilde x=S\widetilde e+\eta$ in the definition of
$\widetilde r_{\mathrm{stat}}$ gives
\[
  \widetilde r_{\mathrm{stat}}
  =\widetilde\Psi_s^{\!\top}\bigl(
    \widetilde\Psi_s\widetilde c
    +S\widetilde e_\perp+\eta\bigr).
\]
Since $\widetilde\Psi_s$ has full column rank,
$\widetilde\Psi_s^{\!\top}\widetilde\Psi_s$ is invertible and
\begin{equation}
  \widetilde c
  =-\widetilde\Psi_s^\dagger(S\widetilde e_\perp+\eta)
  +(\widetilde\Psi_s^{\!\top}\widetilde\Psi_s)^{-1}
   \widetilde r_{\mathrm{stat}}.
  \label{eq:coeff}
\end{equation}
Substituting~\eqref{eq:coeff} into this decomposition
yields~\eqref{eq:error-identity}.

The bound follows from the properties of the operator
$\widetilde{\mathcal P}
:=\widetilde\Psi\,\widetilde\Psi_s^\dagger S$.
Since
$\widetilde\Psi_s^\dagger\widetilde\Psi_s=I_{\widetilde p}$,
$\widetilde{\mathcal P}^2=\widetilde{\mathcal P}$,
so $\widetilde{\mathcal P}$ is idempotent.
Moreover, $\widetilde{\mathcal P}\neq0$ because $\widetilde p\geq1$.
Rewriting~\eqref{eq:error-identity} gives
\[
  \widetilde e
  =(I-\widetilde{\mathcal P})\widetilde e_\perp
  -\widetilde\Psi\,\widetilde\Psi_s^\dagger\eta
  +\widetilde\Psi\,
   (\widetilde\Psi_s^{\!\top}\widetilde\Psi_s)^{-1}
   \widetilde r_{\mathrm{stat}}.
\]
Since $\widetilde\Psi$ has orthonormal columns and
$SS^{\!\top}=I_q$,
\[
  \lVert\widetilde{\mathcal P}\rVert_2
  =\lVert\widetilde\Psi_s^\dagger S\rVert_2
  =\lVert\widetilde\Psi_s^\dagger\rVert_2
  =\widetilde\Lambda,
\]
where the second equality uses
$\lVert\widetilde\Psi_s^\dagger Sv\rVert_2
=\lVert\widetilde\Psi_s^\dagger w\rVert_2$
for $v=S^{\!\top}w$ with $\lVert v\rVert_2=\lVert w\rVert_2$
(since $SS^{\!\top}=I_q$).
For any nonzero idempotent $P$,
$\lVert I-P\rVert_2\leq\lVert P\rVert_2$,
so $\lVert I-\widetilde{\mathcal P}\rVert_2\leq\widetilde\Lambda$.
Since
$\lVert(\widetilde\Psi_s^{\!\top}\widetilde\Psi_s)^{-1}\rVert_2
=\widetilde\Lambda^2$,
the triangle inequality proves~\eqref{eq:error-bound}.
At an interior stationary point of Eq.~\eqref{eq:gappy-nls}, the first-order
condition is
\[
  0=-D_\theta\mathcal F(\widetilde\theta;y)^{\!\top}
    S^{\!\top}(y-S\widetilde x).
\]
Since $\widetilde\Psi$ spans
$\operatorname{range}(D_\theta\mathcal F(\widetilde\theta;y))$,
this condition implies $\widetilde r_{\mathrm{stat}}=0$ and proves the
stationary-point bound.
The factor $\widetilde\Lambda$ is independent of the chosen orthonormal
basis of the local tangent space.
\end{proof}

For QDEIM sampling, the following result gives an explicit bound for
the factor $\widetilde\Lambda$ in~\eqref{eq:error-bound}.  The bound
depends on the QDEIM stability constant and the misalignment of the
local tangent space with the POD subspace.
Here $S_q$ denotes the QDEIM sampling matrix associated with an
orthonormal POD basis $U_q\in\mathbb{R}^{n\times q}$, and
\[
  \Pi_q:=U_qU_q^{\!\top},\qquad
  B_q:=S_qU_q,\qquad
  \kappa_q:=\lVert B_q^{-1}\rVert_2
  =\frac{1}{\sigma_{\min}(B_q)}.
\]
The QDEIM construction ensures that $B_q$ is nonsingular and
$\kappa_q$ is finite.
For the orthonormal basis $\widetilde\Psi$ introduced above, where
$1\leq\widetilde p\leq q$, this misalignment is measured by
\[
  \tau_q:=\lVert(I-\Pi_q)\widetilde\Psi\rVert_2\in[0,1].
\]

\begin{proposition}
\label{prop:stability}
Let the sampling operator $S$ be the QDEIM sampling matrix $S_q$
associated with $U_q$, and assume $\kappa_q\tau_q<1$. Then
\begin{equation}
  \widetilde\Lambda
  \leq
  \Gamma_q(\tau_q)
  :=\frac{\kappa_q}{
    \sqrt{1-\tau_q^2}
    -\sqrt{\kappa_q^2-1}\,\tau_q}\,.
  \label{eq:lambda-bound}
\end{equation}
When $\tau_q=0$, this reduces to
$\widetilde\Lambda\leq\kappa_q$.
\end{proposition}

\begin{proof}
Decompose $\widetilde\Psi$ into its projection onto
$\operatorname{range}(U_q)$ and the corresponding orthogonal residual:
\[
  \widetilde\Psi
  =\Pi_q\widetilde\Psi+(I-\Pi_q)\widetilde\Psi
  =U_qM+E.
\]
Here, $M:=U_q^{\!\top}\widetilde\Psi$ contains the coordinates of the
projected component in the basis $U_q$, whereas
$E:=(I-\Pi_q)\widetilde\Psi$ is the orthogonal residual, with
$\lVert E\rVert_2=\tau_q$. Since $U_qM$ and $E$ are
orthogonal and $\widetilde\Psi$ has orthonormal columns,
\[
  I_{\widetilde p}
  =\widetilde\Psi^{\!\top}\widetilde\Psi
  =M^{\!\top}M+E^{\!\top}E.
\]
Consequently,
$M^{\!\top}M=I_{\widetilde p}-E^{\!\top}E$, which gives
\begin{equation}
  \sigma_{\min}(M)=\sqrt{1-\tau_q^2}.
  \label{eq:M-sigma}
\end{equation}
Since $S_qS_q^{\!\top}=I_q$ and $U_q^{\!\top}U_q=I_q$,
we have $\|S_q\|_2 = \|U_q\|_2 = 1$. Consequently,
\[
  \sigma_i(B_q)
  \leq\lVert B_q\rVert_2
  \leq\lVert S_q\rVert_2\lVert U_q\rVert_2
  =1
\]
for every $i$, and hence $\kappa_q\geq1$.
The identity $S_qS_q^{\!\top}=I_q$ also gives
\[
  \bigl[B_q^{-1}S_q(I-\Pi_q)\bigr]
  \bigl[B_q^{-1}S_q(I-\Pi_q)\bigr]^{\!\top}
  =B_q^{-1}(I_q-B_qB_q^{\!\top})B_q^{-\!\top}
  =B_q^{-1}B_q^{-\!\top}-I_q.
\]
Because $\sigma_i(B_q)\leq1$, all eigenvalues of
$B_q^{-1}B_q^{-\!\top}-I_q$ are nonnegative. This matrix is therefore
positive semidefinite. Its largest eigenvalue is $\kappa_q^2-1$, so
\begin{equation}
  \lVert B_q^{-1}S_q(I-\Pi_q)\rVert_2=\sqrt{\kappa_q^2-1}.
  \label{eq:offspace-norm}
\end{equation}
Now $B_q^{-1}S_q\widetilde\Psi=M+B_q^{-1}S_qE$.
Using the variational characterization of the smallest singular value,
we obtain
\[
  \sigma_{\min}(B_q^{-1}S_q\widetilde\Psi)
  \geq\sigma_{\min}(M)-\lVert B_q^{-1}S_qE\rVert_2
  \geq\sqrt{1-\tau_q^2}
      -\sqrt{\kappa_q^2-1}\,\tau_q.
\]
The condition $\kappa_q\tau_q<1$ ensures positivity.
Since $S_q\widetilde\Psi=B_q(B_q^{-1}S_q\widetilde\Psi)$,
\[
  \sigma_{\min}(S_q\widetilde\Psi)
  \geq\sigma_{\min}(B_q)\,
       \sigma_{\min}(B_q^{-1}S_q\widetilde\Psi)
  \geq
  \frac{
    \sqrt{1-\tau_q^2}
    -\sqrt{\kappa_q^2-1}\,\tau_q
  }{\kappa_q}\,.
\]
Taking reciprocals yields~\eqref{eq:lambda-bound}.
\end{proof}

\begin{corollary}
\label{cor:convergence}
Fix~$y$ and suppose $S=S_q$.
Assume $\Theta$ is compact,
$G(\theta;y):=D_\theta\mathcal F(\theta;y)$
is continuous with constant rank
$1\leq p:=\operatorname{rank}G(\theta;y)\leq q$ on~$\Theta$,
and $\kappa_q\tau_q(\theta;y)<1$ for every $\theta\in\Theta$,
where
$\tau_q(\theta;y)
:=\lVert(I-\Pi_q)G(\theta;y)G(\theta;y)^\dagger\rVert_2$.
The maximum misalignment over~$\Theta$ is
$\tau_{q,*}:=\max_{\theta\in\Theta}\tau_q(\theta;y)$, and
$\Lambda_{\mathcal T}(\theta;y)$ denotes the corresponding tangent
space factor.
Then
\begin{equation}
  \sup_{\theta\in\Theta}\Lambda_{\mathcal T}(\theta;y)
  \leq\Gamma_q(\tau_{q,*})<\infty,
  \label{eq:uniform-bound}
\end{equation}
and \Cref{thm:local-error} gives
\begin{equation}
  \lVert\widetilde e\rVert_2
  \leq\Gamma_q(\tau_{q,*})\bigl(
  \lVert\widetilde e_\perp\rVert_2
  +\lVert\eta\rVert_2\bigr)
  +\Gamma_q(\tau_{q,*})^2\,
  \lVert\widetilde r_{\mathrm{stat}}\rVert_2
  \label{eq:global-bound}
\end{equation}
for all $\widetilde\theta\in\Theta$.
More generally, consider a sequence of problems indexed by~$j$,
each satisfying the above assumptions with stability constant
$\Gamma_q(\tau_{q,*,j})$.
If\; $\sup_j\Gamma_q(\tau_{q,*,j})<\infty$,\;
$\lVert\widetilde e_{\perp,j}\rVert_2\to0$,\;
$\lVert\eta_j\rVert_2\to0$,\;
and\;
$\lVert\widetilde r_{\mathrm{stat},j}\rVert_2\to0$,\;
then $\lVert\widetilde e_j\rVert_2\to0$.
Here $\widetilde e_{\perp,j}$ is the local normal component at the
recovered point. It is not, by itself, the global PMD approximation
error.
\end{corollary}

\begin{proof}
Under the constant-rank assumption, the Moore--Penrose inverse
$G(\theta;y)^\dagger$ depends continuously on~$\theta$, so
$\tau_q(\theta;y)$ is continuous on the compact set~$\Theta$.
The maximum~$\tau_{q,*}$ is therefore attained at some
$\theta_*\in\Theta$, and evaluating the pointwise condition there
gives $\kappa_q\tau_{q,*}<1$.
For any orthonormal basis $\Psi(\theta;y)$ of the tangent space,
$\lVert(I-\Pi_q)\Psi(\theta;y)\rVert_2
=\lVert(I-\Pi_q)G(\theta;y)G(\theta;y)^\dagger\rVert_2
=\tau_q(\theta;y)$.
\Cref{prop:stability} and the monotonicity of $\Gamma_q$
on $[0,1/\kappa_q)$ yield~\eqref{eq:uniform-bound}.
The global bound and convergence follow by substituting
into~\eqref{eq:error-bound}.
\end{proof}

The same identities recover the classical Gappy POD estimate when the
nonlinear lift is omitted.

\begin{proposition}
\label{prop:gappy-pod}
Remove the nonlinear lift and consider
$\mathcal F_{\mathrm{POD}}(a)=\bar x+\Phi a$, where
$\Phi\in\mathbb{R}^{n\times r}$ has orthonormal columns
($\Phi^{\!\top}\Phi=I_r$) and $\operatorname{rank}(S\Phi)=r$.
Let $\widehat a$ satisfy the normal equations
$\Phi^{\!\top}S^{\!\top}(y-S(\bar x+\Phi\widehat a))=0$,
and set $\widehat x:=\bar x+\Phi\widehat a$,
$u_\Phi:=(I-\Phi\Phi^{\!\top})(x-\bar x)$.  Then
\begin{align}
  x-\widehat x
  &=u_\Phi-\Phi(S\Phi)^\dagger(Su_\Phi+\eta),
  \label{eq:pod-identity}\\[4pt]
  \lVert x-\widehat x\rVert_2^2
  &=\lVert u_\Phi\rVert_2^2
   +\lVert(S\Phi)^\dagger(Su_\Phi+\eta)\rVert_2^2,
  \label{eq:pod-pythagoras}\\[4pt]
  \lVert x-\widehat x\rVert_2
  &\leq
  \frac{\lVert u_\Phi\rVert_2+\lVert\eta\rVert_2}
       {\sigma_{\min}(S\Phi)}.
  \label{eq:pod-bound}
\end{align}
Since $D_a^2\mathcal F_{\mathrm{POD}}\equiv0$, the curvature
vanishes.
\end{proposition}

\begin{proof}
The tangent space is $\operatorname{range}(\Phi)$ with orthonormal basis
$\Psi=\Phi$.  Since
$\widehat x-\bar x\in\operatorname{range}(\Phi)$,
\[
  \widetilde e_\perp
  =(I-\Phi\Phi^{\!\top})(x-\widehat x)
  =(I-\Phi\Phi^{\!\top})(x-\bar x)
  =u_\Phi.
\]
The normal equations give $\widetilde r_{\mathrm{stat}}=0$.
Substituting $\Psi=\Phi$, $\widetilde e_\perp=u_\Phi$, and
$\widetilde r_{\mathrm{stat}}=0$
into~\eqref{eq:error-identity}
proves~\eqref{eq:pod-identity}.
The two terms on the right are orthogonal
(the first lies in $\operatorname{range}(\Phi)^\perp$, the second in
$\operatorname{range}(\Phi)$),
which proves~\eqref{eq:pod-pythagoras}.
The bound~\eqref{eq:pod-bound} follows from the stationary
specialization of~\eqref{eq:error-bound} with
$\widetilde\Lambda=\sigma_{\min}(S\Phi)^{-1}$.
\end{proof}

\section{Numerical results}\label{sec:numerical-results}
The proposed Gappy PMD method and the DPS point selection strategy are evaluated on three benchmark incompressible viscous flow problems: flow past a cylinder, the lid-driven cavity flow, and the backward-facing step flow. All three test cases are governed by the incompressible Navier--Stokes equations:
\begin{align}
\frac{\partial \mathbf{u}}{\partial t} + (\mathbf{u} \cdot \nabla)\mathbf{u}
&= -\frac{1}{\rho}\nabla p + \nu \nabla^2 \mathbf{u}, \\
\nabla \cdot \mathbf{u} &= 0,
\end{align}
where $\mathbf{u} = \mathbf{u}(\mathbf{x}, t)$ denotes the velocity field, $p = p(\mathbf{x}, t)$ the pressure field, $\rho$ the fluid density, and $\nu$ the kinematic viscosity.

The numerical simulation data for all three test cases are generated using the open-source finite-element computational fluid dynamics solver Fluidity. For each problem, a total of 3,200 solution snapshots are collected over the simulation time horizon. The first 2,000 snapshots constitute the offline training set ($m=2000$), used to construct the basis functions for the reduced-order model. The subsequent 200 snapshots serve as a validation set ($|\mathcal{V}|=200$) for optimization tasks such as the DPS point selection and hyperparameter tuning. The remaining 1,000 snapshots form the online test set on which reconstruction accuracy is evaluated. The reconstruction quality is quantified by the relative $L^2$ error between the reconstructed field $\widehat{x}$ and the reference solution $x$, defined as
\begin{equation}
e_{\mathrm{rel}} = \frac{\lVert x - \widehat{x}\rVert_2}{\lVert x\rVert_2}.
\end{equation}

The reconstruction accuracy of Gappy PMD is compared against Gappy POD. To ensure a fair comparison, the POD basis dimension for Gappy POD
is set equal to the total PMD dimension $r+r_1$
for each test case. Both methods share the same set of sampling points determined by QDEIM from the POD basis, so that the selected points are naturally adapted to Gappy POD. To further demonstrate the advantage of DPS, the proposed point selection strategy is tested within both the Gappy PMD and Gappy POD frameworks, each compared against QDEIM. To assess robustness to measurement noise, each test case is additionally evaluated under corrupted measurements of the form
\begin{equation}
y = S x + \eta, \quad
\eta \sim \mathcal{N}(\mathbf{0}, \sigma^2 I_q), \quad
\sigma = \frac{\eta_\%}{100}\,\sigma_{\mathrm{train}}, \quad
\sigma_{\mathrm{train}} = \sqrt{\frac{1}{qm}\sum_{i=1}^{q}\sum_{j=1}^{m}(Y_{ij} - \bar{Y})^2},
\end{equation}
where $S x \in \mathbb{R}^q$ denotes the noiseless sampled state at the $q$ selected points, $S \in \mathbb{R}^{q \times n}$ is the sampling operator, $\eta_\%$ is the prescribed noise level, $Y = SX \in \mathbb{R}^{q \times m}$ is the sampled training matrix, and $\bar{Y}$ the mean of its entries. 
For each noise level, the reported error is averaged over five independent noise
realizations.

In addition to reconstruction accuracy, we compare the online computational time of Gappy POD and Gappy PMD. As reported in Table~\ref{tab:online_time}, Gappy PMD requires a higher online computational cost than Gappy POD, owing to the nonlinear optimization and the subsequent reconstruction from nonlinear coordinates. Nevertheless, the wall-clock time per snapshot remains on the order of milliseconds, which is acceptable in practice. All timings are measured on an Intel Core i5-12500H processor.

\begin{table}[t]
  \centering
  \caption{Problem setup and average online computational time per snapshot
    for the three test cases
    ($n$: number of spatial degrees of freedom;
     $q$: number of sampling points;
     $r$: linear subspace dimension;
     $r_1$: nonlinear manifold dimension).}
  \label{tab:online_time}
  \begin{tabular}{l cc ccc cc}
    \toprule
    & & &
    \multicolumn{3}{c}{Reduced dimension} &
    \multicolumn{2}{c}{Online wall-clock time (ms/snapshot)} \\
    \cmidrule(lr){4-6} \cmidrule(lr){7-8}
    & & &
    \multirow{2}{*}{Gappy POD} &
    \multicolumn{2}{c}{Gappy PMD} &
    \multirow{2}{*}{Gappy POD} &
    \multirow{2}{*}{Gappy PMD} \\
    \cmidrule(lr){5-6}
    Case & $n$ & $q$ & & $r$ & $r_1$ & & \\
    \midrule
    Flow past a cylinder       & 11{,}930 & 6 & 4 & 2 & 2 & 0.0642 & 1.7784 \\
    Lid-driven cavity flow     & 10{,}201 & 8 & 6 & 2 & 4 & 0.0744 & 5.1260 \\
    Backward-facing step flow  & 18{,}058 & 9 & 7 & 2 & 5 & 0.0847 & 5.5743 \\
    \bottomrule
  \end{tabular}
\end{table}

\subsection{Flow past a cylinder}

\begin{figure}[htbp]
   \centering
   \includegraphics[width=0.85\textwidth]{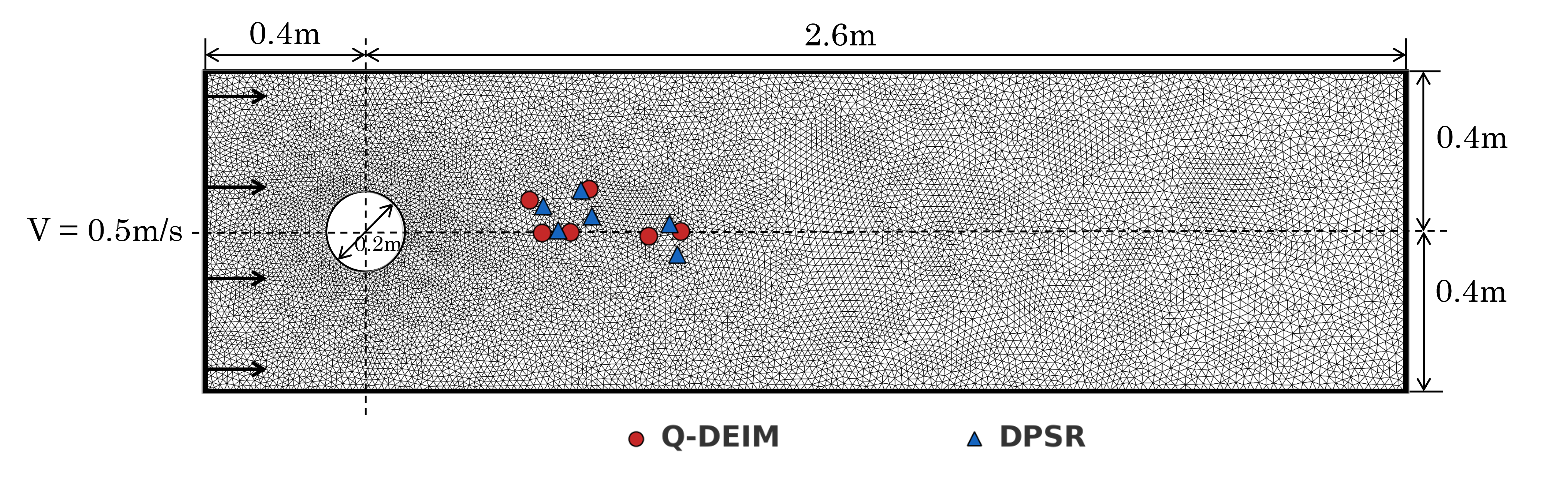}
   \caption{Computational domain, boundary conditions and unstructured triangular mesh for
   flow past a cylinder. Red circles mark the six sampling points selected by QDEIM and blue
   triangles those selected by DPS.}
   \label{fig:cylinder_domain}
\end{figure}

We first consider two-dimensional flow past a cylinder at 
Reynolds number $\mathrm{Re} = \rho U D / \mu = 100$, where 
$\rho = 1 \, \mathrm{kg/m^3}$ is the fluid density, 
$U = 0.5 \, \mathrm{m/s}$ is the inlet velocity, 
$D = 0.2 \, \mathrm{m}$ is the cylinder diameter, and 
$\mu = 0.001 \, \mathrm{Pa \cdot s}$ is the dynamic viscosity. 
At this Reynolds number, the flow develops periodic vortex shedding 
in the wake region.
The computational domain is discretized using an unstructured 
triangular mesh with 11,930 nodes, as shown in Figure~\ref{fig:cylinder_domain}. A uniform 
inflow velocity of $0.5\,\mathrm{m/s}$ is imposed on the left 
boundary, and no-slip conditions are applied on the upper and 
lower walls as well as the cylinder surface. The flow is simulated 
using Fluidity for a total duration of $32$\,s with a time step of 
$\Delta t = 0.01$\,s, producing 3,200 snapshots. The streamwise
velocity $u_x$ is used for reconstruction.
For Gappy PMD, the linear subspace dimension is set to 2, and the 
nonlinear manifold dimension is also set to 2. The number of sampling 
points is set to 6, which oversamples relative to the manifold 
dimension to improve the numerical stability of the reconstruction.
The sampling locations selected by QDEIM and DPS for this configuration are also indicated in Figure~\ref{fig:cylinder_domain}.

\begin{figure}[htbp]
   \centering
   \begin{tabular}{cc}
      \includegraphics[width=0.48\textwidth]{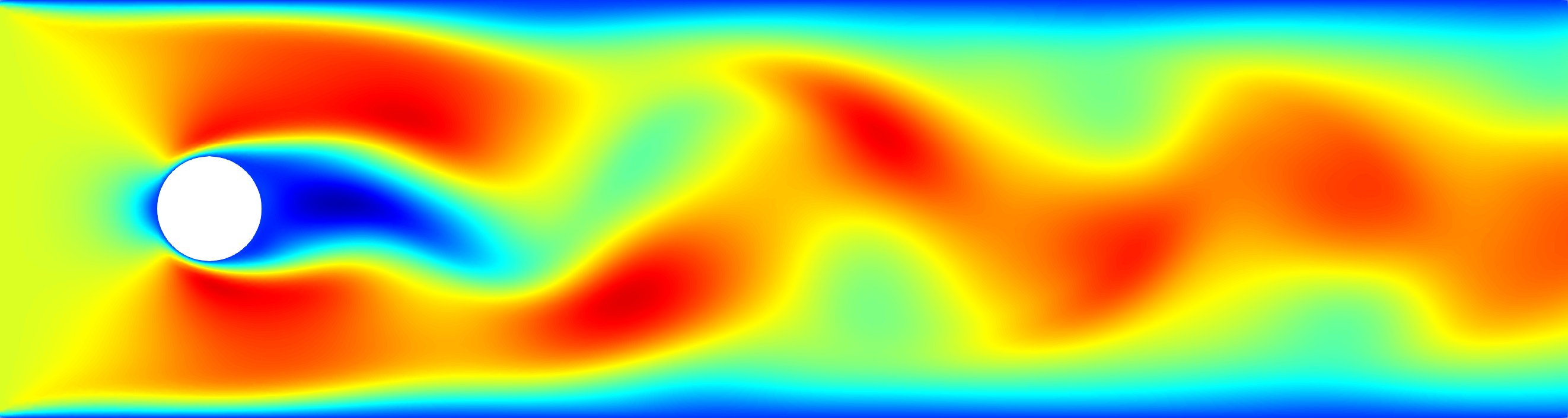}
      & \includegraphics[width=0.48\textwidth]{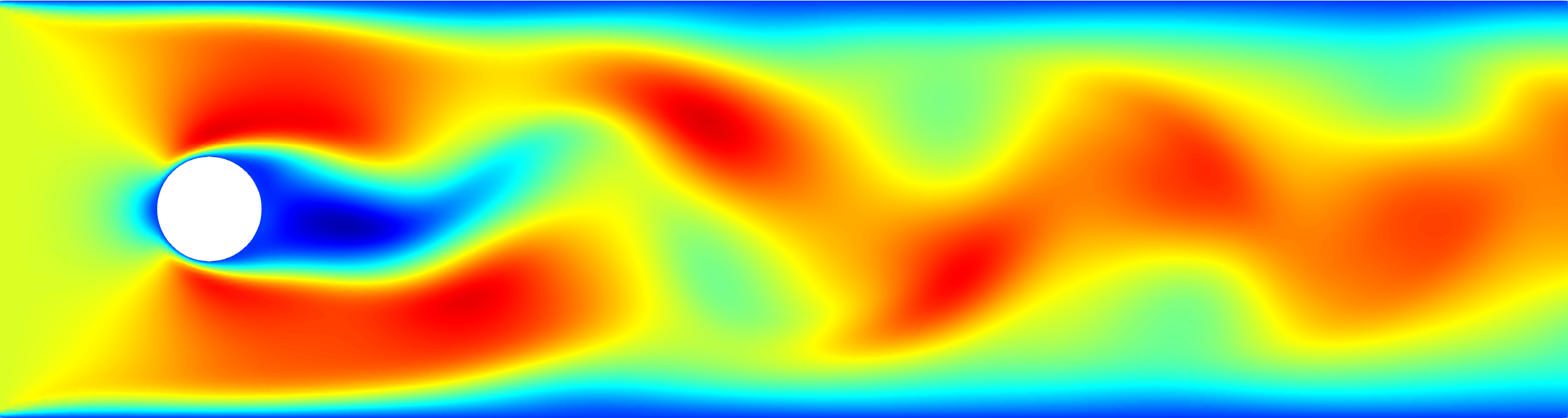}              \\
      (a) Reference solution, $t = 24\,\mathrm{s}$
      & (b) Reference solution, $t = 28\,\mathrm{s}$                           \\[6pt]
      \includegraphics[width=0.48\textwidth]{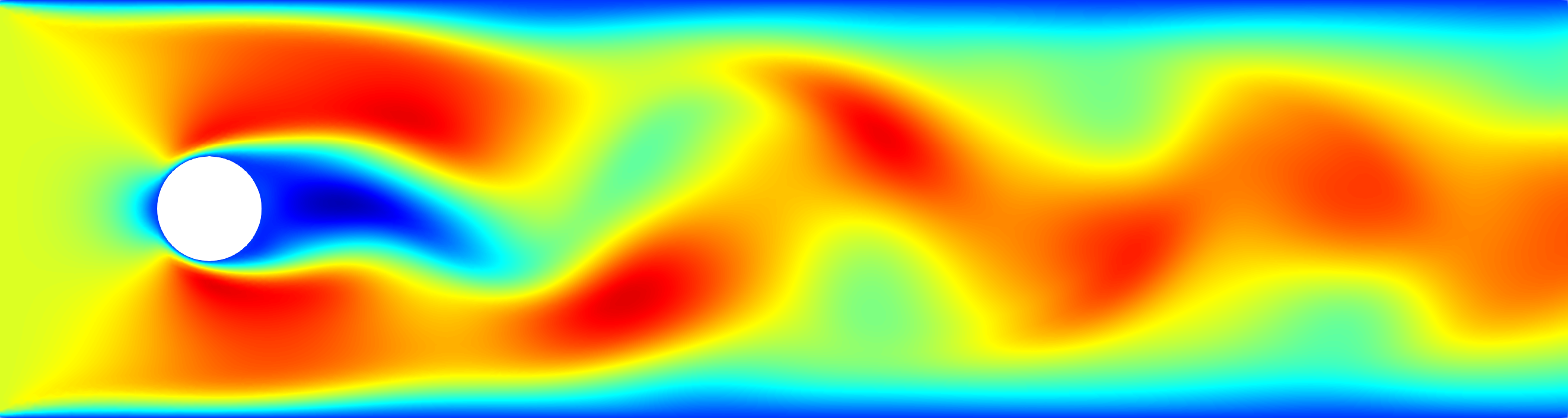}
      & \includegraphics[width=0.48\textwidth]{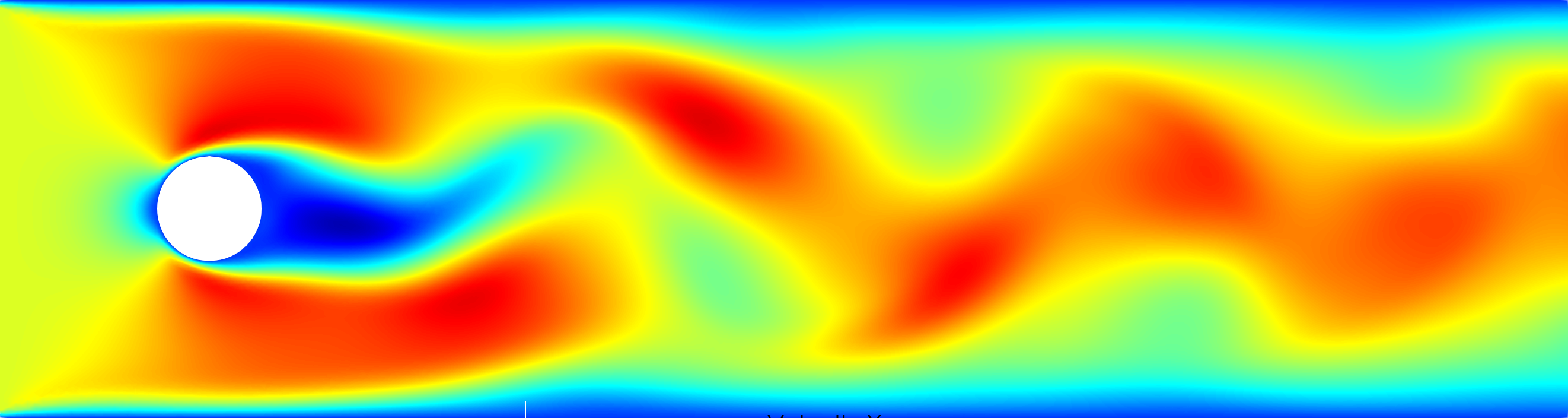}     \\
      (c) Gappy PMD + QDEIM, $t = 24\,\mathrm{s}$
      & (d) Gappy PMD + QDEIM, $t = 28\,\mathrm{s}$                            \\[6pt]
      \includegraphics[width=0.48\textwidth]{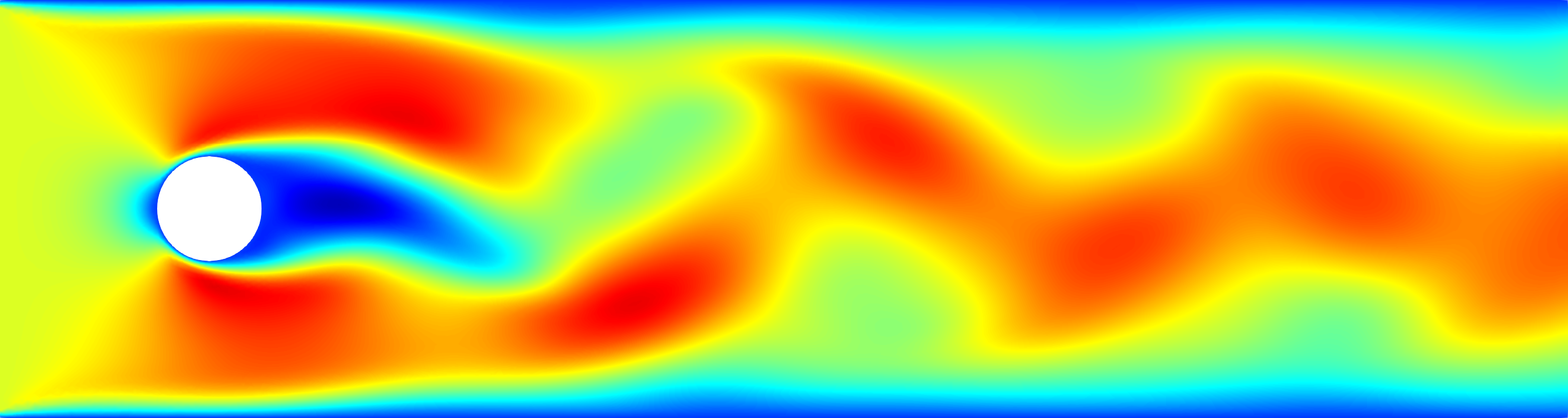}
      & \includegraphics[width=0.48\textwidth]{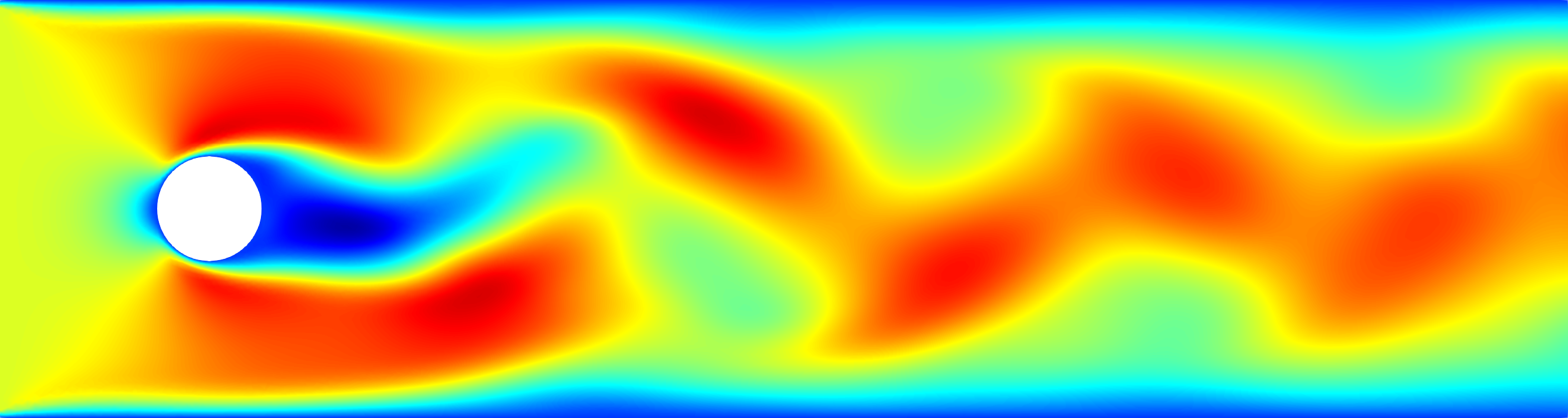}     \\
      (e) Gappy POD + QDEIM, $t = 24\,\mathrm{s}$
      & (f) Gappy POD + QDEIM, $t = 28\,\mathrm{s}$                            \\[2pt]
      \includegraphics[width=0.48\textwidth]{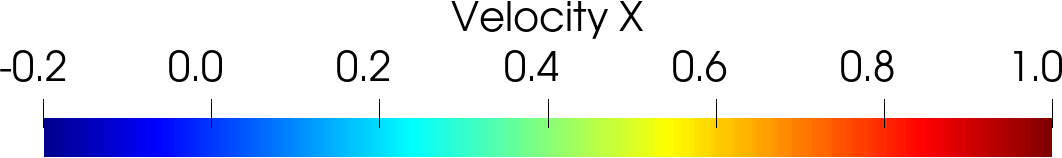}
      & \includegraphics[width=0.48\textwidth]{cylinder_colorbar1.png}           \\[6pt]
      \includegraphics[width=0.48\textwidth]{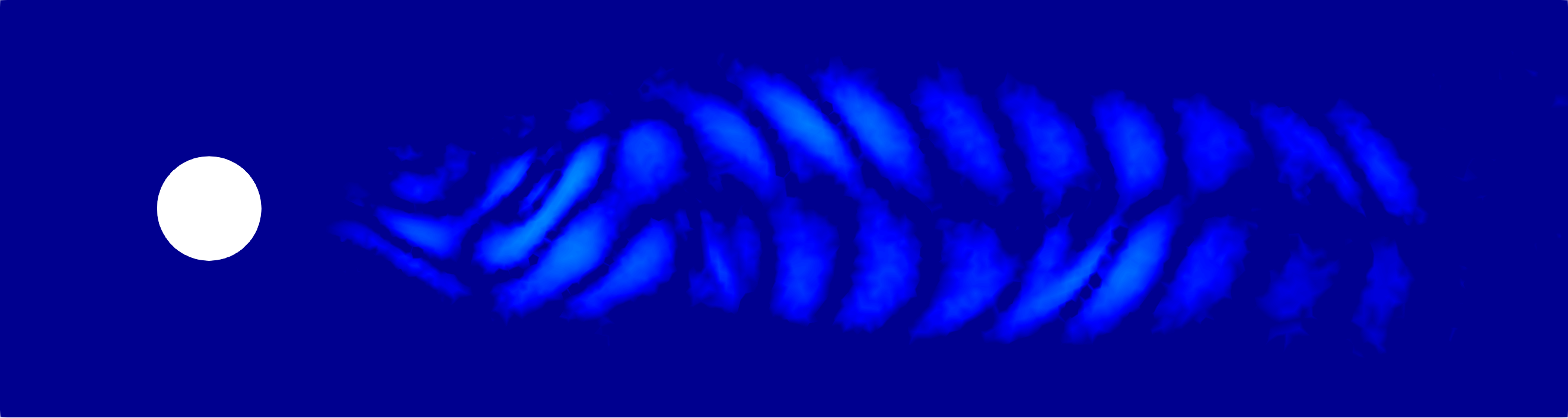}
      & \includegraphics[width=0.48\textwidth]{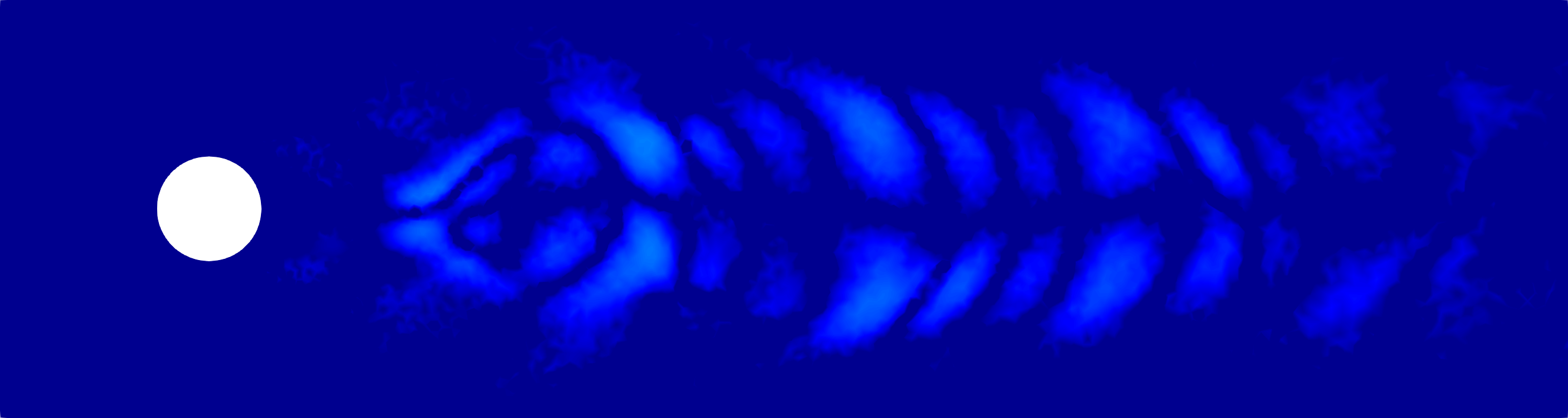} \\
      (g) Error of Gappy PMD + QDEIM, $t = 24\,\mathrm{s}$
      & (h) Error of Gappy PMD + QDEIM, $t = 28\,\mathrm{s}$                   \\[6pt]
      \includegraphics[width=0.48\textwidth]{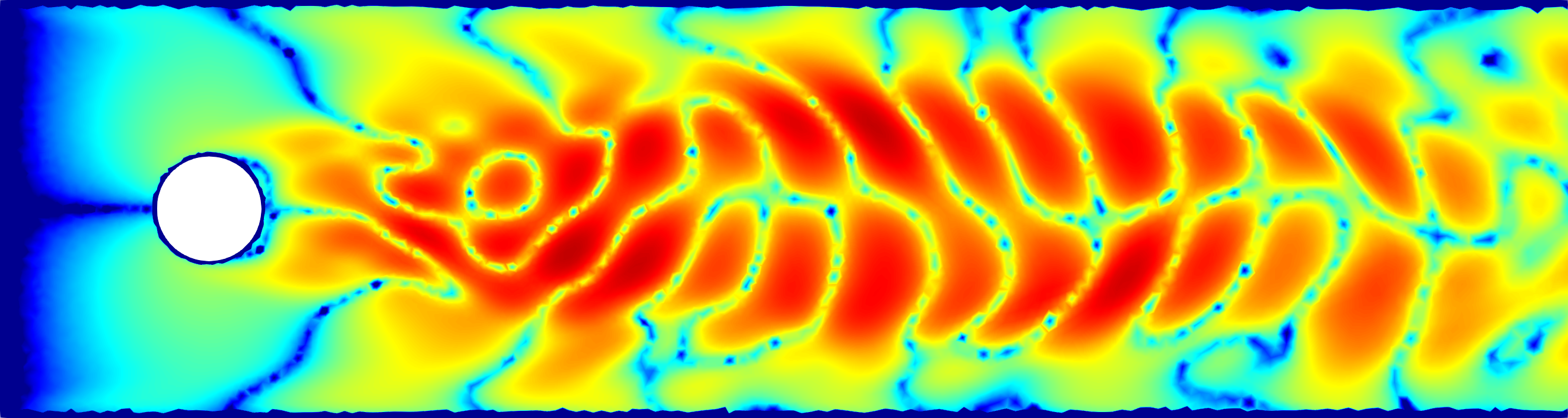}
      & \includegraphics[width=0.48\textwidth]{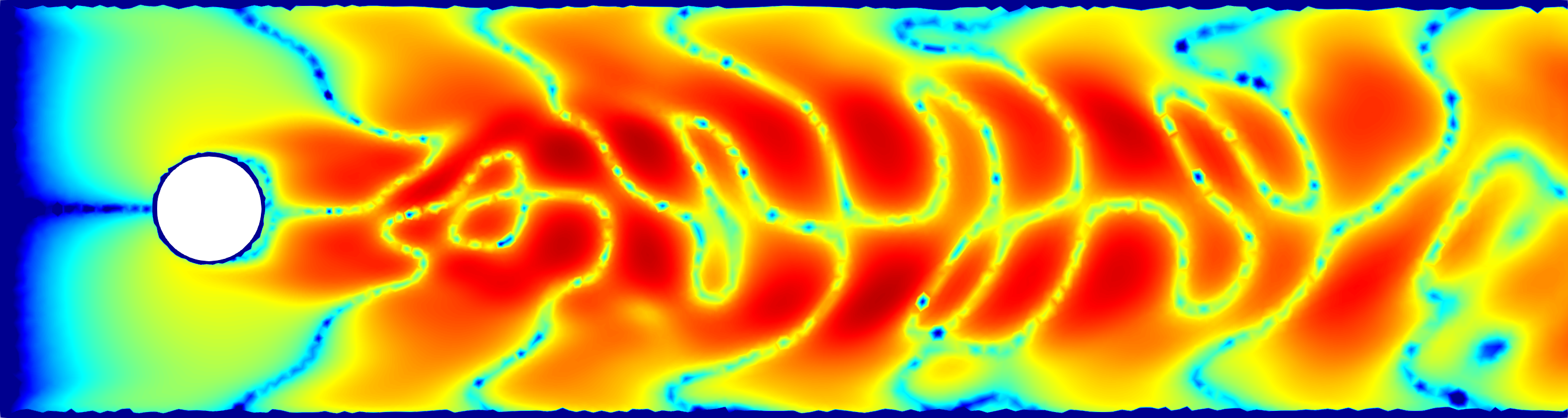} \\
      (i) Error of Gappy POD + QDEIM, $t = 24\,\mathrm{s}$
      & (j) Error of Gappy POD + QDEIM, $t = 28\,\mathrm{s}$                   \\[2pt]
      \includegraphics[width=0.48\textwidth]{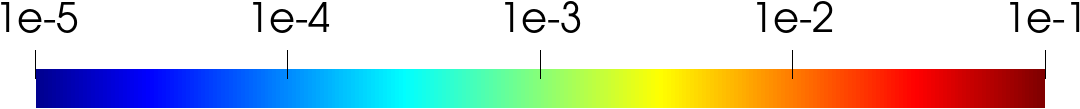}
      & \includegraphics[width=0.48\textwidth]{cylinder_err_colorbar1.png}       \\
   \end{tabular}
   \caption{Reconstructed streamwise velocity $u_x$ and the corresponding absolute error
   fields for flow past a cylinder, comparing Gappy PMD and Gappy POD.}
   \label{fig:cylinder_recon_comparison}
\end{figure}

Figure~\ref{fig:cylinder_recon_comparison} compares Gappy PMD and Gappy POD using the
same six sampling points, at two instants of the vortex shedding cycle.
Both methods reproduce the vortex shedding pattern of the reference solution.
The difference between them appears only in the corresponding error fields, which are
plotted on a logarithmic scale.
The error of Gappy PMD remains low over the whole domain and rises only where the vortex
cores pass. The Gappy POD error is larger by orders of magnitude and follows the vortex
street along the entire wake, with the highest values in the near wake and in the
separated shear layers.
These error fields therefore indicate that Gappy PMD captures the wake structures
considerably more accurately than Gappy POD.

\begin{figure}[htbp]
   \centering
   \begin{tabular}{cc}
      \includegraphics[width=0.48\textwidth]{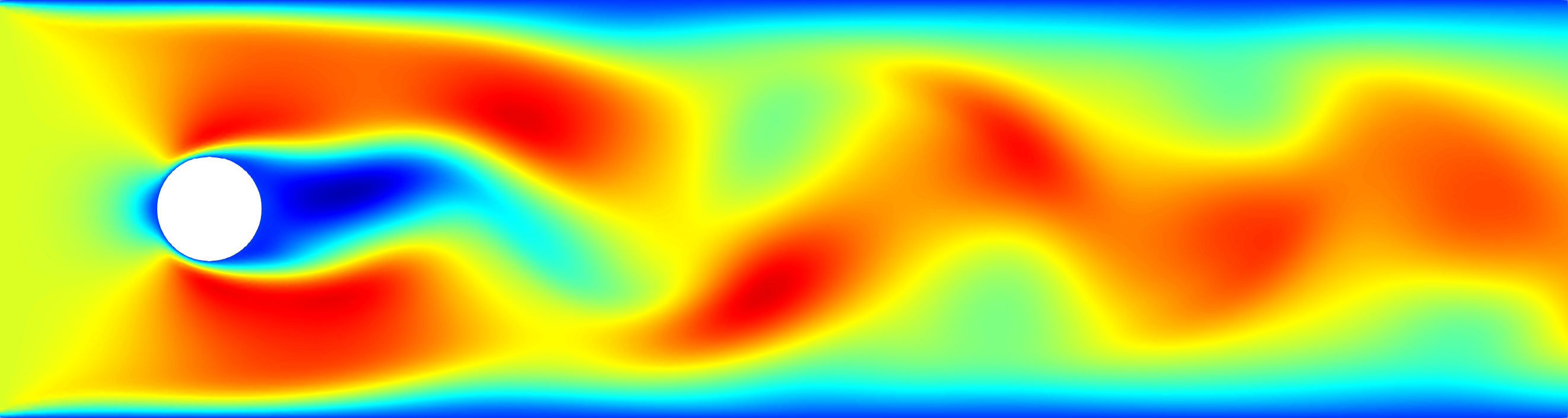}
      & \includegraphics[width=0.48\textwidth]{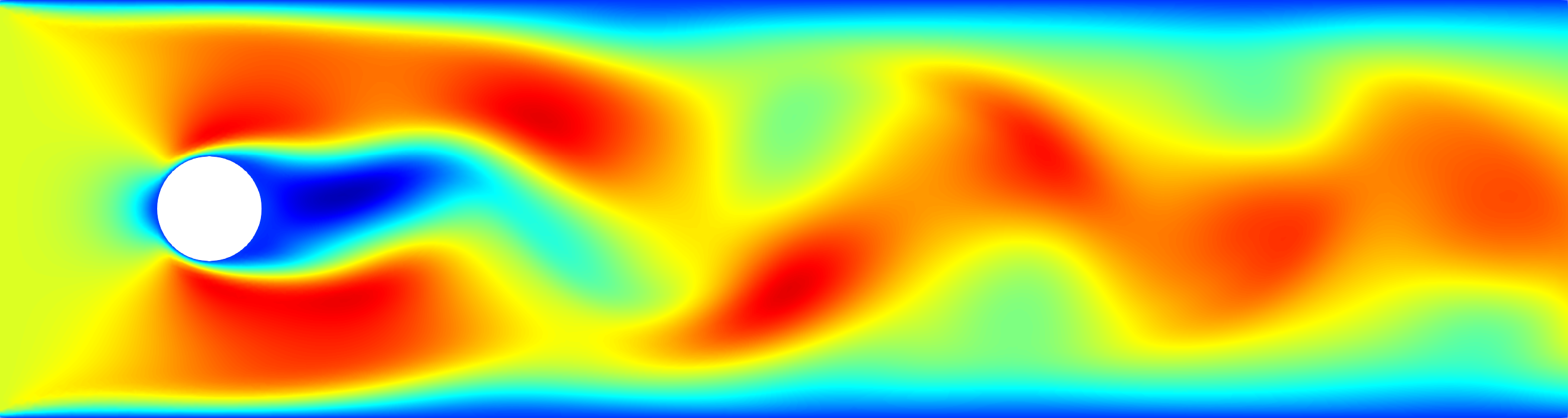}              \\
      (a) Reference solution, $t = 26\,\mathrm{s}$
      & (b) Reference solution, $t = 29\,\mathrm{s}$                           \\[6pt]
      \includegraphics[width=0.48\textwidth]{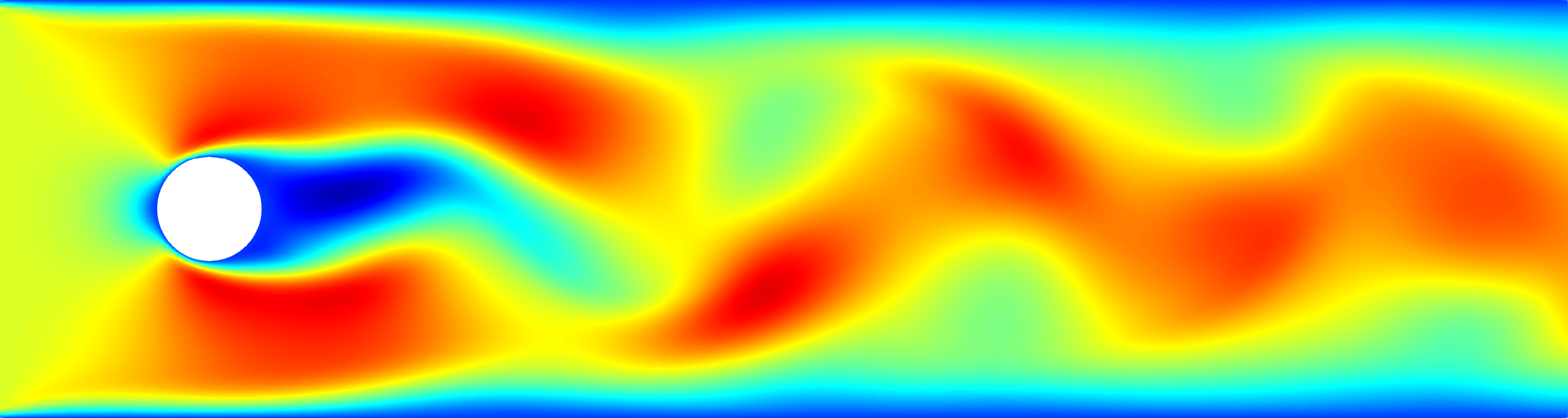}
      & \includegraphics[width=0.48\textwidth]{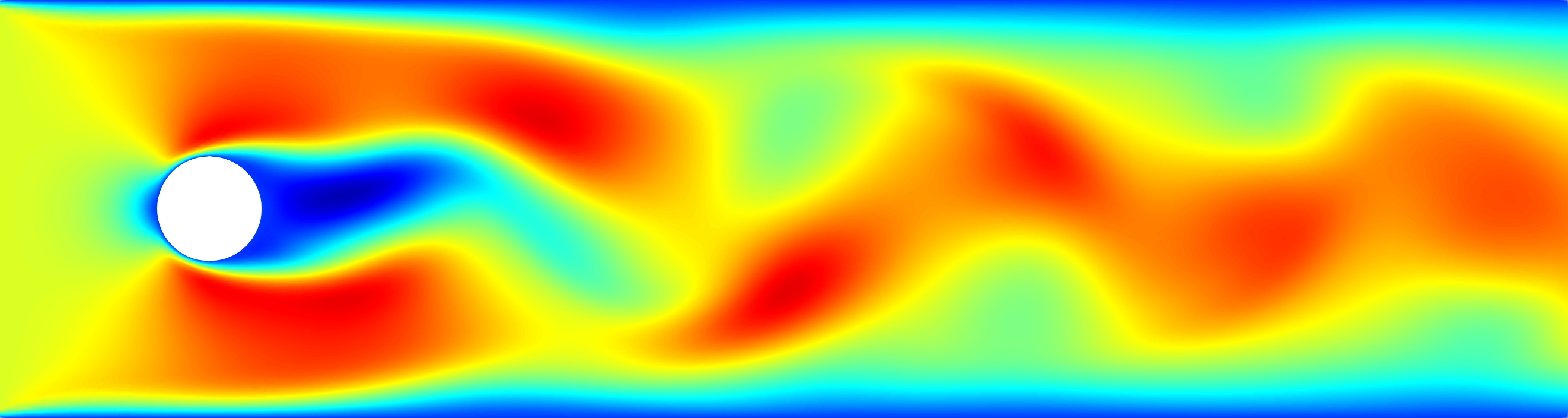}      \\
      (c) Gappy PMD + DPS, $t = 26\,\mathrm{s}$
      & (d) Gappy PMD + DPS, $t = 29\,\mathrm{s}$                             \\[6pt]
      \includegraphics[width=0.48\textwidth]{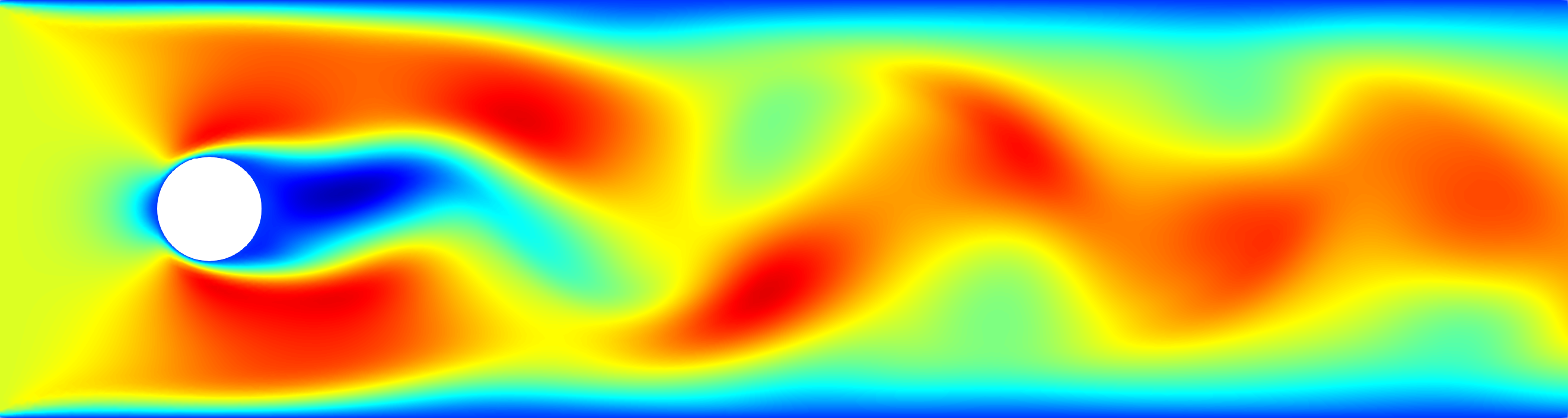}
      & \includegraphics[width=0.48\textwidth]{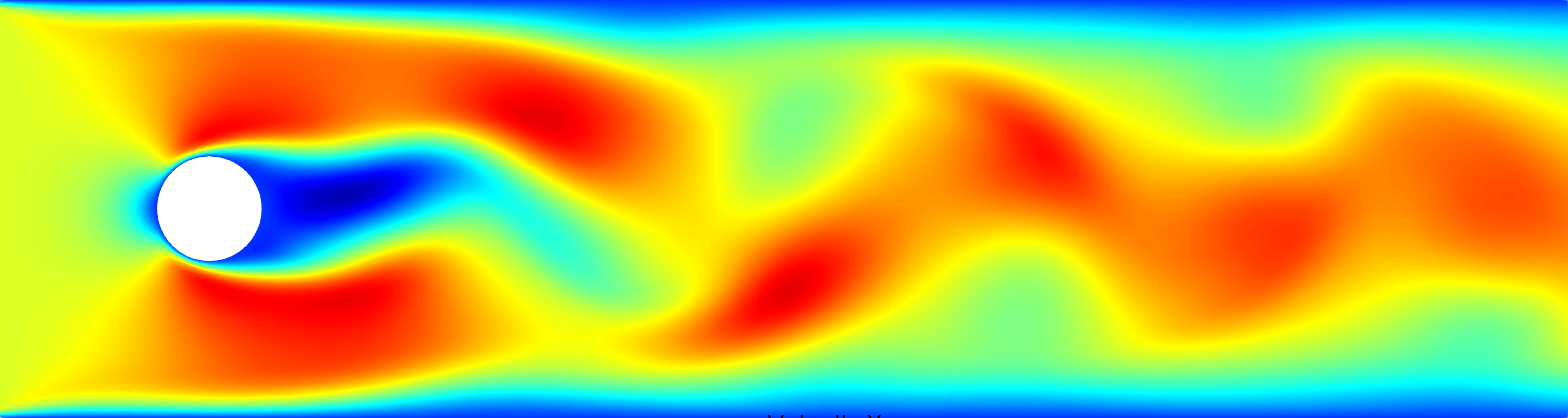}     \\
      (e) Gappy PMD + QDEIM, $t = 26\,\mathrm{s}$
      & (f) Gappy PMD + QDEIM, $t = 29\,\mathrm{s}$                            \\[2pt]
      \includegraphics[width=0.48\textwidth]{cylinder_colorbar1.png}
      & \includegraphics[width=0.48\textwidth]{cylinder_colorbar1.png}           \\[6pt]
      \includegraphics[width=0.48\textwidth]{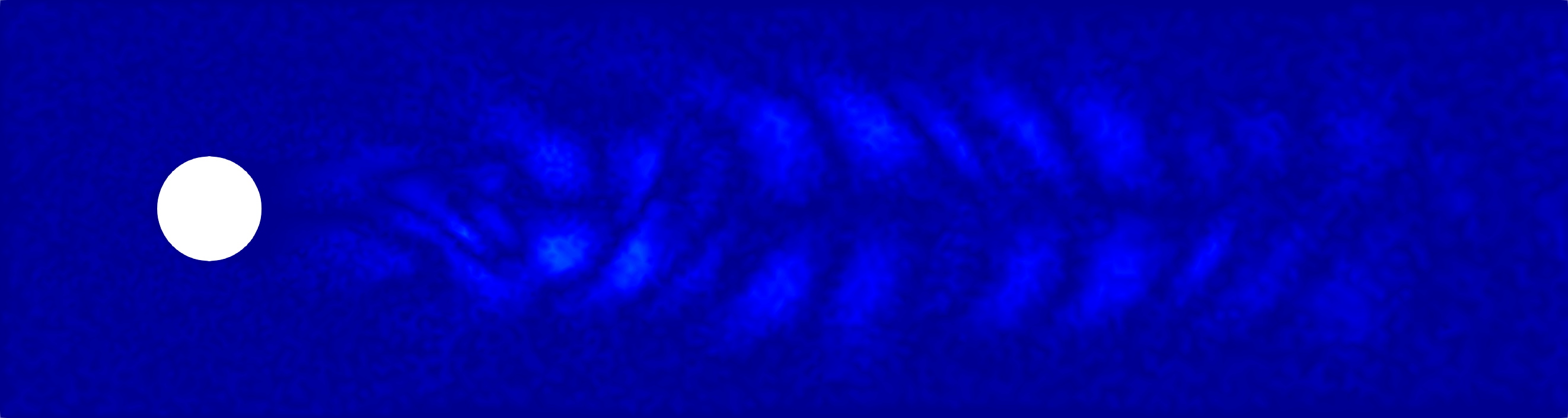}
      & \includegraphics[width=0.48\textwidth]{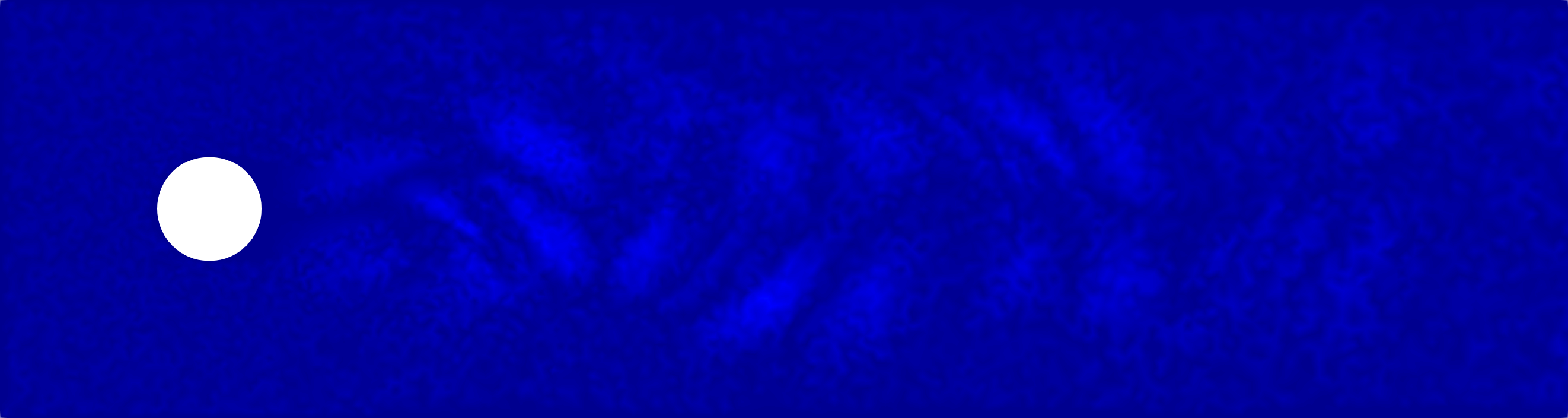}  \\
      (g) Error of Gappy PMD + DPS, $t = 26\,\mathrm{s}$
      & (h) Error of Gappy PMD + DPS, $t = 29\,\mathrm{s}$                    \\[6pt]
      \includegraphics[width=0.48\textwidth]{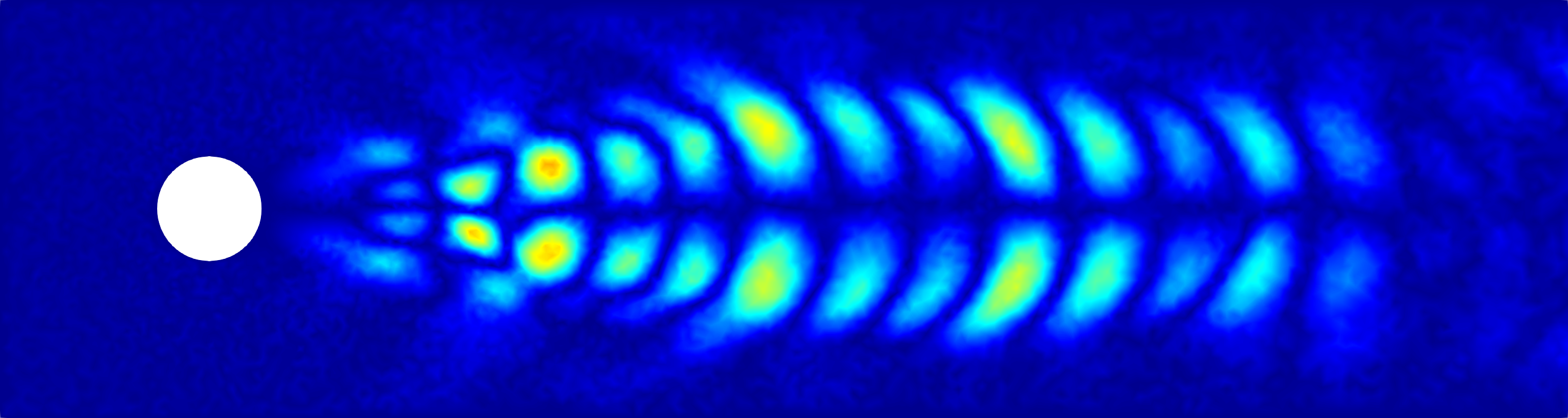}
      & \includegraphics[width=0.48\textwidth]{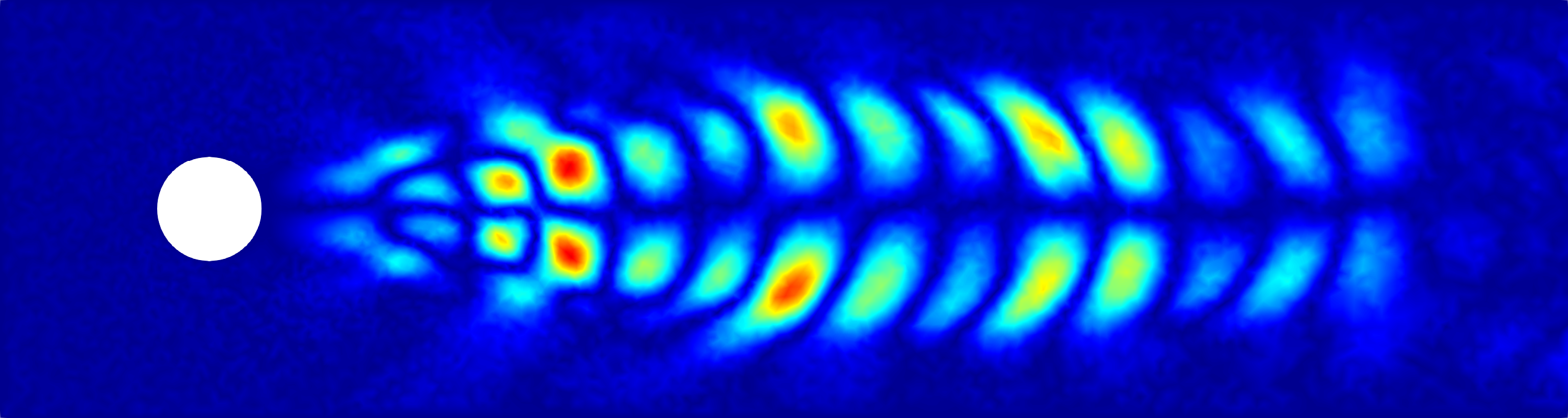} \\
      (i) Error of Gappy PMD + QDEIM, $t = 26\,\mathrm{s}$
      & (j) Error of Gappy PMD + QDEIM, $t = 29\,\mathrm{s}$                   \\[2pt]
      \includegraphics[width=0.48\textwidth]{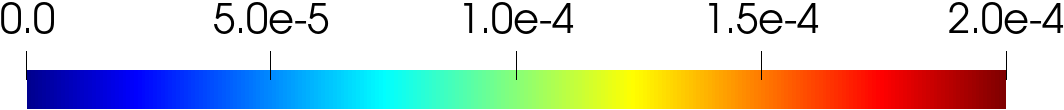}
      & \includegraphics[width=0.48\textwidth]{cylinder_err_colorbar2.png}       \\
   \end{tabular}
   \caption{Reconstructed streamwise velocity $u_x$ and the corresponding absolute error
   fields for flow past a cylinder, comparing QDEIM and DPS sampling under Gappy PMD.}
   \label{fig:cylinder_sampling_comparison}
\end{figure}

Figure~\ref{fig:cylinder_sampling_comparison} compares the QDEIM and DPS sampling
points under the same Gappy PMD reconstruction.
With QDEIM sampling the error is organized by the vortex street, with pronounced error
cores in the near wake and alternating bands extending downstream.
With DPS sampling the error is much weaker and the pattern far less pronounced.
The weaker and less structured error under DPS indicates that, for the same number of
sampling points, its locations are better suited to the Gappy PMD reconstruction.

\begin{figure}[htbp]
   \centering
   \includegraphics[width=0.6\textwidth]{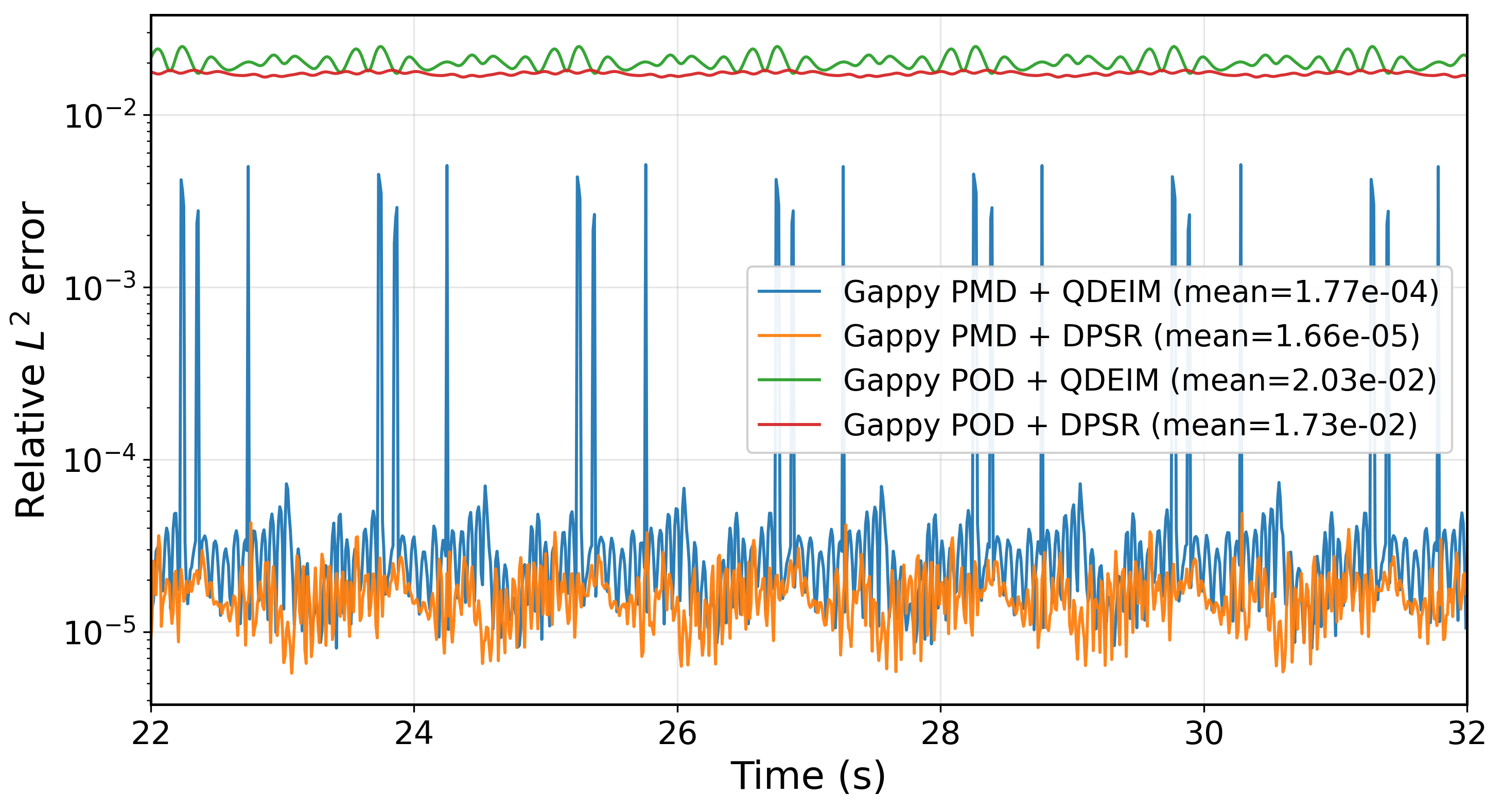}
   \caption{Relative $L^2$ reconstruction error over the test snapshots for flow past a
   cylinder, for the four combinations of reconstruction method and sampling strategy.}
   \label{fig:cylinder_error_curve}
\end{figure}

Figure~\ref{fig:cylinder_error_curve} shows the relative $L^2$ error over the test
snapshots for the four combinations of reconstruction method and sampling strategy.
Gappy PMD is more accurate than Gappy POD by orders of magnitude throughout the test
window, regardless of which sampling strategy is used.
For both reconstruction methods, DPS sampling gives a lower error than QDEIM sampling.
This advantage is most apparent in the shape of the curves. Under QDEIM sampling the
Gappy PMD error develops pronounced spikes that recur at the vortex shedding period,
whereas under DPS sampling no comparable excursions occur.
DPS sampling therefore yields a reconstruction that is robust over the entire shedding
cycle rather than accurate only at favorable phases.

\begin{figure}[htbp]
   \centering
   \includegraphics[width=0.65\textwidth]{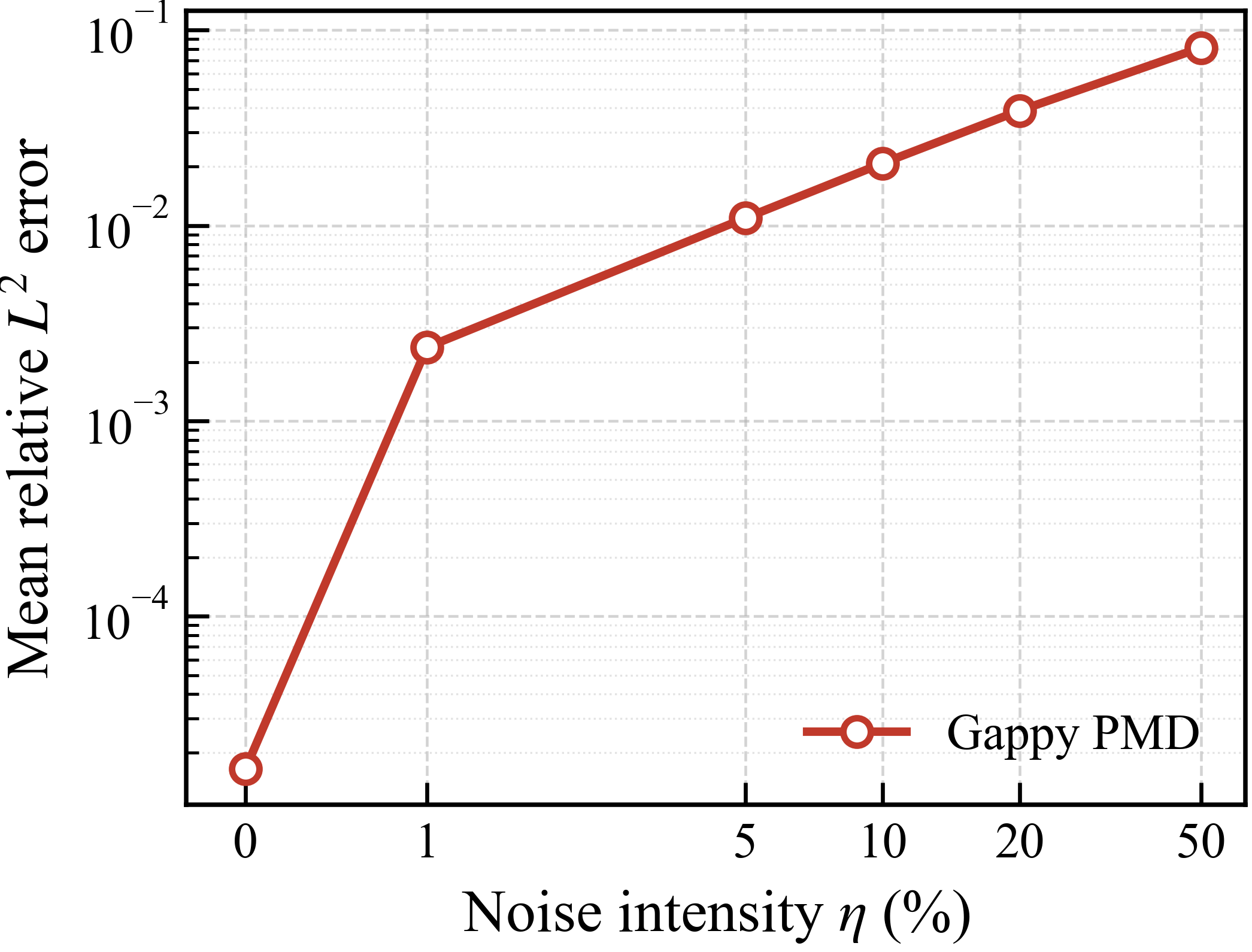}
   \caption{Effect of the noise level $\eta_\%$ on the mean relative $L^2$ error of Gappy PMD
   with DPS sampling, averaged over the test snapshots and over five independent noise
   realizations.}
   \label{fig:cylinder_noise_robustness}
\end{figure}

Figure~\ref{fig:cylinder_noise_robustness} shows the mean relative $L^2$ error of Gappy
PMD with DPS sampling as a function of the noise level.
The error rises sharply at $\eta_\% = 1$, beyond which the noise rather than the
reconstruction itself governs the accuracy.
It then grows approximately linearly with the noise level and remains below 10\% even at
$\eta_\% = 50$.
The reconstruction therefore does not break down even under severe corruption of the
measurements.

\subsection{Lid-driven cavity flow}

\begin{figure}[htbp]
   \centering
   \includegraphics[width=0.6\textwidth]{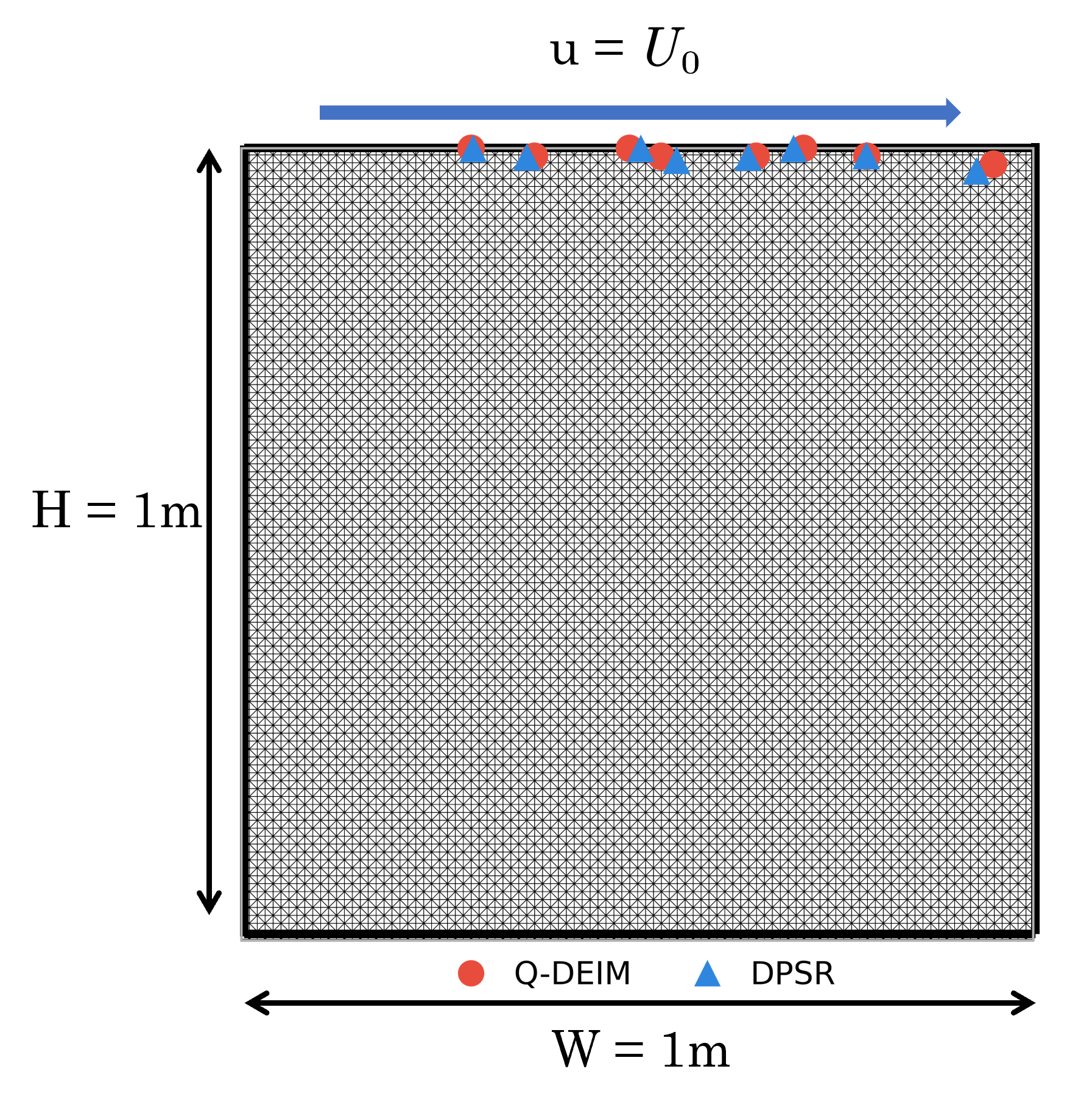}
   \caption{Computational domain, boundary conditions and mesh for the lid-driven cavity
   flow. Red circles mark the sampling points selected by QDEIM and blue triangles those
   selected by DPS.}
   \label{fig:cavity_domain}
\end{figure}

The second test case is a two-dimensional lid-driven cavity flow. The computational
domain is a $1\,\mathrm{m}\times 1\,\mathrm{m}$ square cavity discretized by a regular
triangular mesh with $101\times 101$ nodes, as shown in
Figure~\ref{fig:cavity_domain}. On the upper moving lid, the prescribed velocity
is given for $0\le x\le 1$ by
\begin{equation}
u_x^{\mathrm{lid}}(x,t)=10\,s(x)\,g(x,t), \qquad u_y^{\mathrm{lid}}(x,t)=0,
\label{eq:cavity-lid-velocity}
\end{equation}
where
\begin{equation}
s(x)=1-\frac{\exp(10(x-0.5))+\exp(-10(x-0.5))}
{\exp(5)+\exp(-5)}
\label{eq:cavity-spatial-envelope}
\end{equation}
and
\begin{equation}
g(x,t)=1+0.18\sin(2\pi x-\omega_1t+0.20)
      +0.09\sin(4\pi x-\omega_2t+0.70).
\label{eq:cavity-temporal-modulation}
\end{equation}
Here $\omega_1=2\pi/5.2$ and
$\omega_2=\tfrac{1}{2}(1+\sqrt{5})\,\omega_1$. No-slip boundary conditions are imposed on the remaining
walls. With the mean lid speed $U=10\,\mathrm{m/s}$ and the cavity length
$L=1\,\mathrm{m}$ as reference scales, the Reynolds number is $\mathrm{Re}\approx
3.0\times 10^3$. The flow snapshots used for reconstruction are sampled every
$0.05\,\mathrm{s}$. For Gappy PMD, the linear subspace dimension, nonlinear
manifold dimension, and number of sampling points are set to 2, 4, and 8,
respectively.

\begin{figure}[htbp]
   \centering
   \setlength{\tabcolsep}{2pt}
   \begin{tabular}{ccc}
      \includegraphics[width=0.32\textwidth]{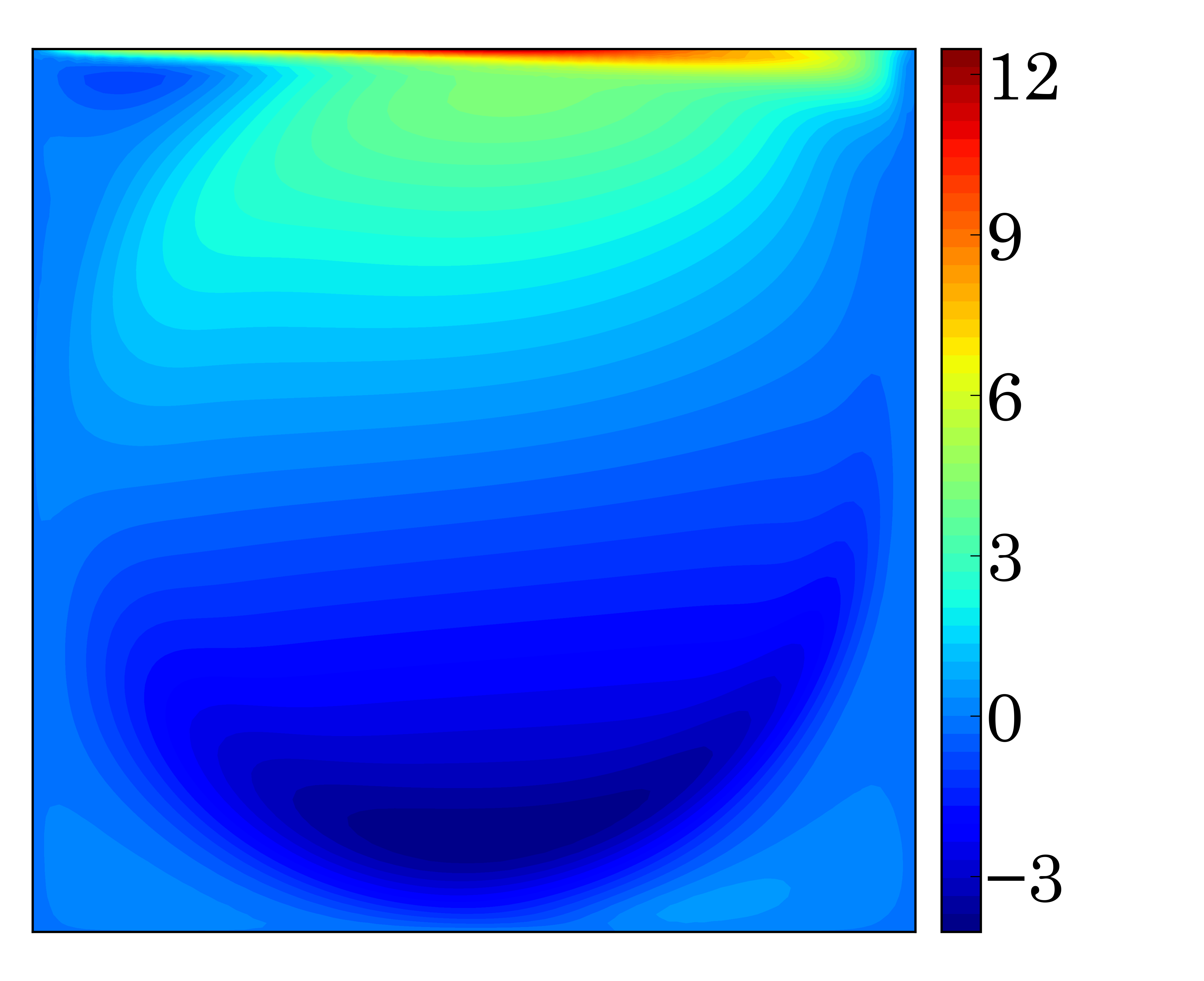}
      & \includegraphics[width=0.32\textwidth]{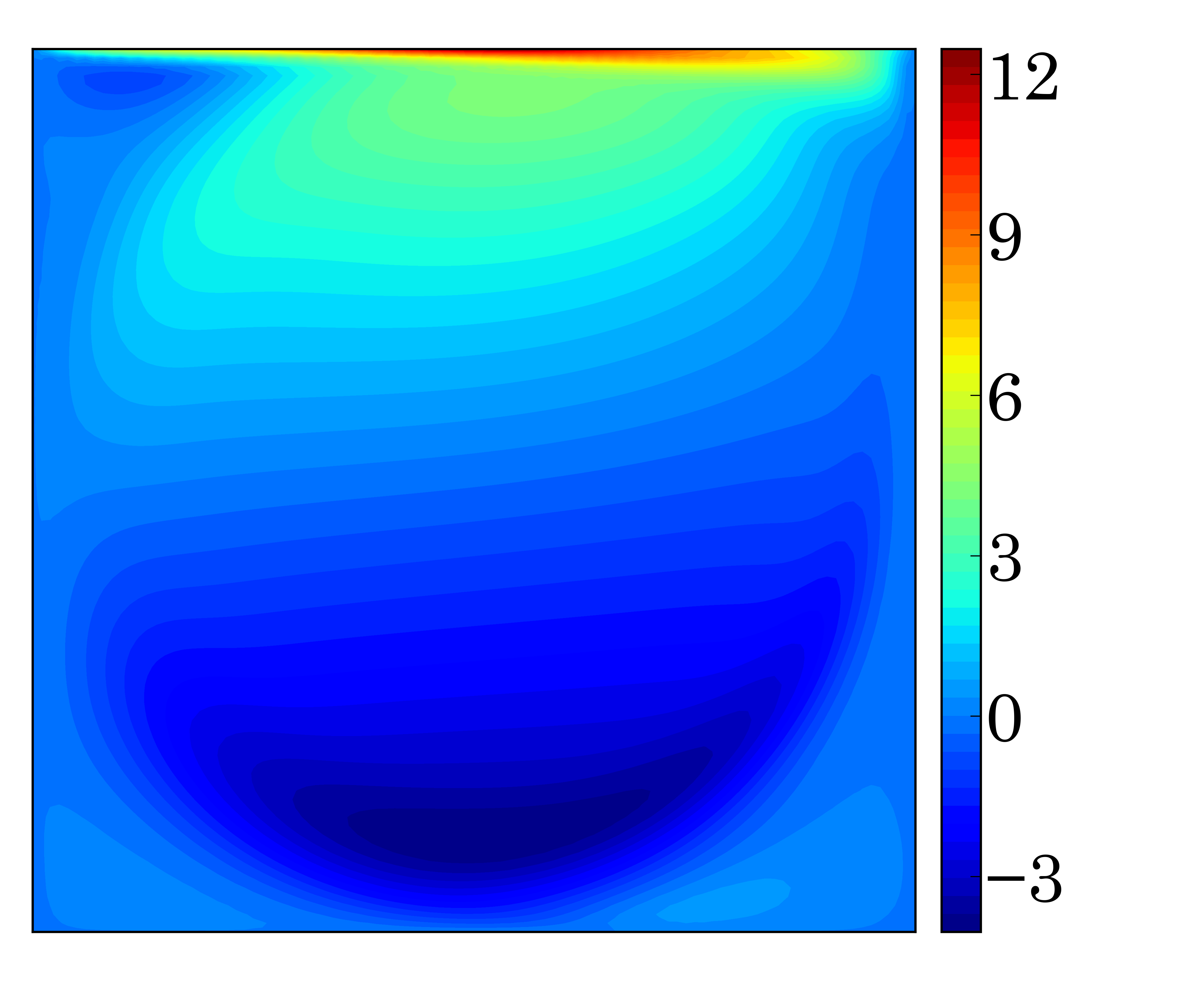}
      & \includegraphics[width=0.32\textwidth]{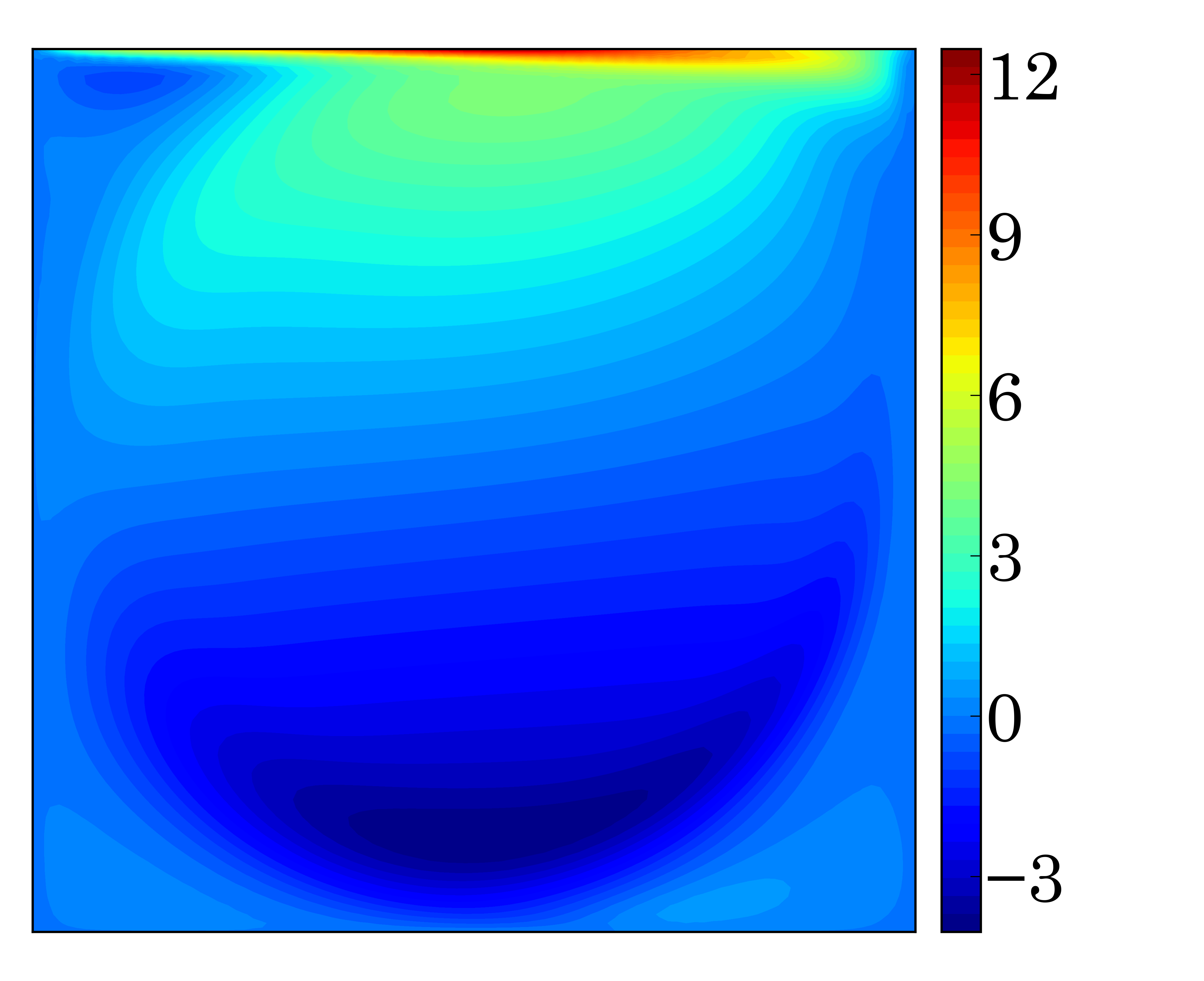}         \\
      \panellabel{(a) Reference solution, $t = 125\,\mathrm{s}$}
      & \panellabel{(b) Gappy PMD + QDEIM, $t = 125\,\mathrm{s}$}
      & \panellabel{(c) Gappy POD + QDEIM, $t = 125\,\mathrm{s}$}                  \\[6pt]
      & \includegraphics[width=0.32\textwidth]{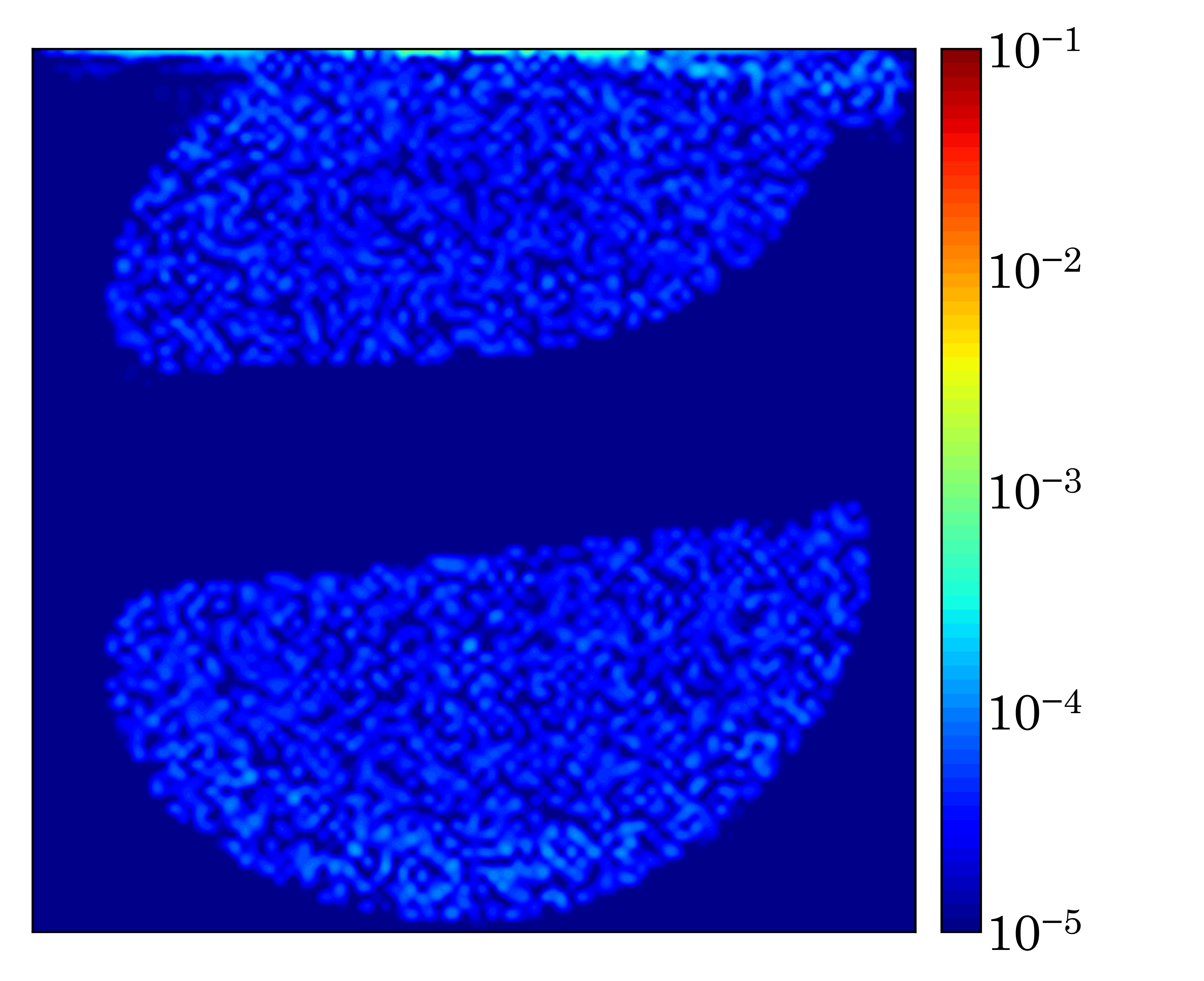}
      & \includegraphics[width=0.32\textwidth]{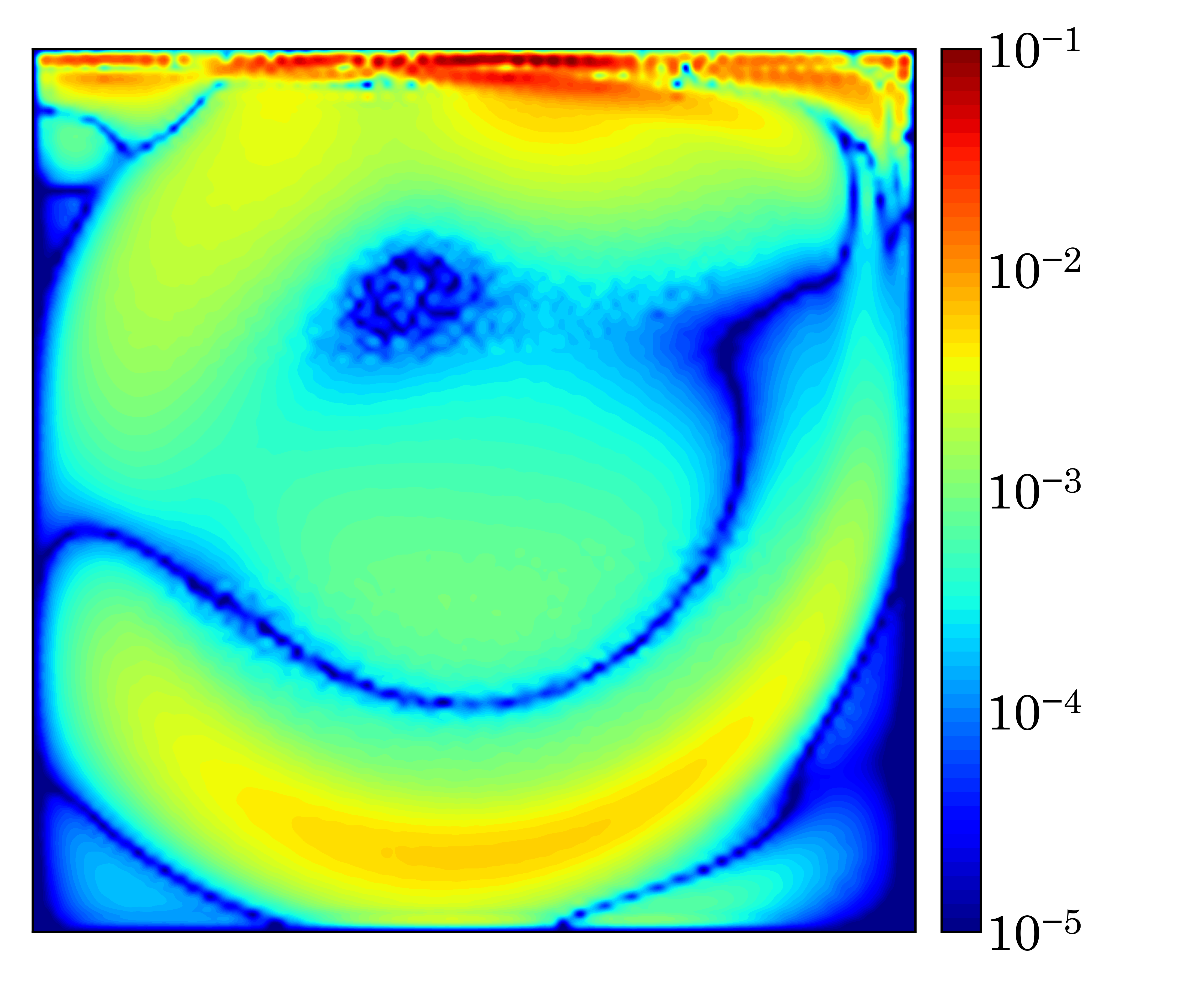}     \\
      & \panellabel{(d) Error of Gappy PMD + QDEIM, $t = 125\,\mathrm{s}$}
      & \panellabel{(e) Error of Gappy POD + QDEIM, $t = 125\,\mathrm{s}$}         \\[6pt]
      \includegraphics[width=0.32\textwidth]{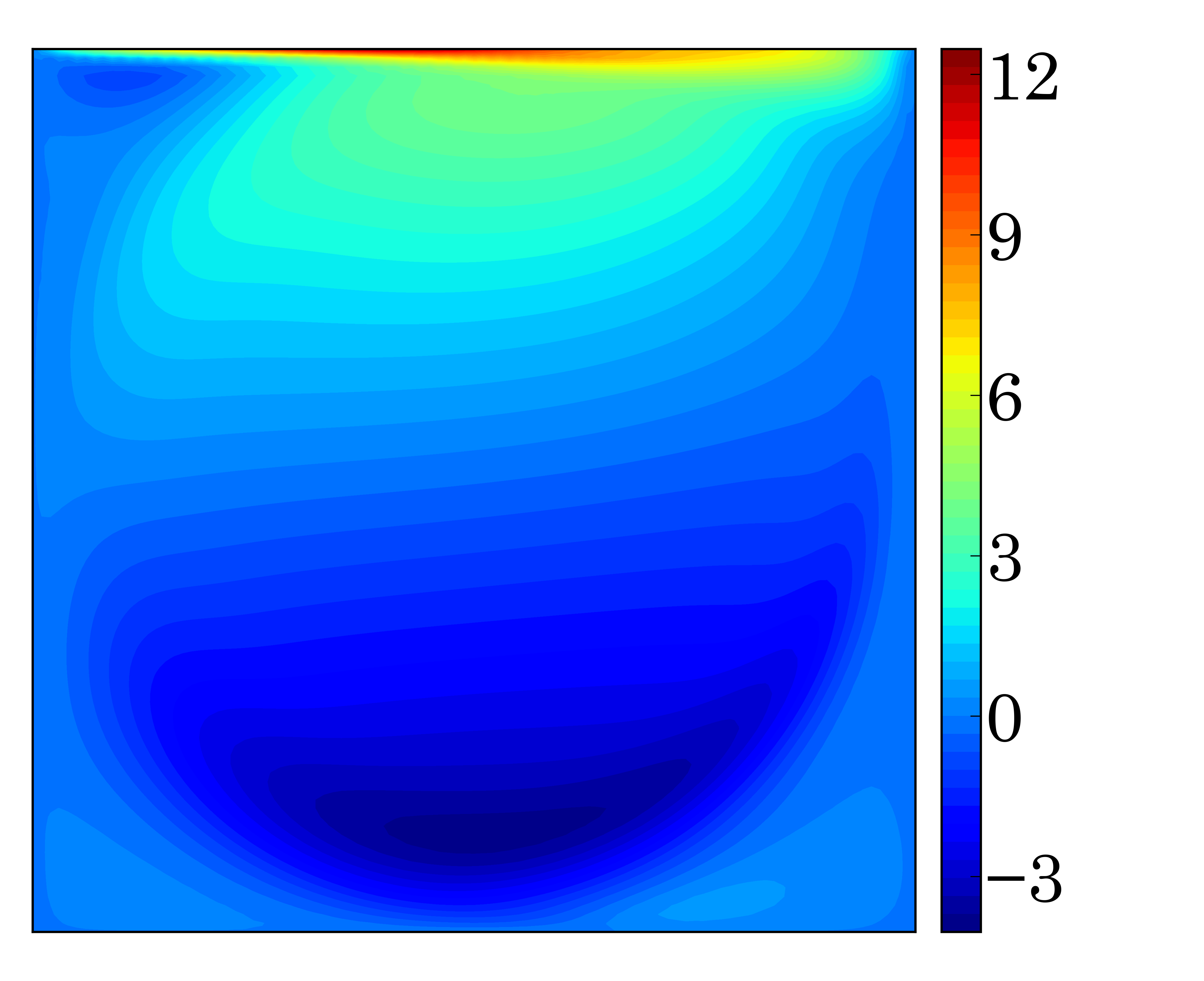}
      & \includegraphics[width=0.32\textwidth]{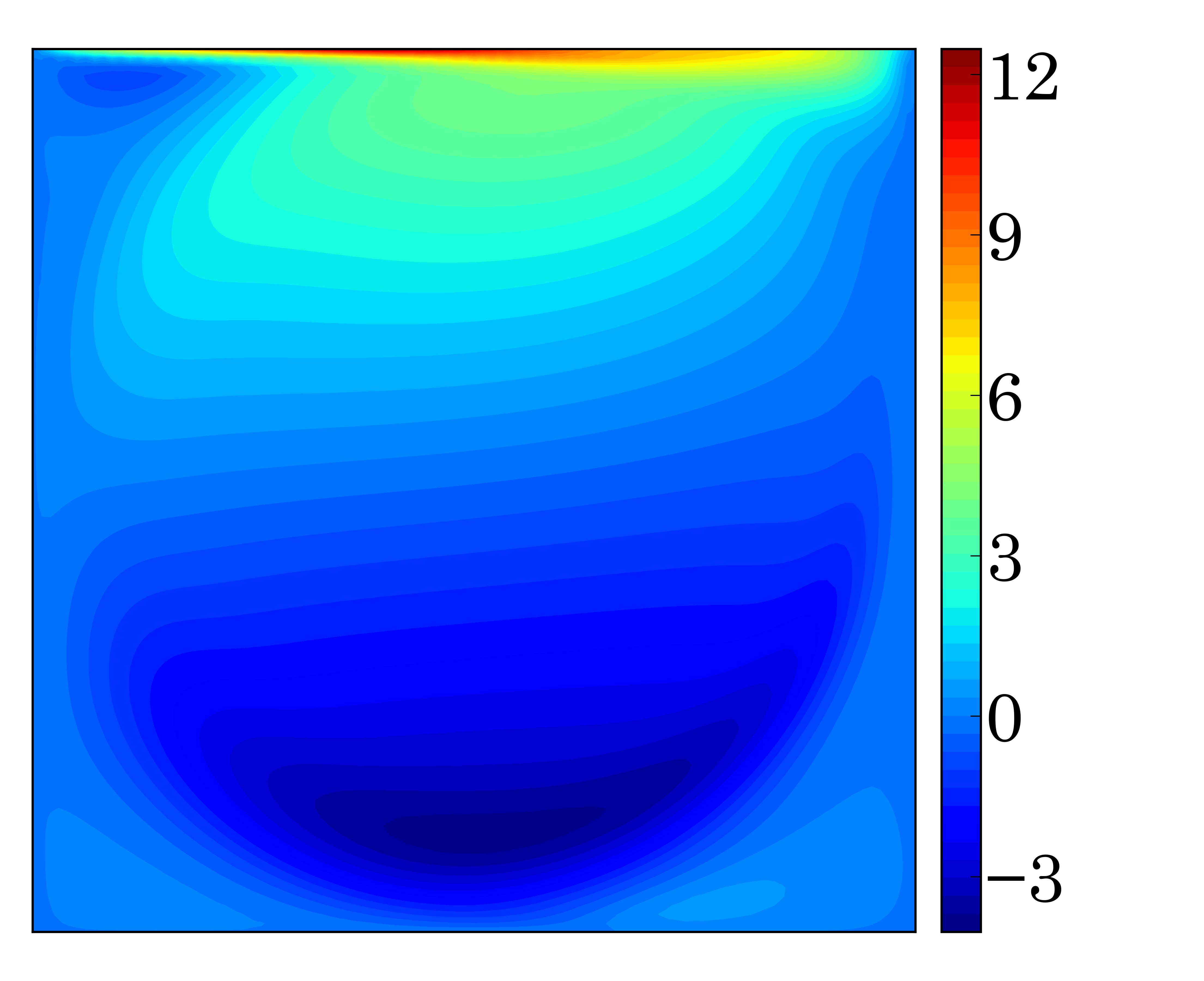}
      & \includegraphics[width=0.32\textwidth]{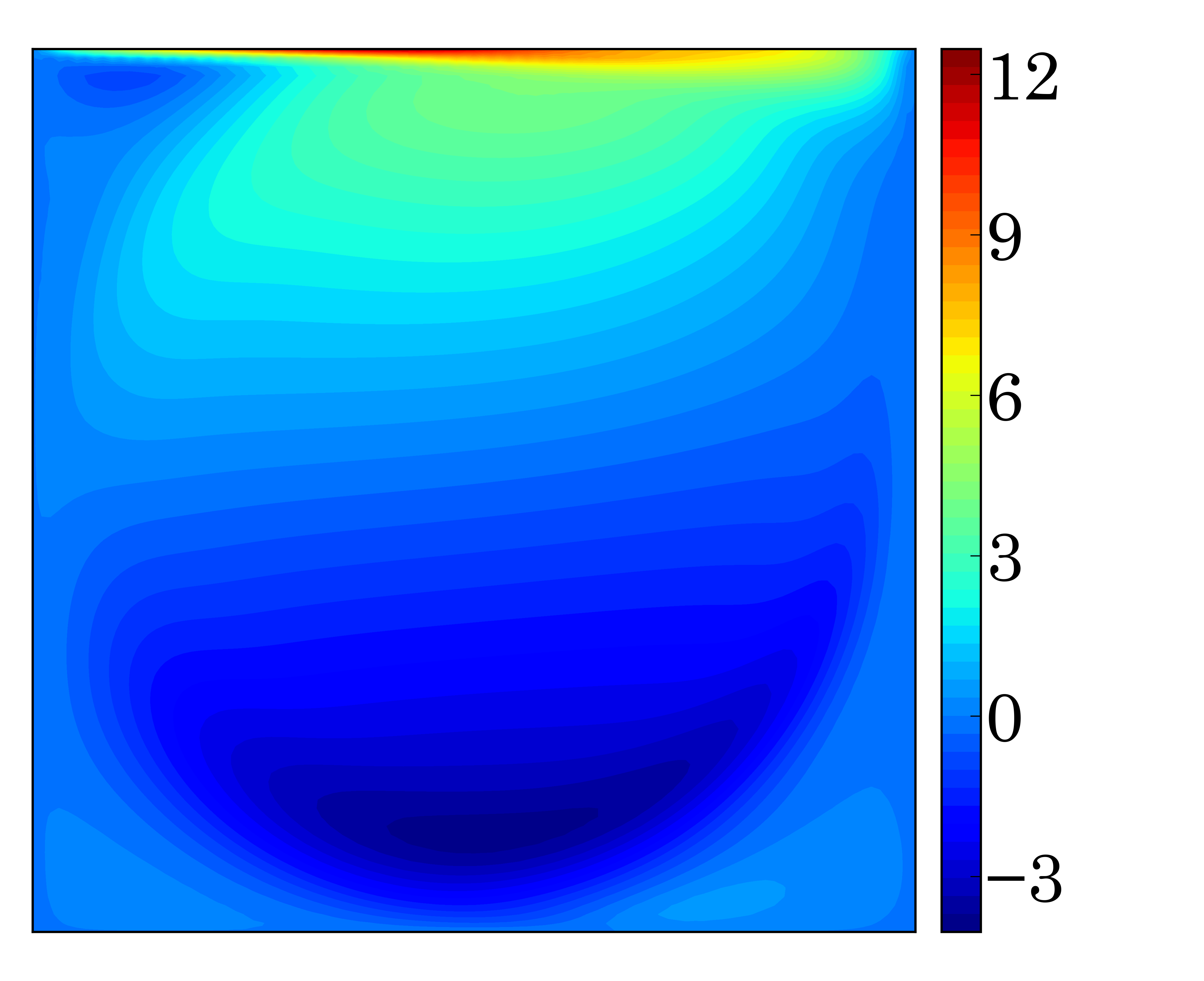}         \\
      \panellabel{(f) Reference solution, $t = 150\,\mathrm{s}$}
      & \panellabel{(g) Gappy PMD + QDEIM, $t = 150\,\mathrm{s}$}
      & \panellabel{(h) Gappy POD + QDEIM, $t = 150\,\mathrm{s}$}                  \\[6pt]
      & \includegraphics[width=0.32\textwidth]{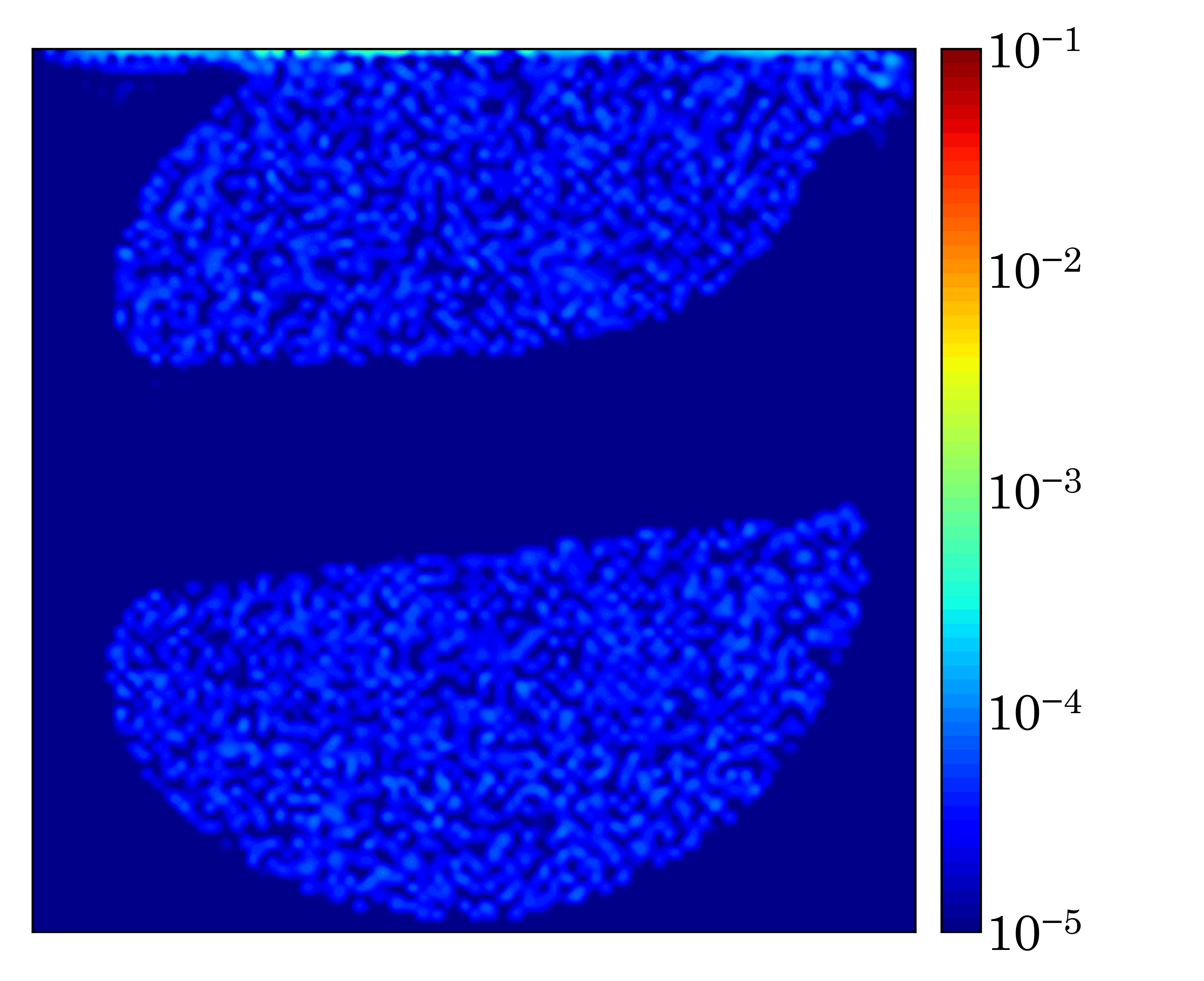}
      & \includegraphics[width=0.32\textwidth]{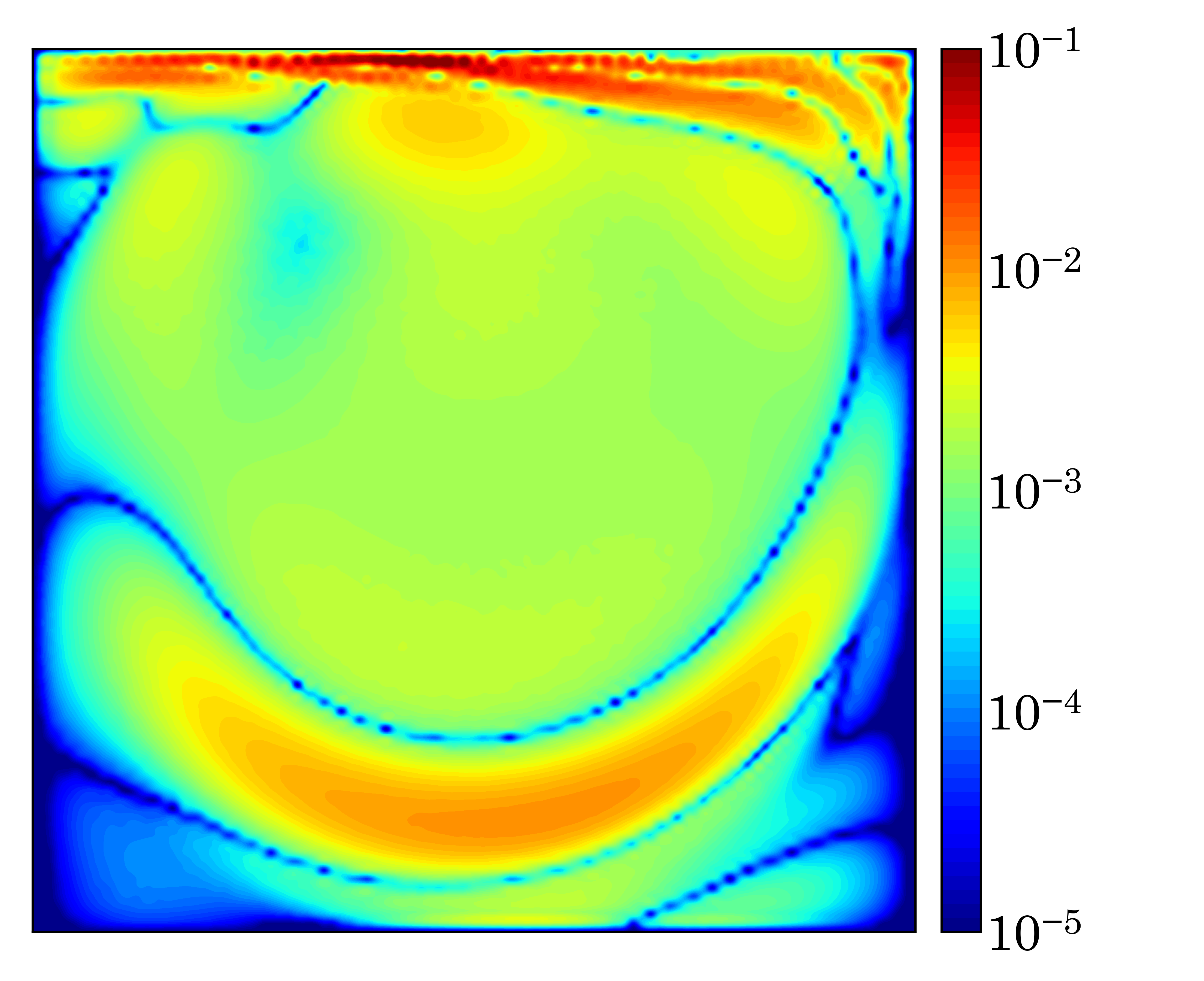}     \\
      & \panellabel{(i) Error of Gappy PMD + QDEIM, $t = 150\,\mathrm{s}$}
      & \panellabel{(j) Error of Gappy POD + QDEIM, $t = 150\,\mathrm{s}$}         \\
   \end{tabular}
   \caption{Reconstructed horizontal velocity $u_x$ and the corresponding absolute error
   fields for the lid-driven cavity flow, comparing Gappy PMD and Gappy POD.}
   \label{fig:cavity_recon_comparison}
\end{figure}

Figure~\ref{fig:cavity_recon_comparison} compares Gappy PMD and Gappy POD using
QDEIM sampling. Both methods reproduce the main structure of the horizontal
velocity field. The difference is more apparent in the error fields. Gappy PMD
gives small errors over most of the cavity, with localized errors in regions
where the reference field has larger spatial variation. In contrast, Gappy POD
produces broader error regions over the domain. Thus, under the same QDEIM
sampling, Gappy PMD is more accurate than Gappy POD.

\begin{figure}[htbp]
   \centering
   \setlength{\tabcolsep}{2pt}
   \begin{tabular}{ccc}
      \includegraphics[width=0.32\textwidth]{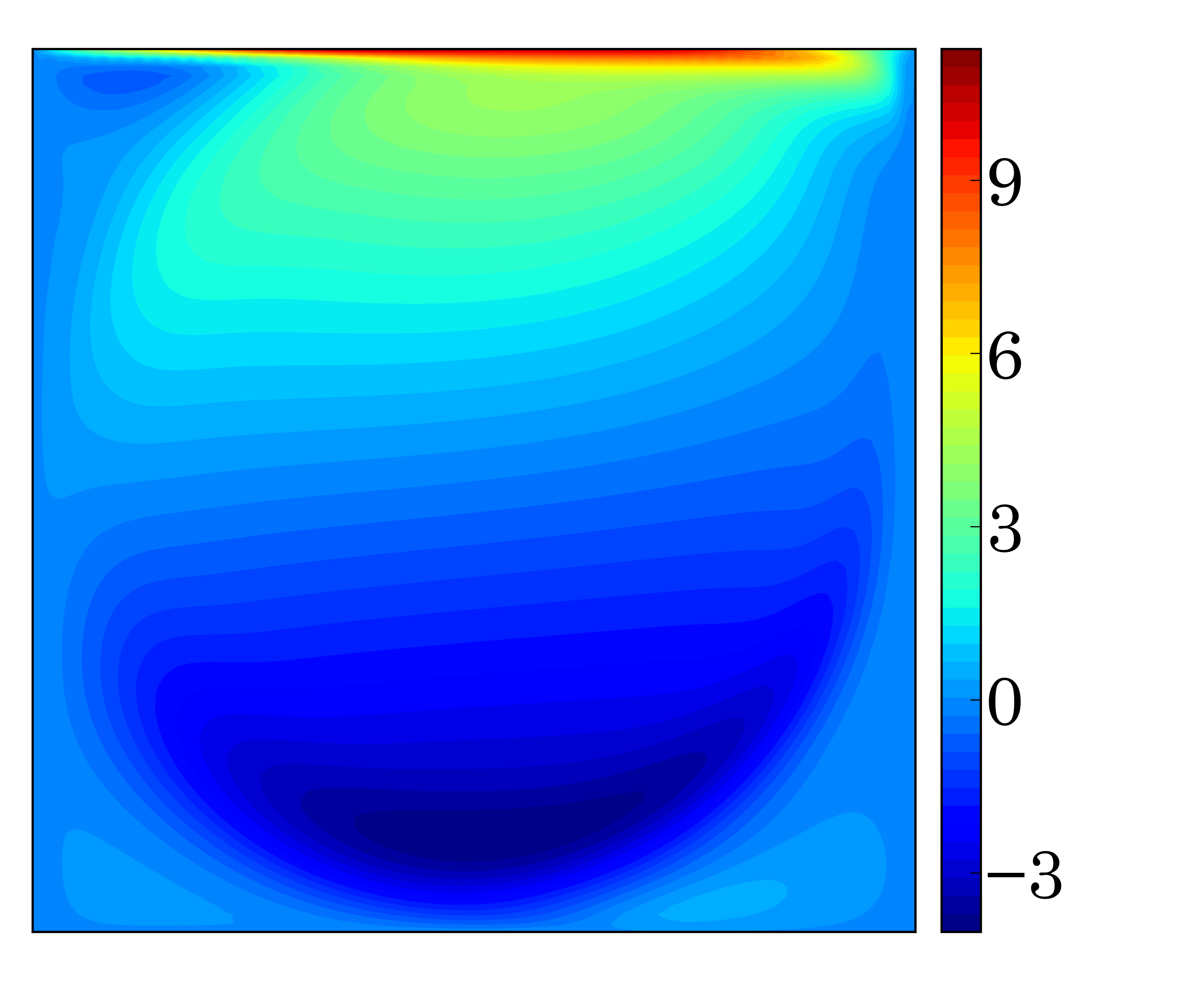}
      & \includegraphics[width=0.32\textwidth]{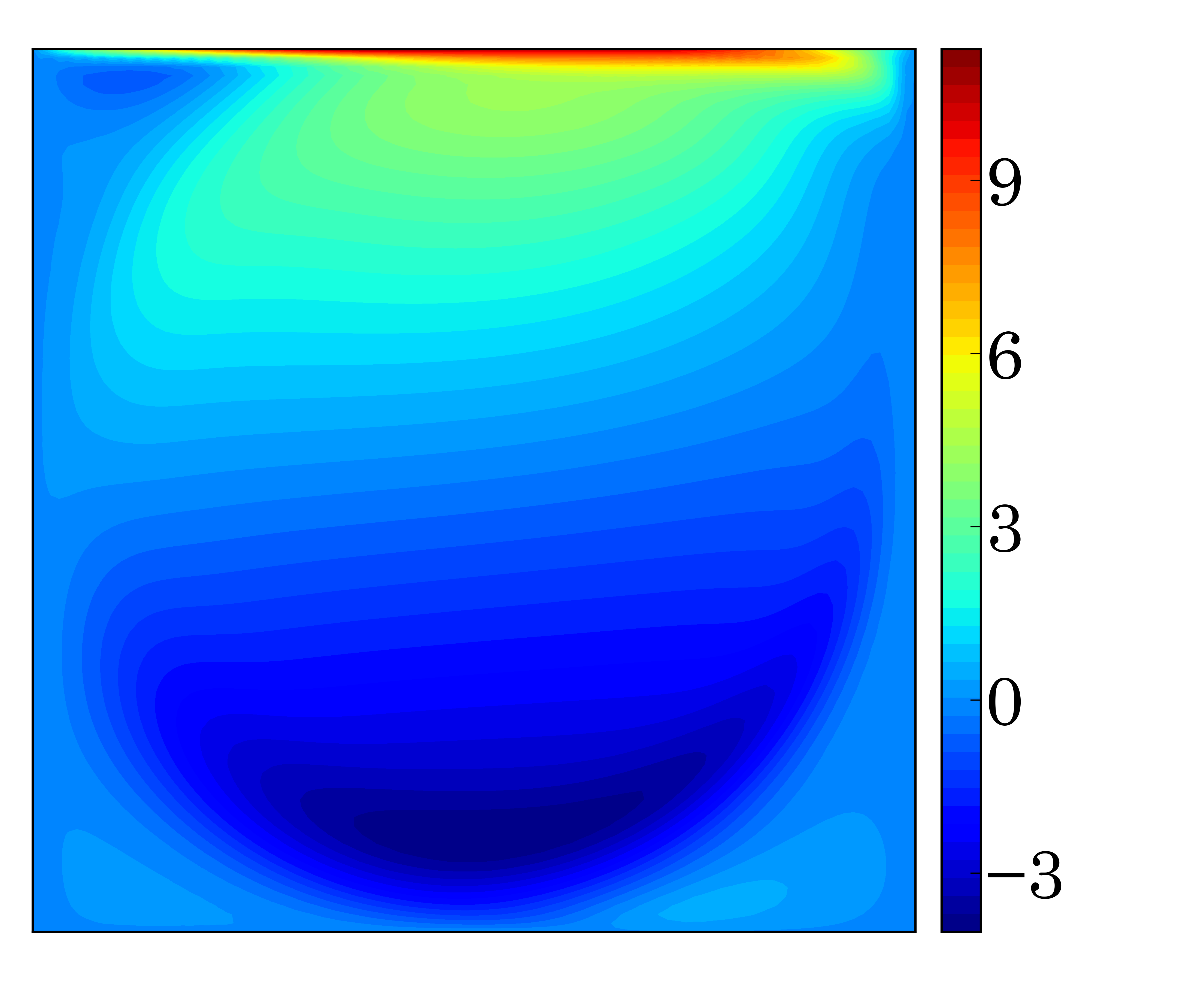}
      & \includegraphics[width=0.32\textwidth]{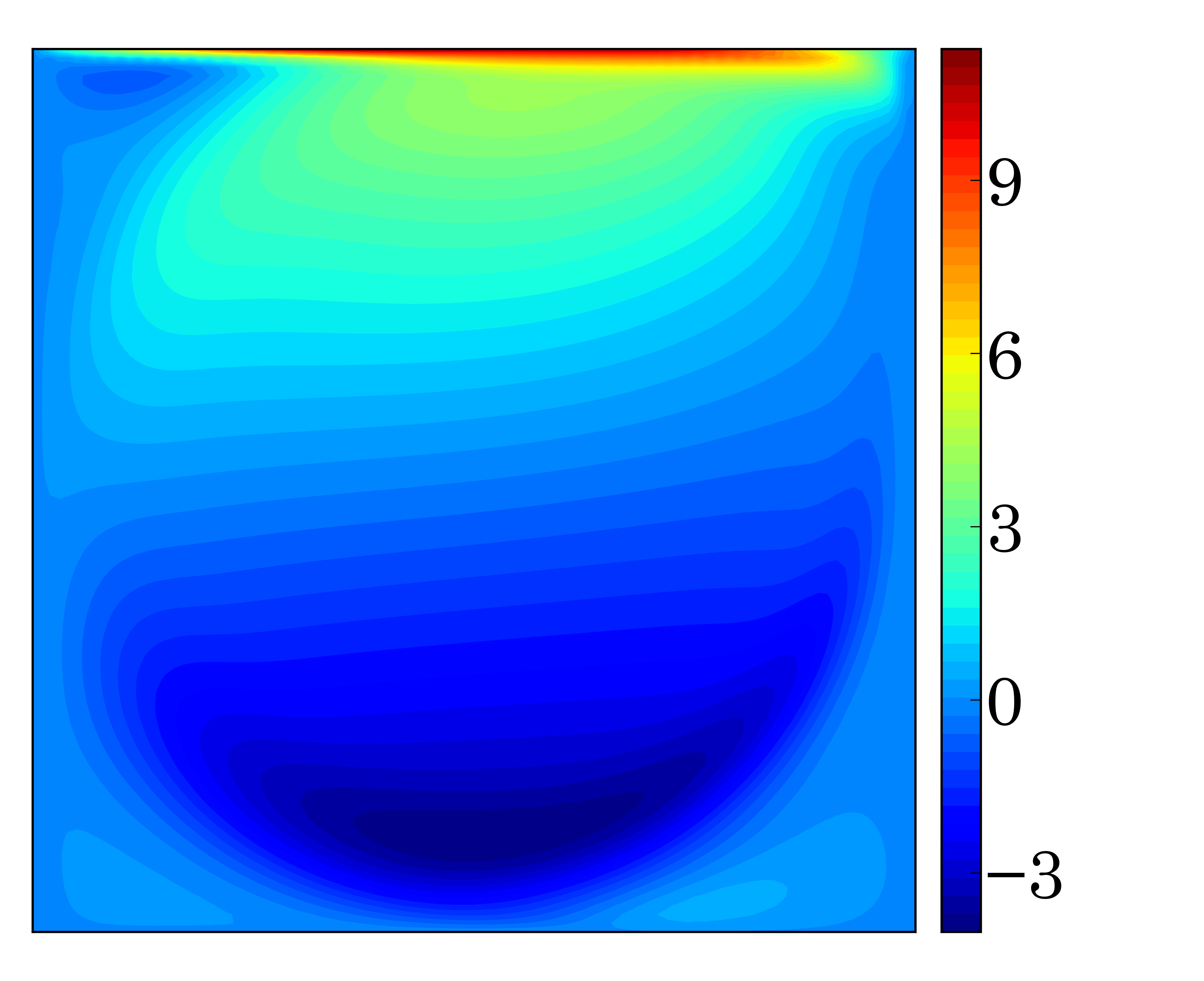}         \\
      \panellabel{(a) Reference solution, $t = 130\,\mathrm{s}$}
      & \panellabel{(b) Gappy PMD + DPS, $t = 130\,\mathrm{s}$}
      & \panellabel{(c) Gappy PMD + QDEIM, $t = 130\,\mathrm{s}$}                  \\[6pt]
      & \includegraphics[width=0.32\textwidth]{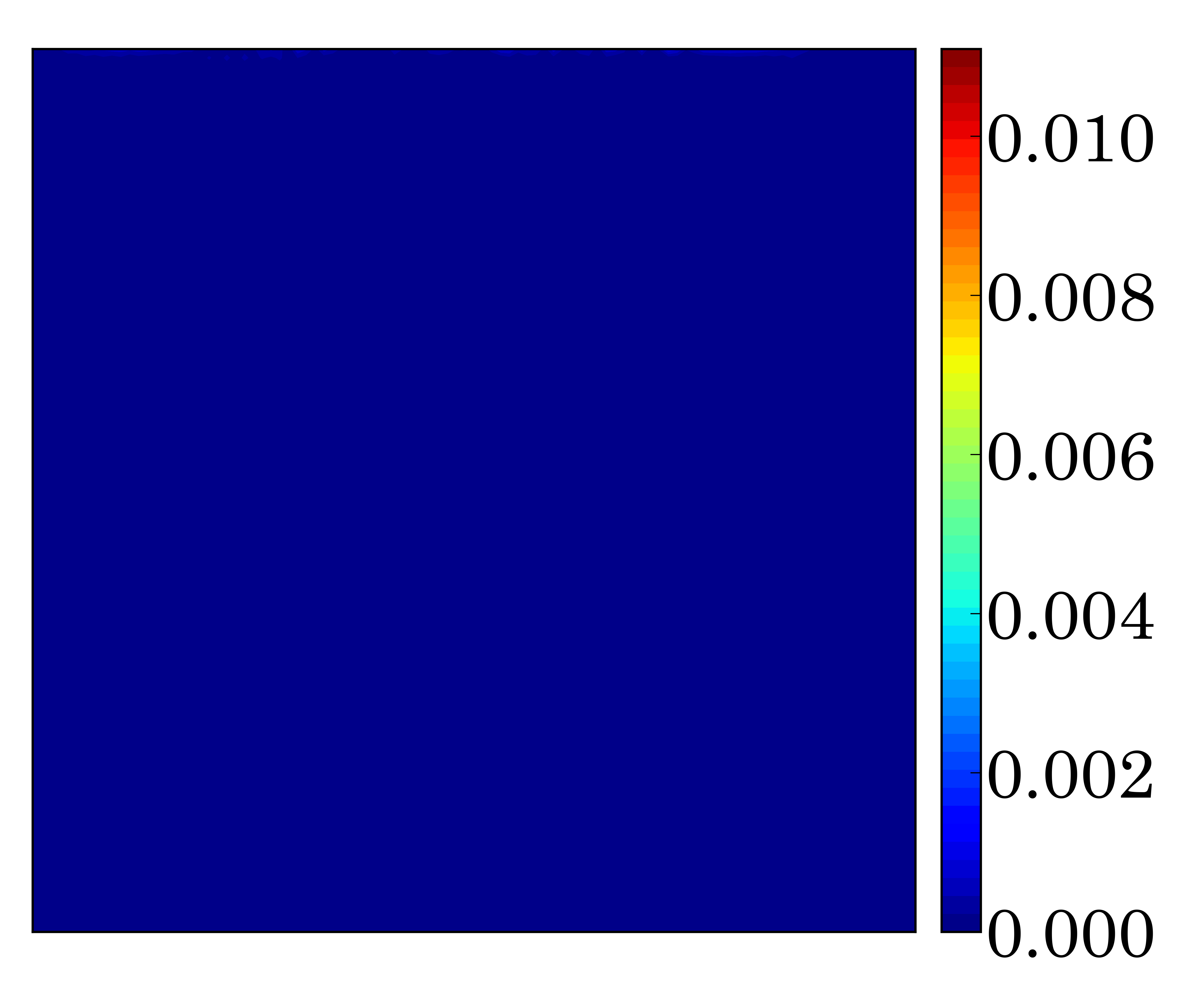}
      & \includegraphics[width=0.32\textwidth]{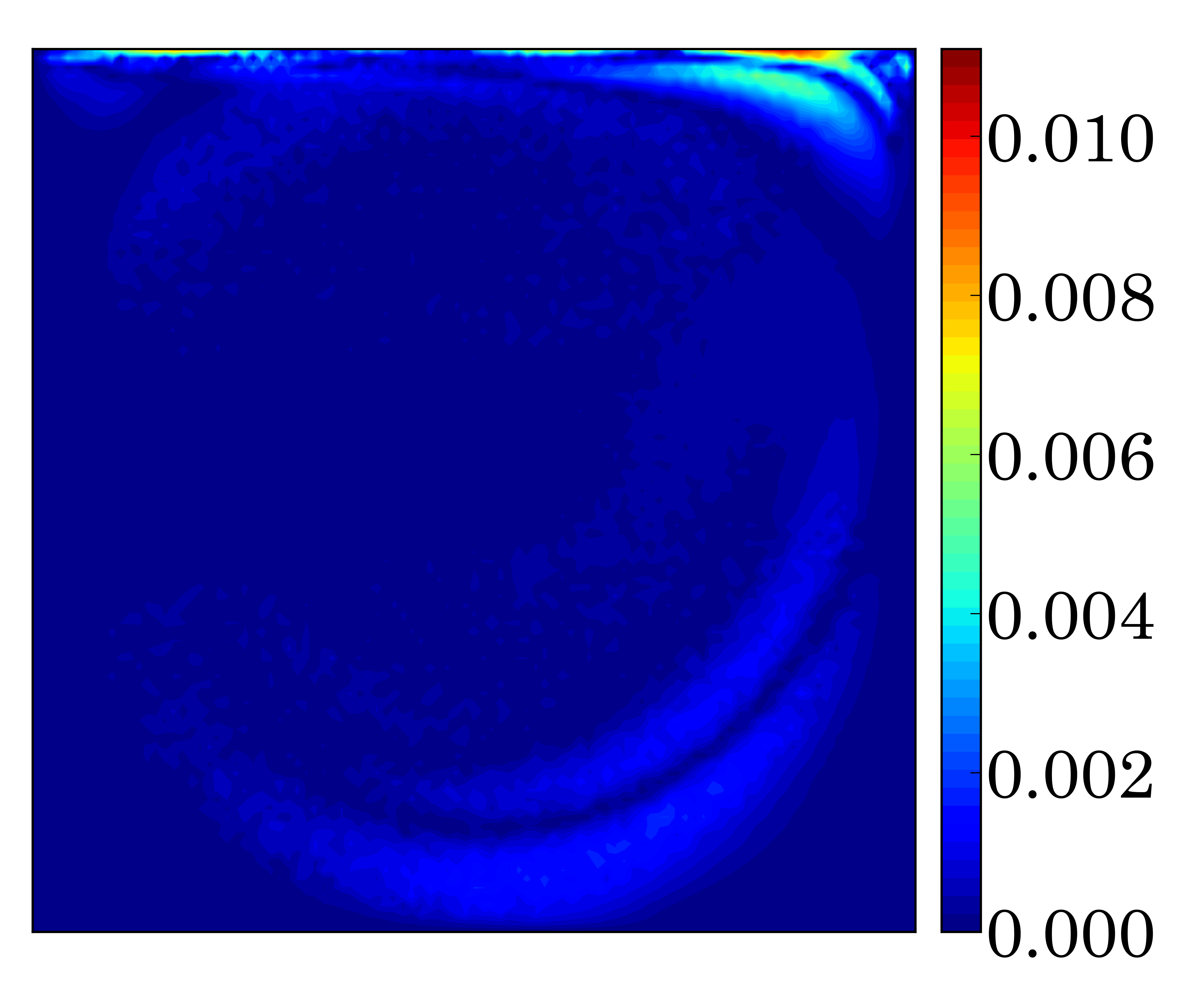}     \\
      & \panellabel{(d) Error of Gappy PMD + DPS, $t = 130\,\mathrm{s}$}
      & \panellabel{(e) Error of Gappy PMD + QDEIM, $t = 130\,\mathrm{s}$}         \\[6pt]
      \includegraphics[width=0.32\textwidth]{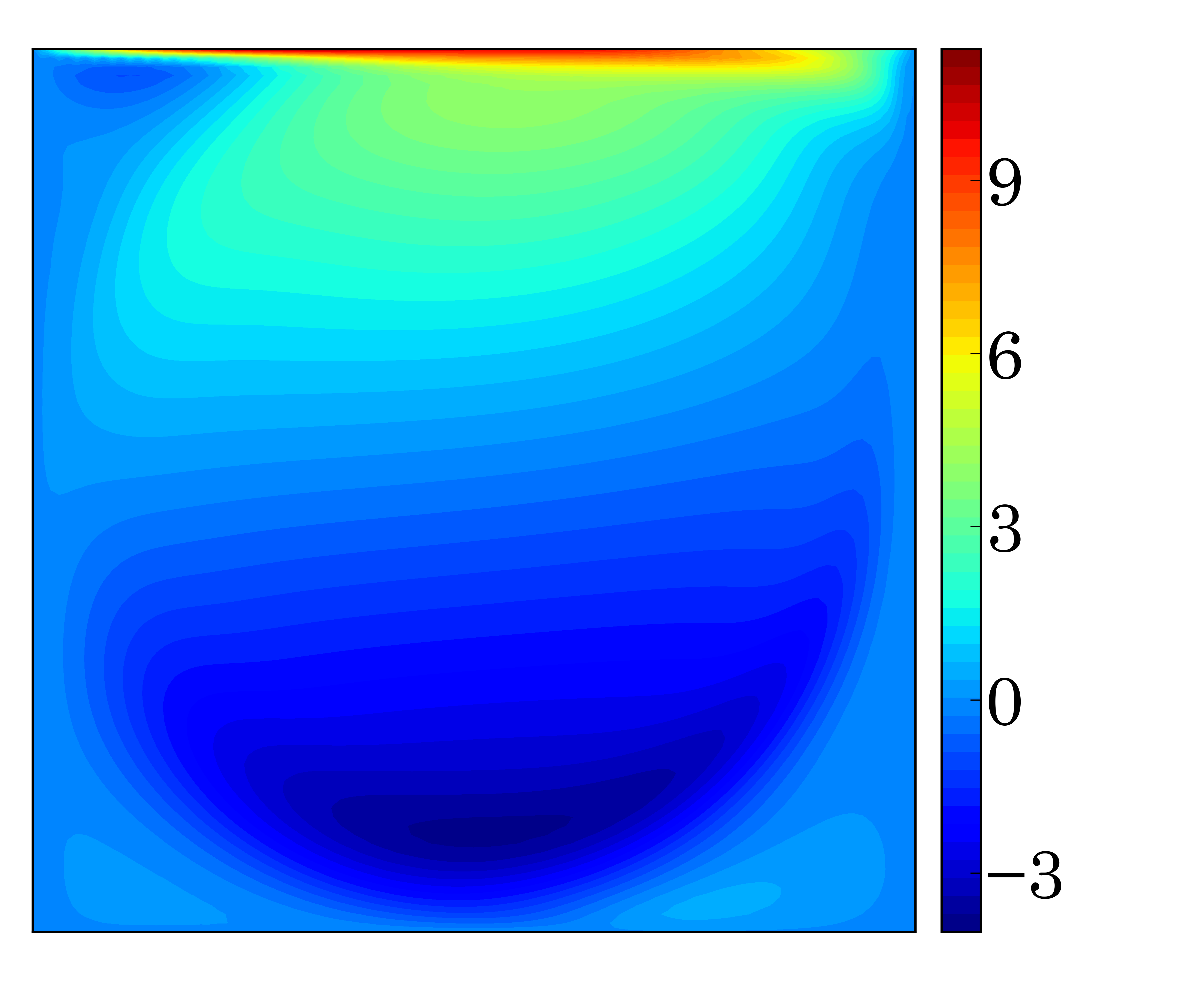}
      & \includegraphics[width=0.32\textwidth]{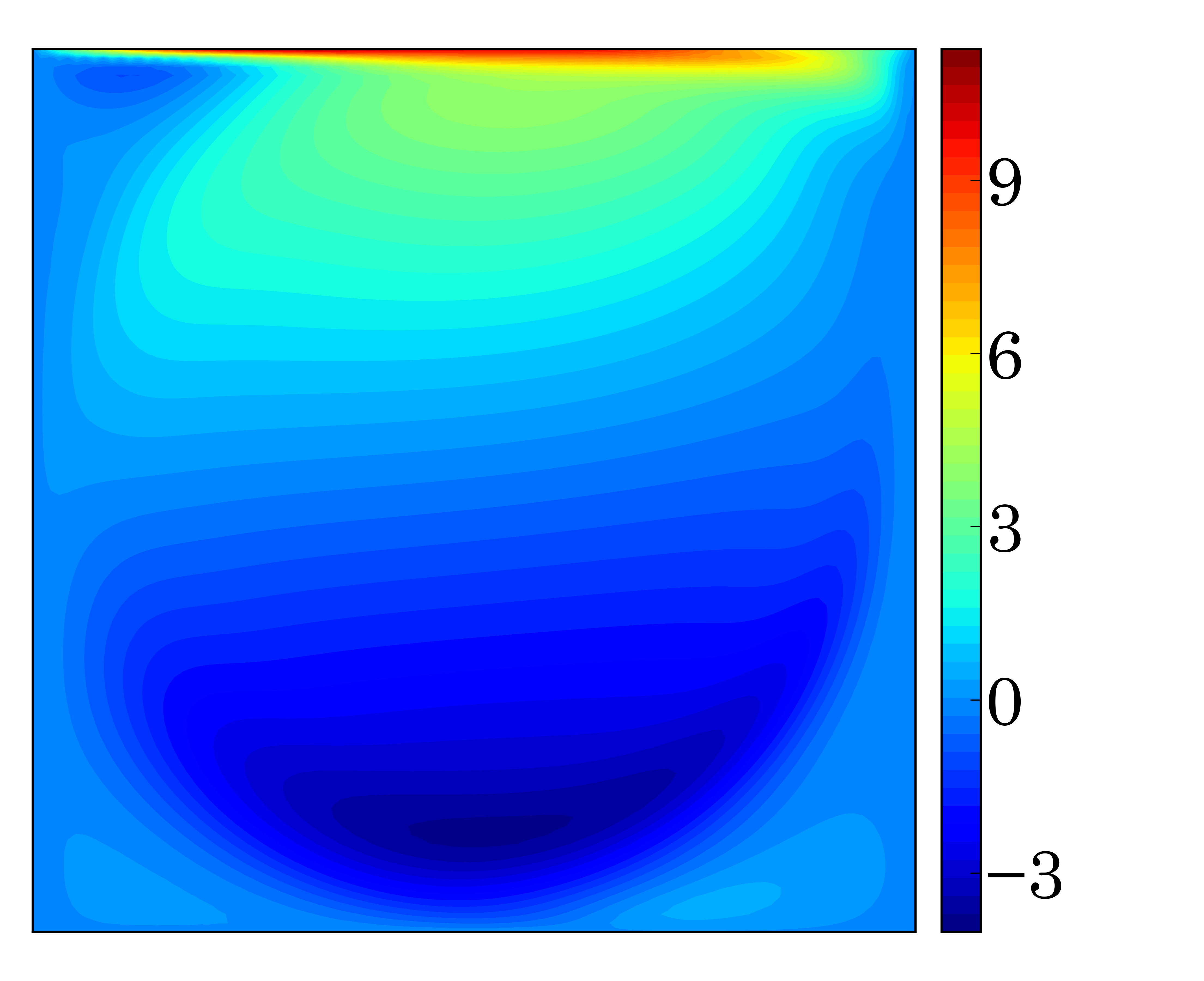}
      & \includegraphics[width=0.32\textwidth]{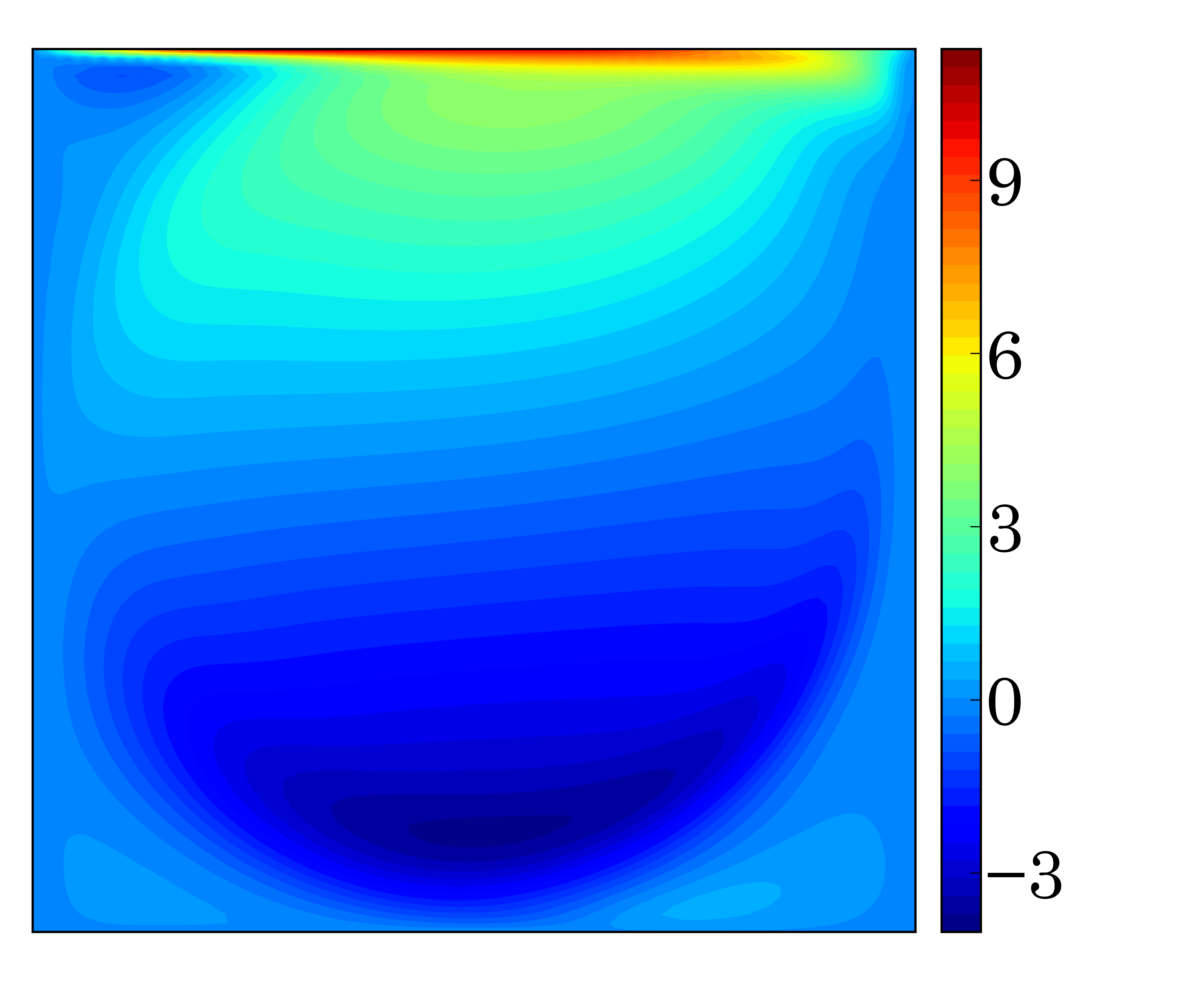}         \\
      \panellabel{(f) Reference solution, $t = 155\,\mathrm{s}$}
      & \panellabel{(g) Gappy PMD + DPS, $t = 155\,\mathrm{s}$}
      & \panellabel{(h) Gappy PMD + QDEIM, $t = 155\,\mathrm{s}$}                  \\[6pt]
      & \includegraphics[width=0.32\textwidth]{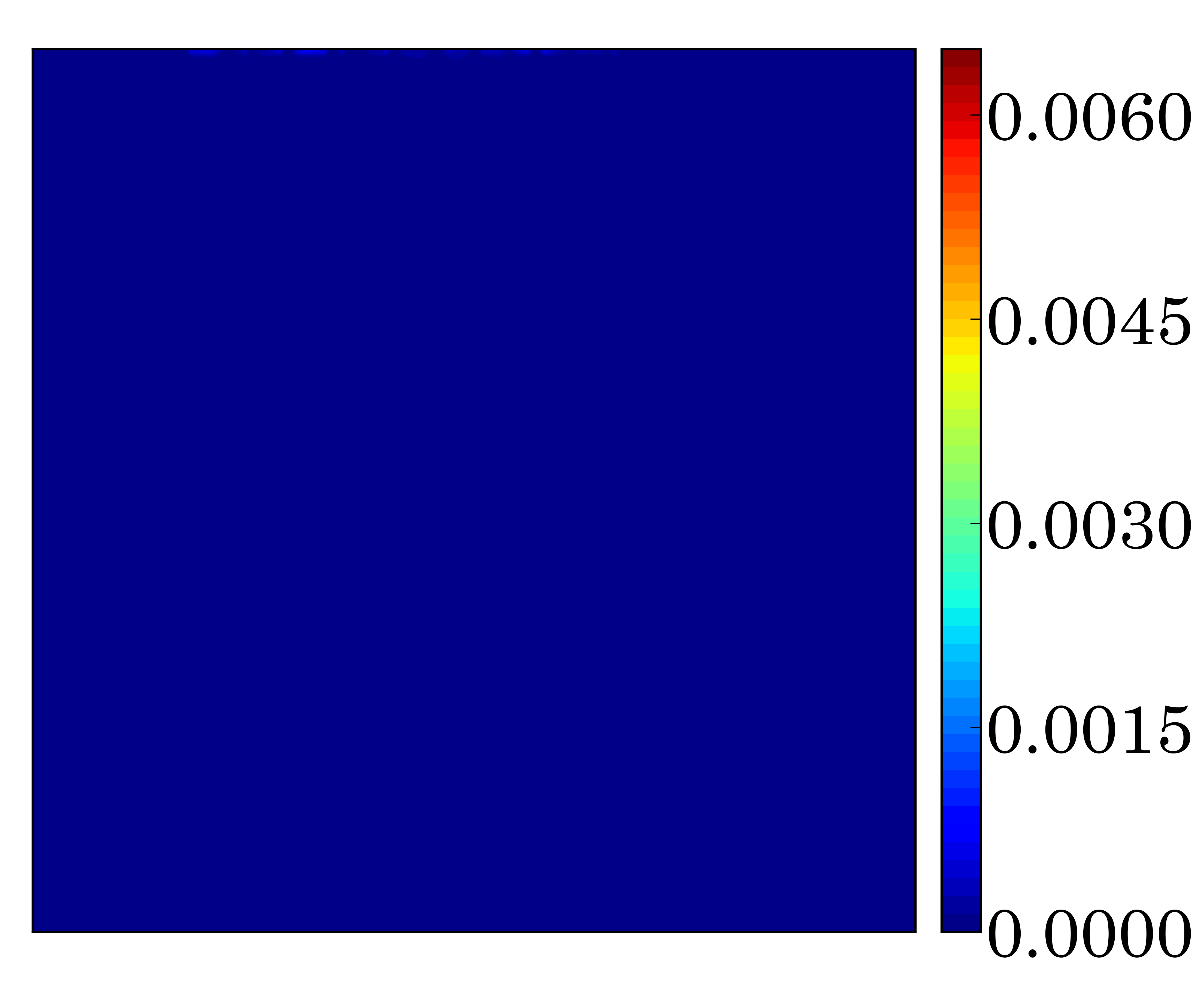}
      & \includegraphics[width=0.32\textwidth]{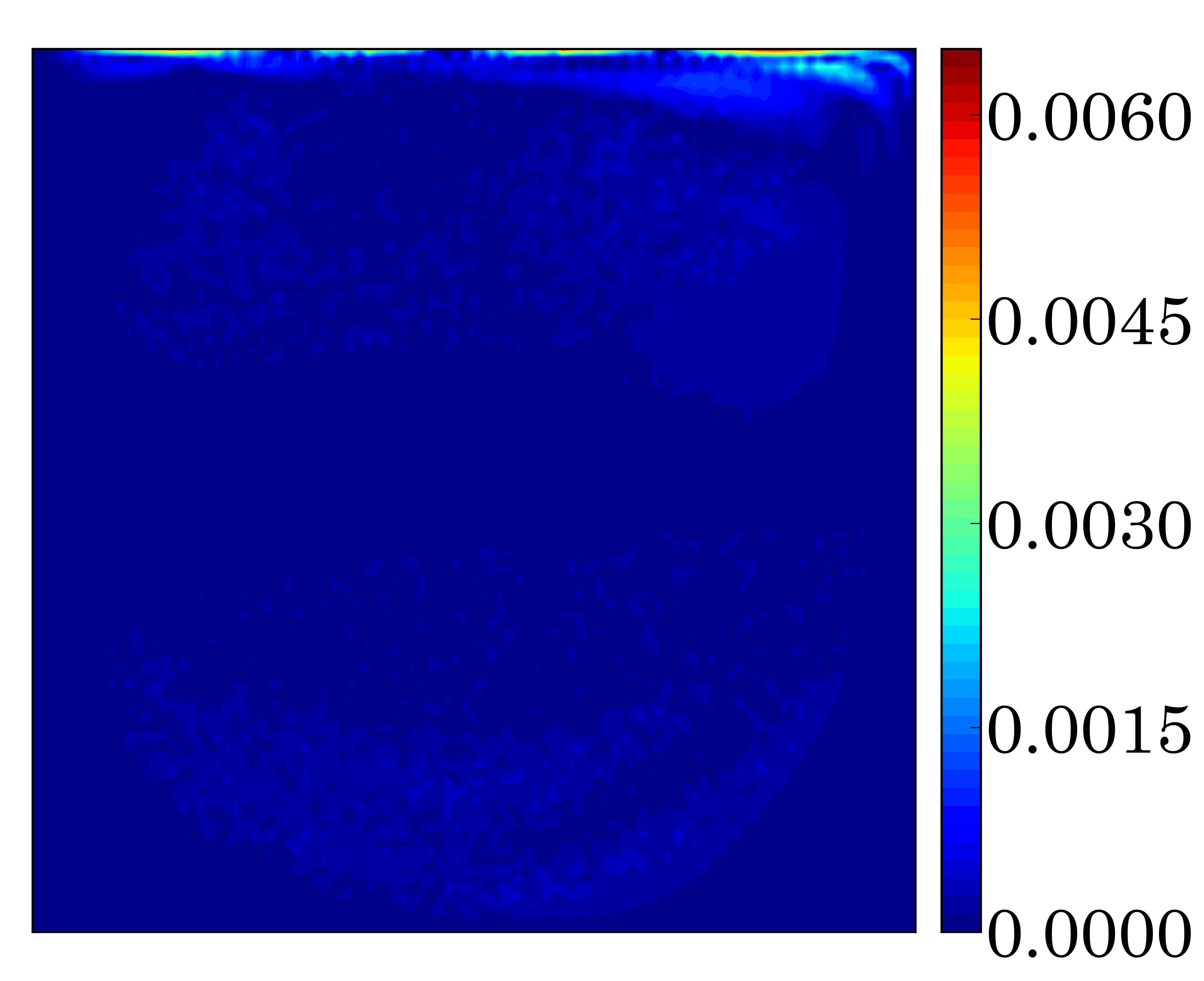}     \\
      & \panellabel{(i) Error of Gappy PMD + DPS, $t = 155\,\mathrm{s}$}
      & \panellabel{(j) Error of Gappy PMD + QDEIM, $t = 155\,\mathrm{s}$}         \\
   \end{tabular}
   \caption{Reconstructed horizontal velocity $u_x$ and the corresponding absolute error
   fields for the lid-driven cavity flow, comparing QDEIM and DPS sampling under Gappy
   PMD.}
   \label{fig:cavity_sampling_comparison}
\end{figure}

Figure~\ref{fig:cavity_sampling_comparison} compares DPS and QDEIM within the
Gappy PMD framework. At the level of the reconstructed fields, the two sampling
strategies give comparable results. The error fields show a clearer distinction,
with DPS producing weaker and less structured errors, especially in regions
where QDEIM gives localized error concentrations. These results demonstrate the
benefit of using the reconstruction error itself, rather than a criterion based
only on the reduced basis, to guide point selection in Gappy PMD.

\begin{figure}[htbp]
   \centering
   \includegraphics[width=0.6\textwidth]{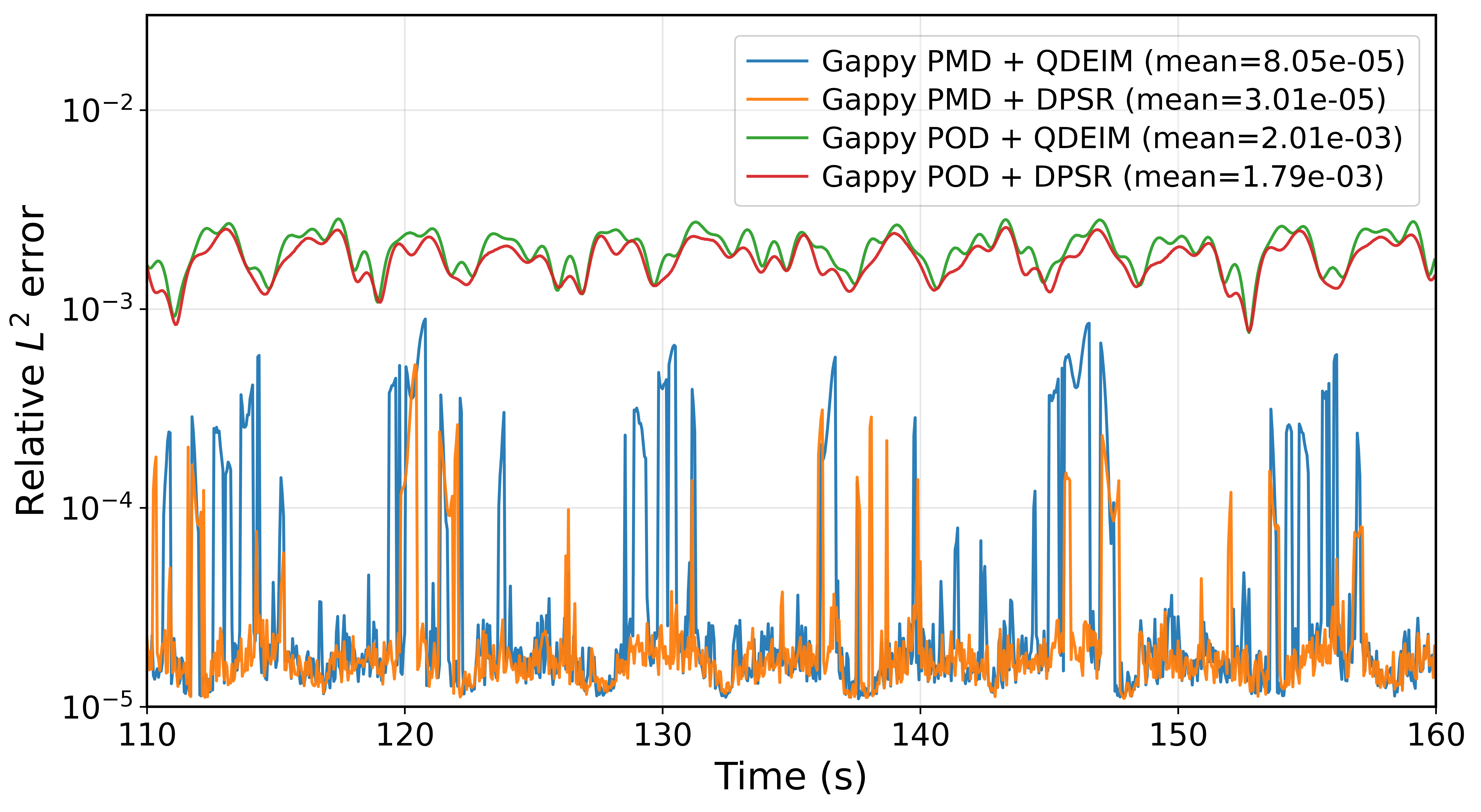}
   \caption{Relative $L^2$ reconstruction error over the test snapshots for the
   lid-driven cavity flow, for the four combinations of reconstruction method and
   sampling strategy.}
   \label{fig:cavity_error_curve}
\end{figure}

The error histories in Figure~\ref{fig:cavity_error_curve} agree with the
instantaneous comparisons. Over the full test interval, the Gappy PMD curves
remain about two orders of magnitude below the corresponding Gappy POD curves.
For Gappy PMD, replacing QDEIM with DPS lowers the mean error by more than half
and also reduces transient error peaks. These results indicate that the
PMD representation provides the main accuracy gain, while DPS gives
an additional improvement by matching the sampling criterion to the
reconstruction objective.

\begin{figure}[htbp]
   \centering
   \includegraphics[width=0.65\textwidth]{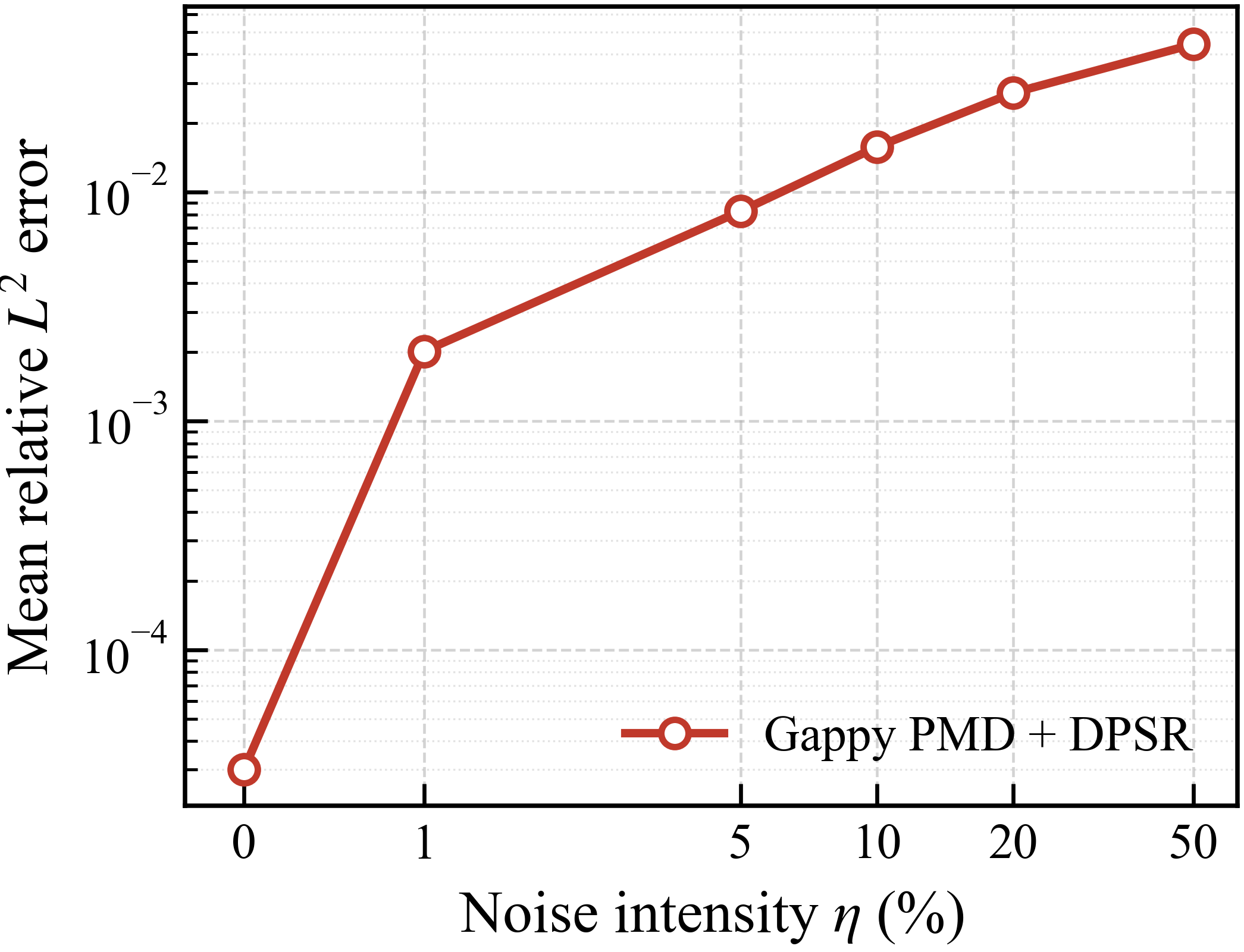}
   \caption{Effect of the noise level $\eta_\%$ on the mean relative $L^2$ error of Gappy PMD
   with DPS sampling for the lid-driven cavity flow, averaged over the test snapshots and
   over five independent noise realizations.}
   \label{fig:cavity_noise_robustness}
\end{figure}

Finally, Figure~\ref{fig:cavity_noise_robustness} reports the response of Gappy
PMD with DPS sampling to noisy measurements. The mean relative error increases
monotonically with the prescribed noise level on logarithmic axes. The increase
is gradual, indicating that the reconstruction error is progressively controlled
by the measurement noise rather than by an abrupt loss of stability. Even for
strongly corrupted measurements, the method maintains a bounded reconstruction
error. This behavior supports the robustness of Gappy PMD for sparse and noisy
measurements.

\subsection{Backward-facing step flow}

\begin{figure}[htbp]
   \centering
   \includegraphics[width=0.85\textwidth]{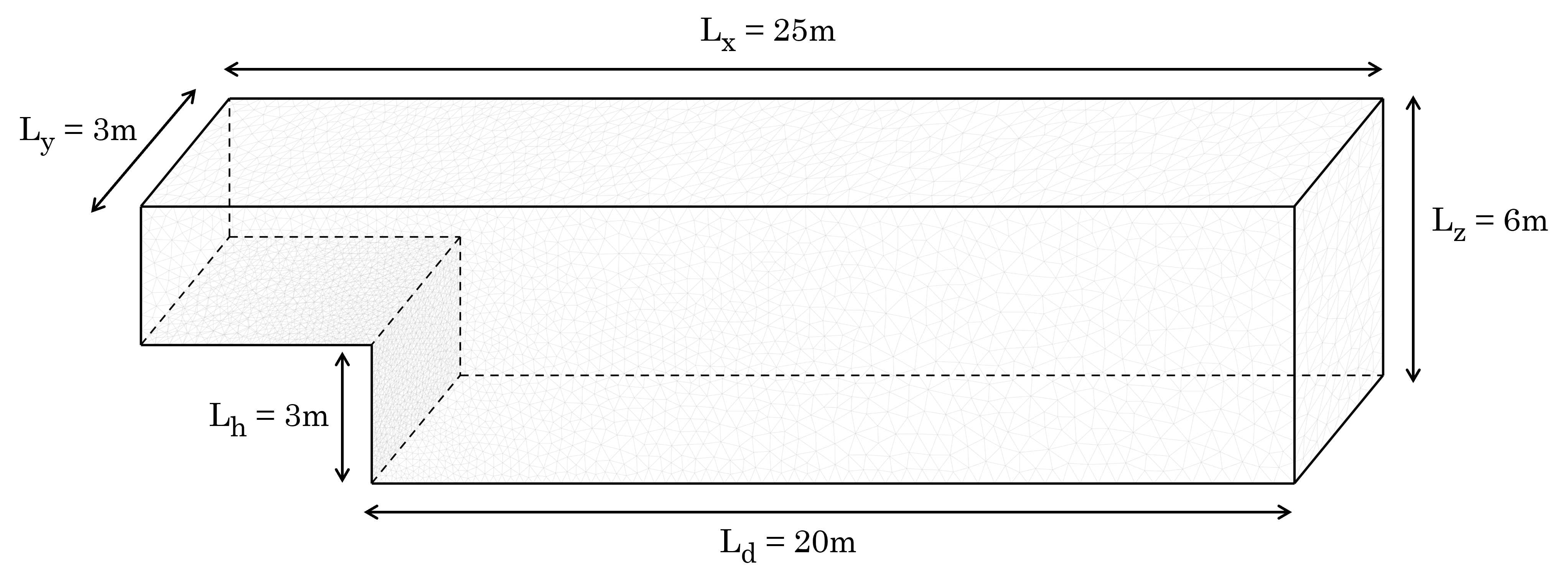}
   \caption{Computational domain, boundary conditions and mesh for the backward-facing step
   flow.}
   \label{fig:backward_domain}
\end{figure}

We finally consider the three-dimensional backward-facing step flow shown in
Figure~\ref{fig:backward_domain}. The step height is $h=3\,\mathrm{m}$, and the
downstream channel has length $20\,\mathrm{m}$, height $6\,\mathrm{m}$, and
spanwise width $3\,\mathrm{m}$. On the inlet, the prescribed velocity is
\begin{equation}
u_x(z,t)=4U_0\bar{z}(1-\bar{z})
\left[1+0.12\sin\left(\frac{2\pi t}{8}\right)\right],
\qquad
\bar{z}=\frac{z-3}{3},\quad 3\le z\le 6,
\end{equation}
with $u_y=u_z=0$ and $U_0=1.0\,\mathrm{m/s}$. The viscosity is set to
$\nu=0.00167\,\mathrm{m^2/s}$. The Reynolds number based on the step height
and the mean inlet velocity is $\mathrm{Re}_h=\frac{U_bh}{\nu}\approx 1200$,
where $U_b=\frac{2U_0}{3}$ for the prescribed parabolic inlet profile. No-slip
conditions are imposed on the solid walls, and free-slip
conditions on the spanwise side walls and the top boundary. Snapshots are
stored every
$0.1\,\mathrm{s}$. The streamwise velocity $u_x$ is reconstructed on an
unstructured finite-element mesh with 18,058 nodes. For Gappy PMD, the linear
subspace dimension is $r=2$ and the nonlinear manifold dimension is $r_1=5$.
The number of sampling points is $q=9$. The Gappy POD basis dimension is set
equal to the total PMD dimension $r+r_1=7$.

\begin{figure}[htbp]
   \centering
   \begin{tabular}{cc}
      \includegraphics[width=0.48\textwidth]{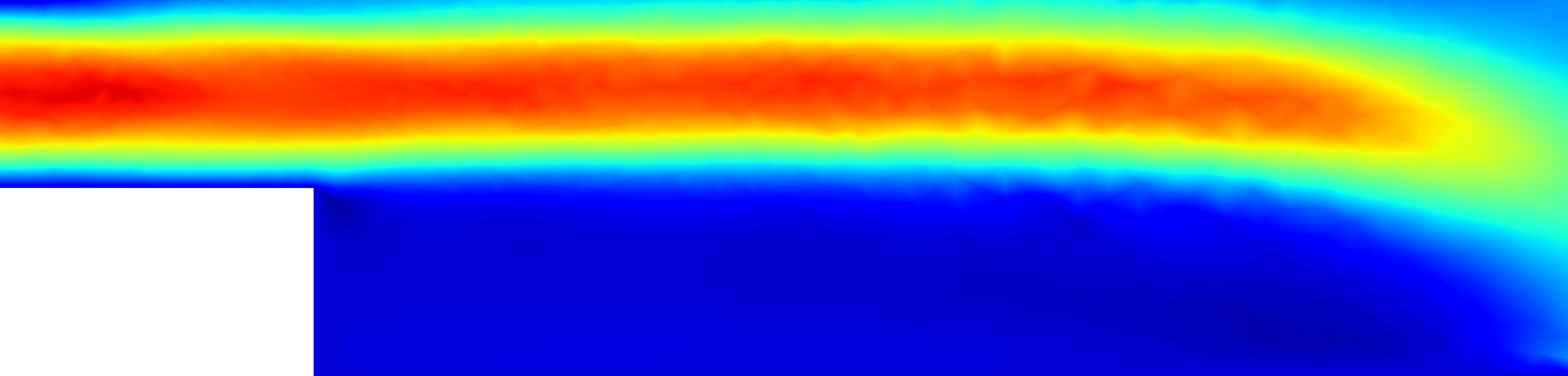}
      & \includegraphics[width=0.48\textwidth]{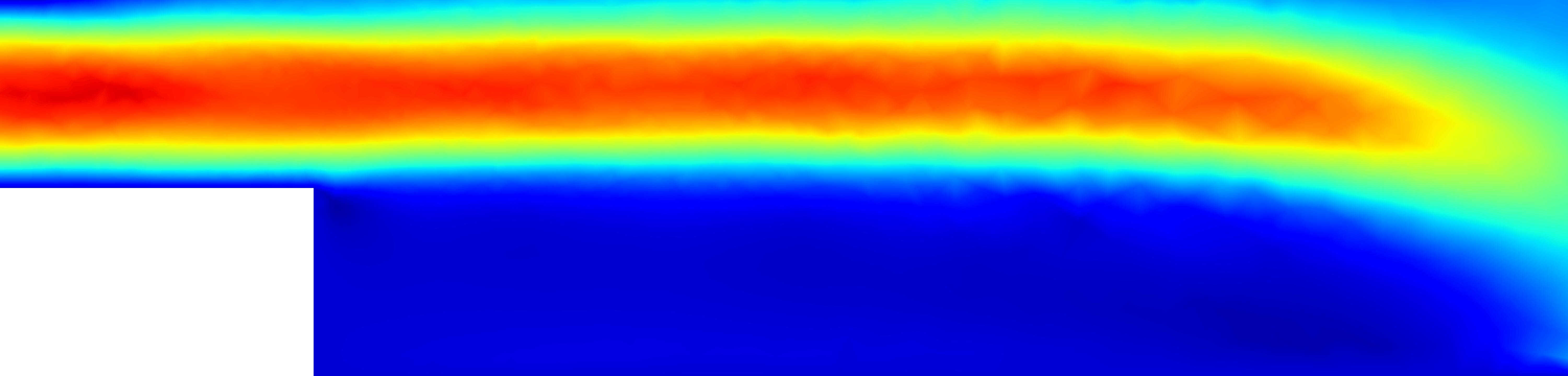}              \\
      (a) Reference solution, $t = 240\,\mathrm{s}$
      & (b) Reference solution, $t = 280\,\mathrm{s}$                           \\[6pt]
      \includegraphics[width=0.48\textwidth]{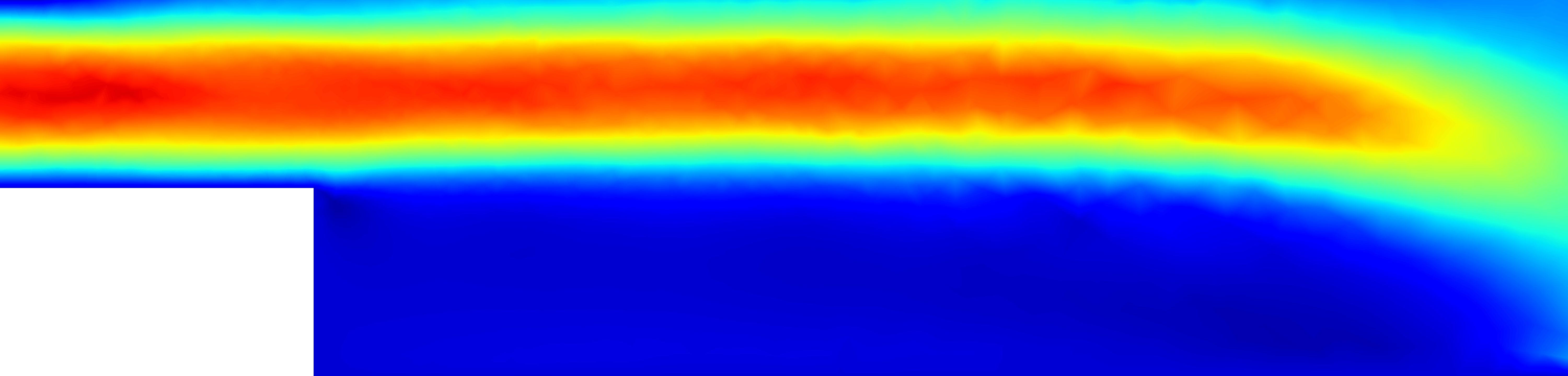}
      & \includegraphics[width=0.48\textwidth]{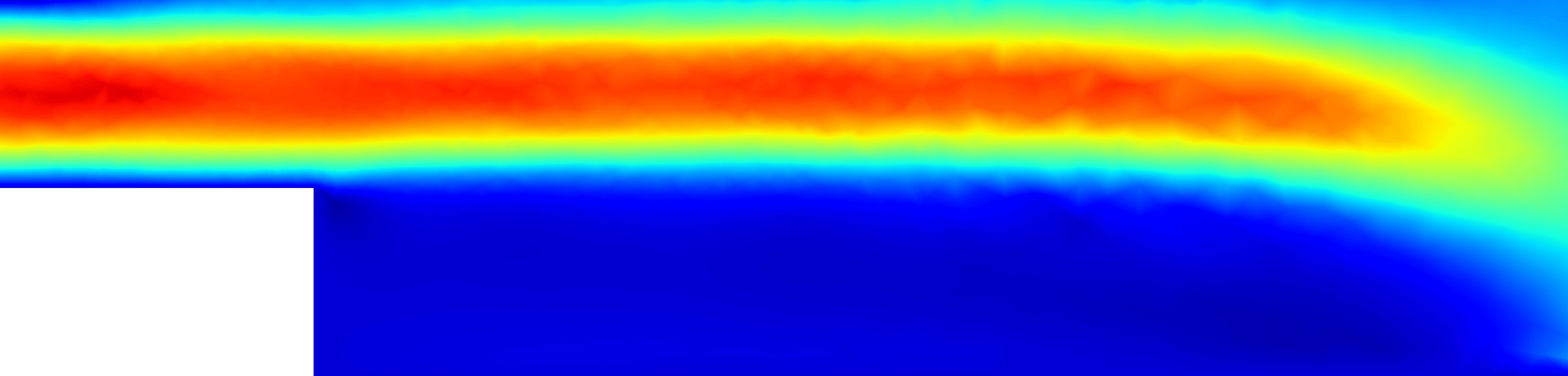}     \\
      (c) Gappy PMD + QDEIM, $t = 240\,\mathrm{s}$
      & (d) Gappy PMD + QDEIM, $t = 280\,\mathrm{s}$                            \\[6pt]
      \includegraphics[width=0.48\textwidth]{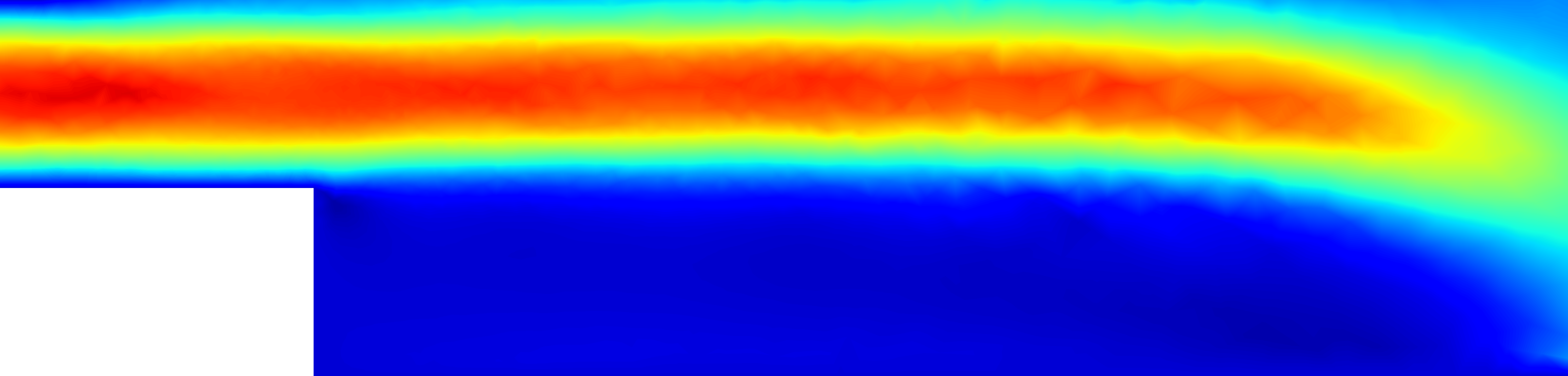}
      & \includegraphics[width=0.48\textwidth]{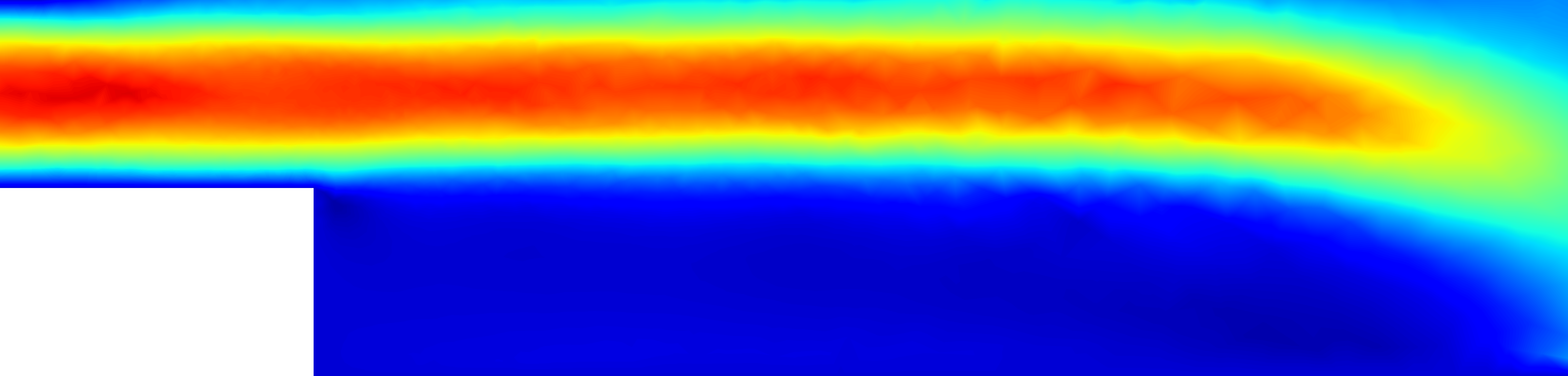}     \\
      (e) Gappy POD + QDEIM, $t = 240\,\mathrm{s}$
      & (f) Gappy POD + QDEIM, $t = 280\,\mathrm{s}$                            \\[2pt]
      \includegraphics[width=0.48\textwidth]{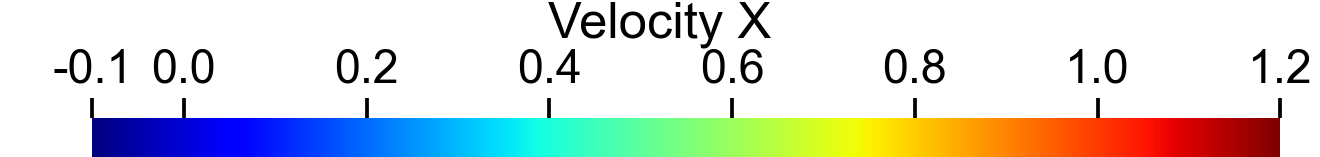}
      & \includegraphics[width=0.48\textwidth]{backward_colorbar1.png}           \\[6pt]
      \includegraphics[width=0.48\textwidth]{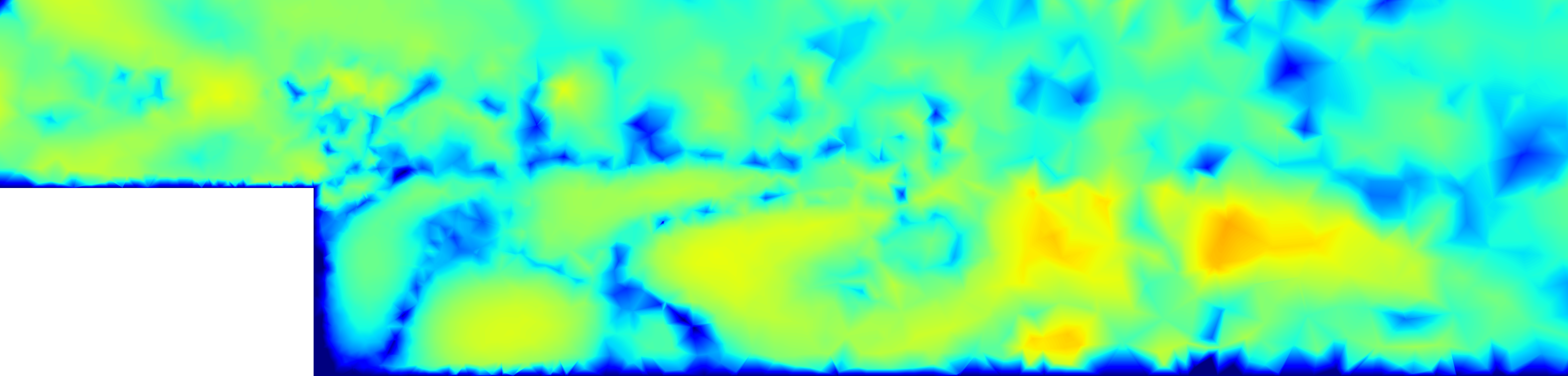}
      & \includegraphics[width=0.48\textwidth]{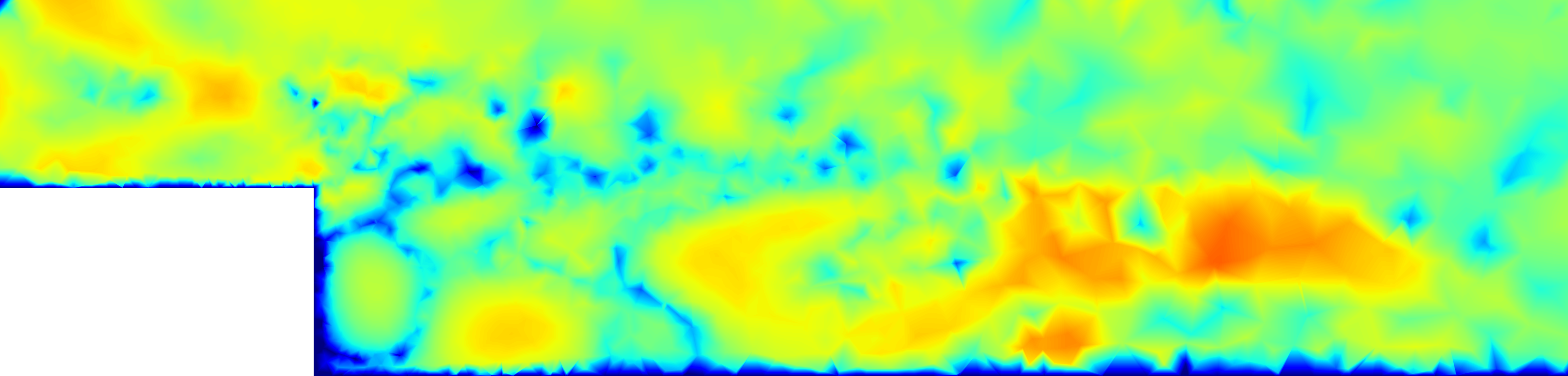} \\
      (g) Error of Gappy PMD + QDEIM, $t = 240\,\mathrm{s}$
      & (h) Error of Gappy PMD + QDEIM, $t = 280\,\mathrm{s}$                   \\[6pt]
      \includegraphics[width=0.48\textwidth]{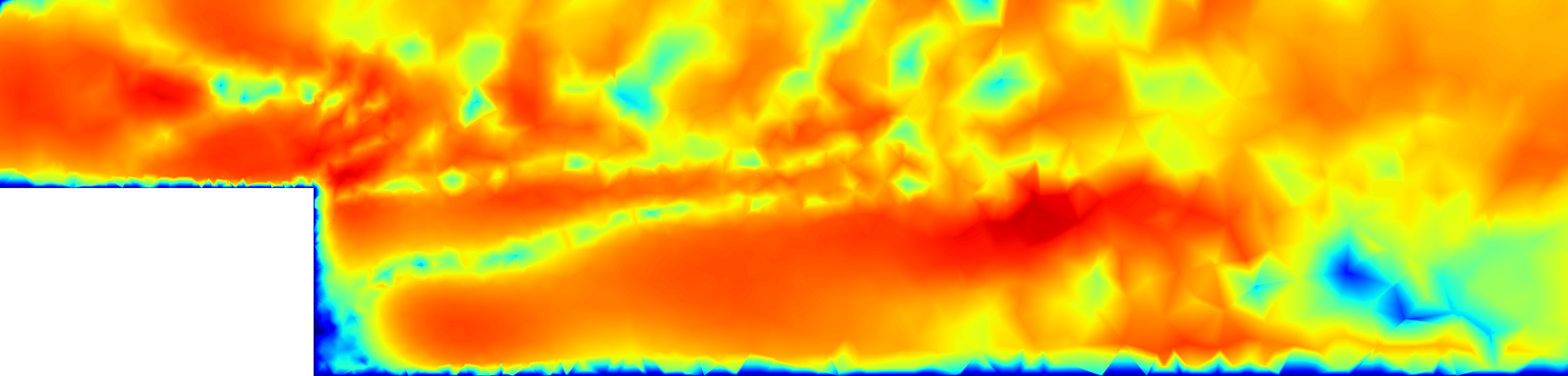}
      & \includegraphics[width=0.48\textwidth]{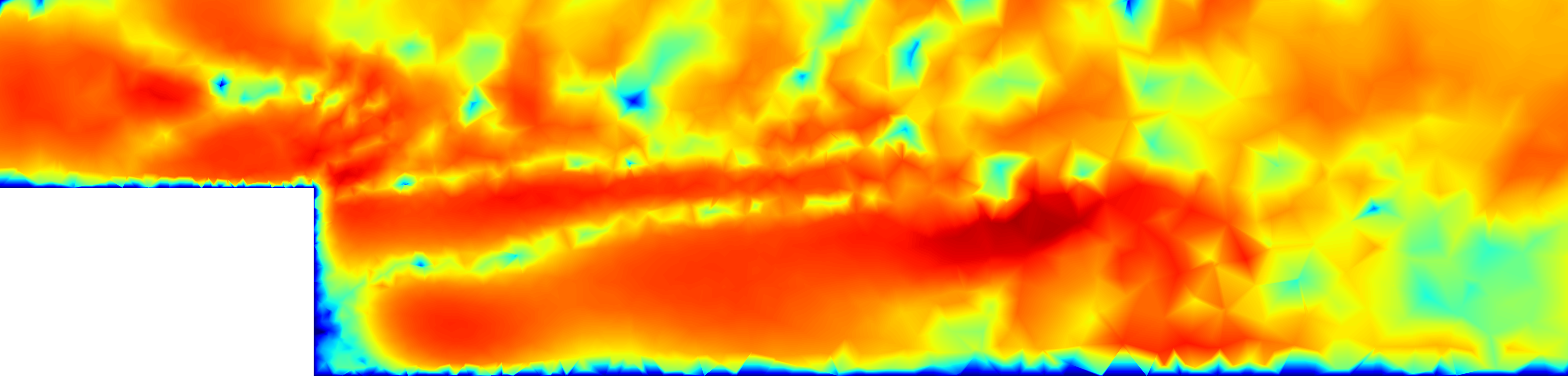} \\
      (i) Error of Gappy POD + QDEIM, $t = 240\,\mathrm{s}$
      & (j) Error of Gappy POD + QDEIM, $t = 280\,\mathrm{s}$                   \\[2pt]
      \includegraphics[width=0.48\textwidth]{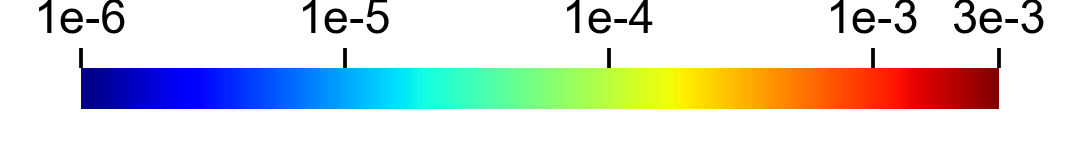}
      & \includegraphics[width=0.48\textwidth]{backward_err_colorbar1.png}       \\
   \end{tabular}
   \caption{Reconstructed streamwise velocity $u_x$ and the corresponding absolute error
   fields for the backward-facing step flow, comparing Gappy PMD and Gappy POD.}
   \label{fig:backward_recon_comparison}
\end{figure}

Figure~\ref{fig:backward_recon_comparison} compares Gappy PMD and Gappy POD
using the same nine QDEIM sampling points. The reconstructed velocity fields
are close to the reference fields at the scale shown. The error fields provide
a more sensitive comparison. Gappy PMD yields small reconstruction errors
throughout most of the domain. In contrast, Gappy POD exhibits elevated errors
over a wider region. This confirms the superior accuracy of the Gappy PMD
reconstruction for this test case.

\begin{figure}[htbp]
   \centering
   \begin{tabular}{cc}
      \includegraphics[width=0.48\textwidth]{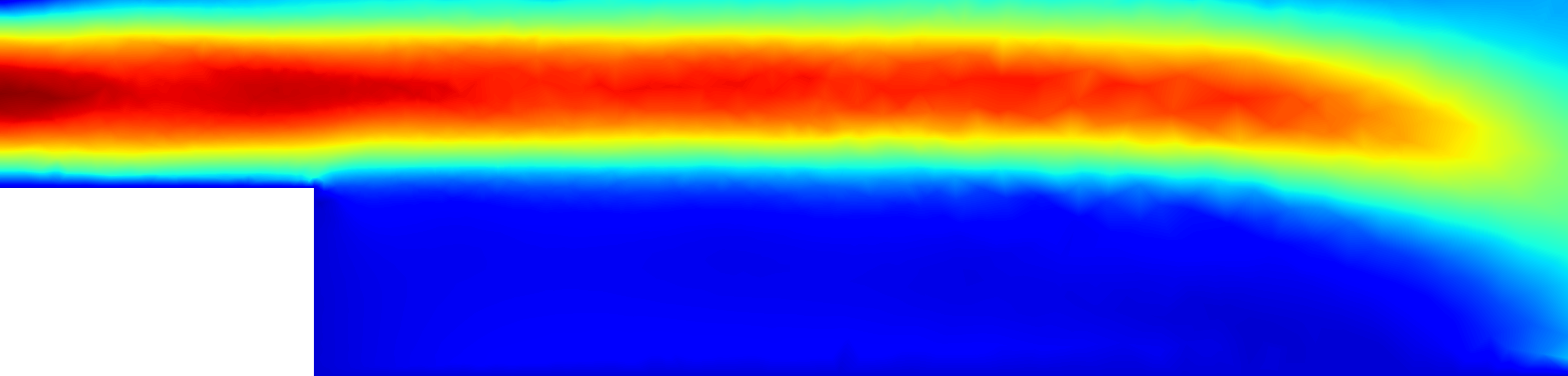}
      & \includegraphics[width=0.48\textwidth]{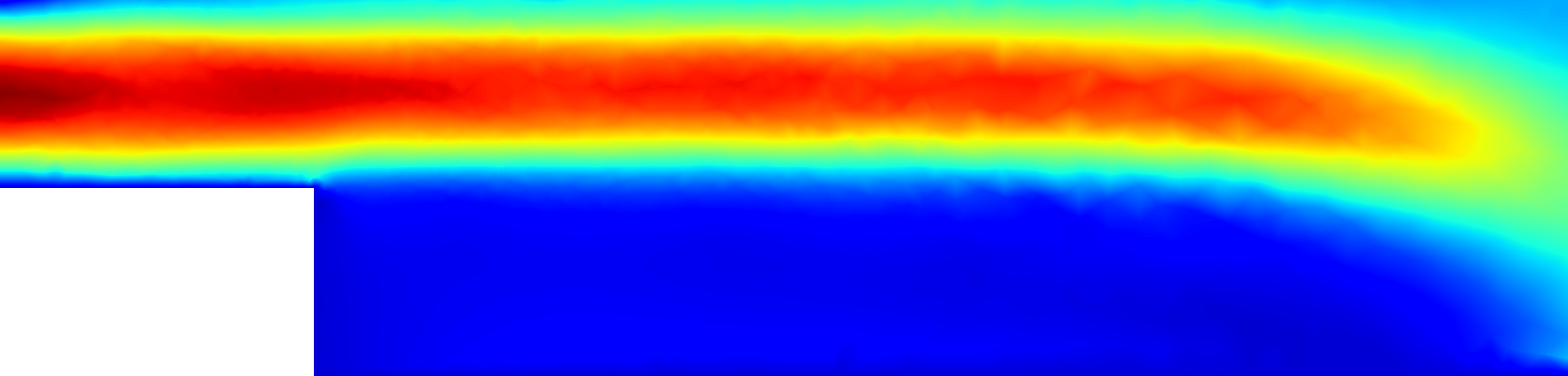}              \\
      (a) Reference solution, $t = 270\,\mathrm{s}$
      & (b) Reference solution, $t = 310\,\mathrm{s}$                           \\[6pt]
      \includegraphics[width=0.48\textwidth]{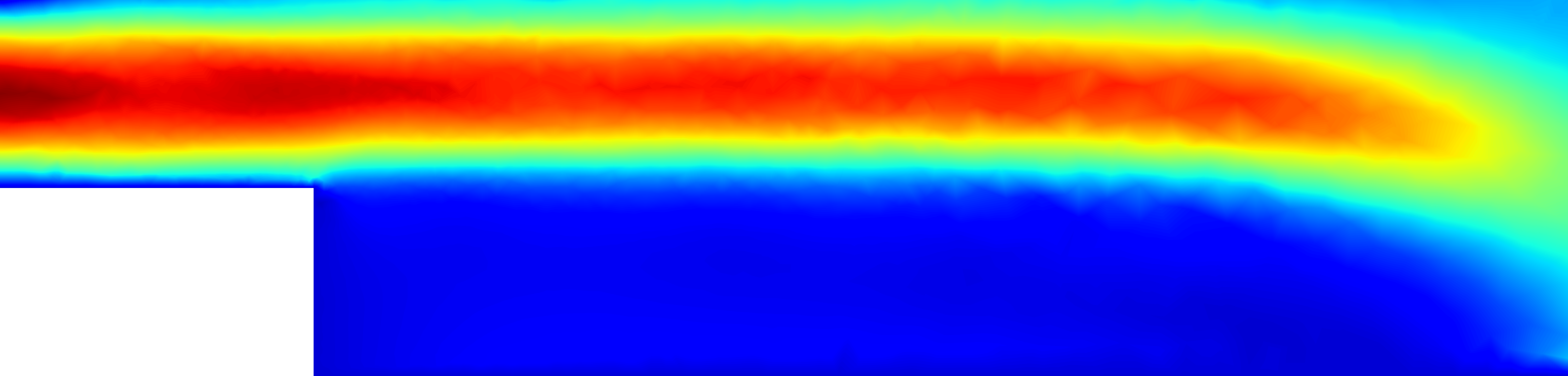}
      & \includegraphics[width=0.48\textwidth]{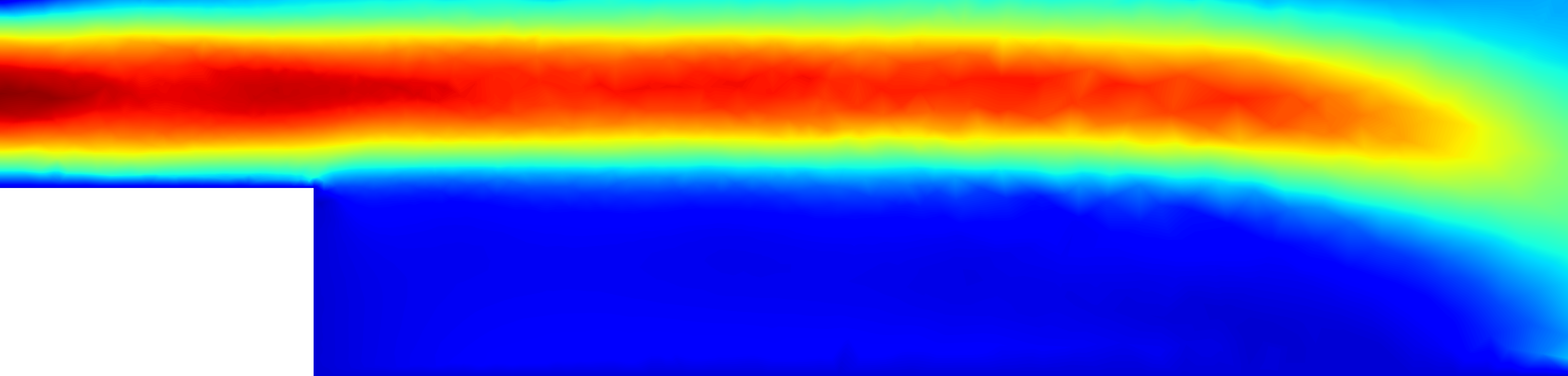}      \\
      (c) Gappy PMD + DPS, $t = 270\,\mathrm{s}$
      & (d) Gappy PMD + DPS, $t = 310\,\mathrm{s}$                             \\[6pt]
      \includegraphics[width=0.48\textwidth]{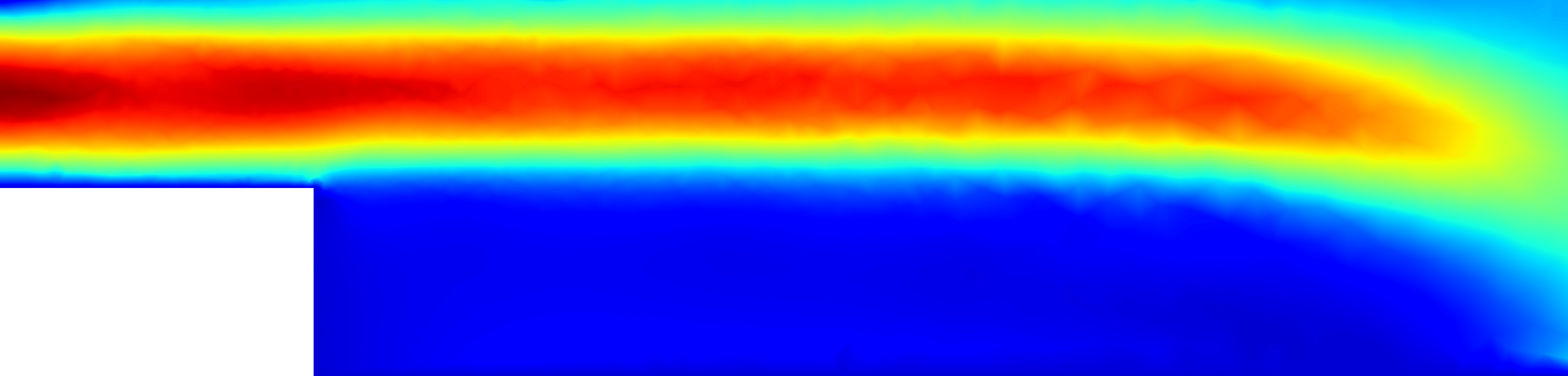}
      & \includegraphics[width=0.48\textwidth]{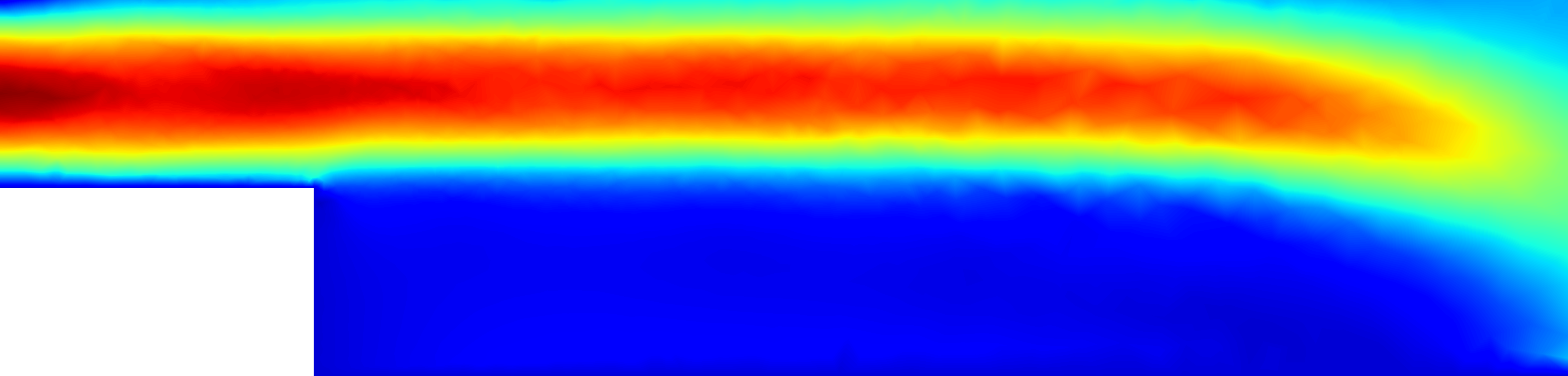}     \\
      (e) Gappy PMD + QDEIM, $t = 270\,\mathrm{s}$
      & (f) Gappy PMD + QDEIM, $t = 310\,\mathrm{s}$                            \\[2pt]
      \includegraphics[width=0.48\textwidth]{backward_colorbar1.png}
      & \includegraphics[width=0.48\textwidth]{backward_colorbar1.png}           \\[6pt]
      \includegraphics[width=0.48\textwidth]{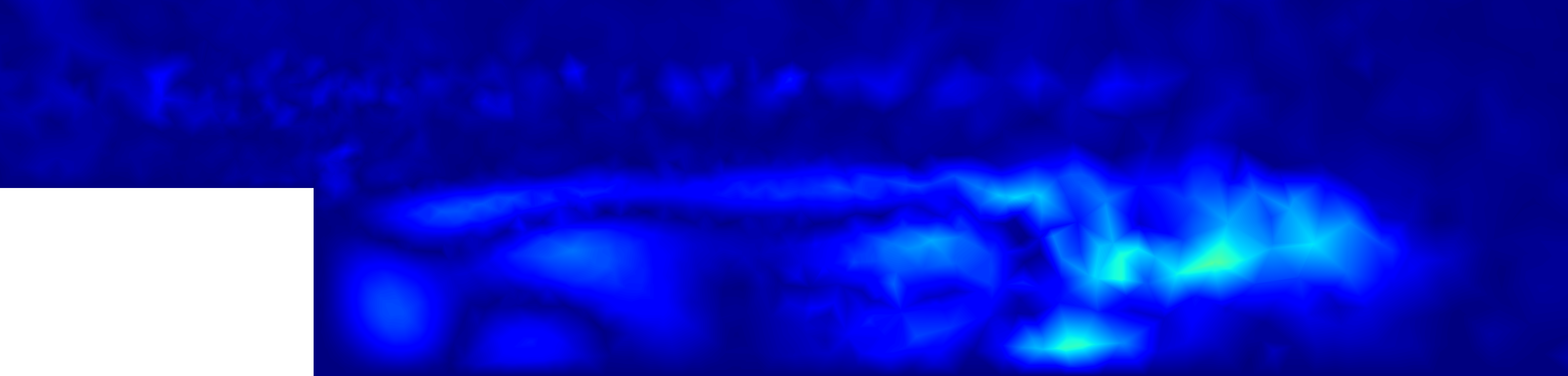}
      & \includegraphics[width=0.48\textwidth]{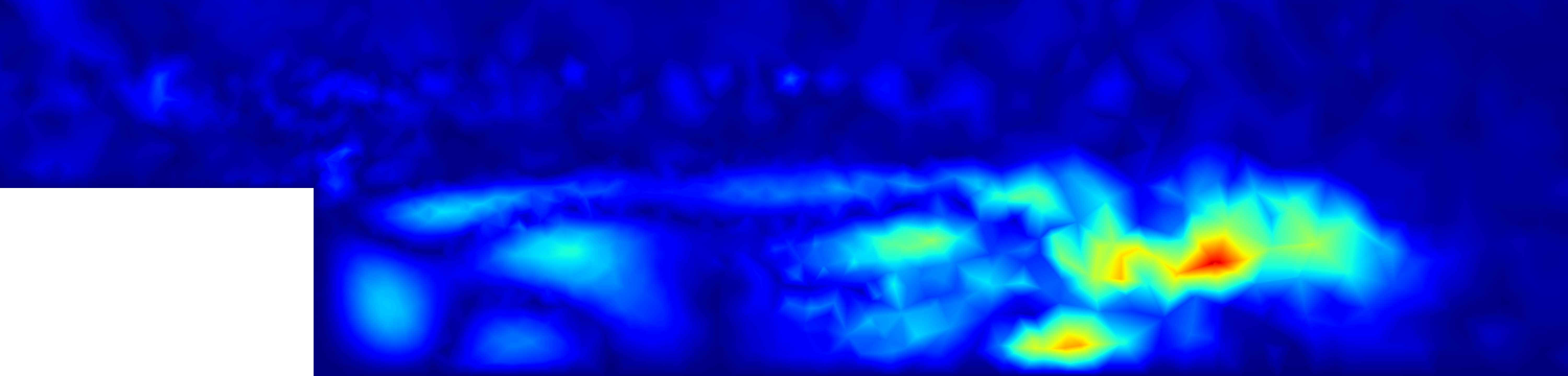}  \\
      (g) Error of Gappy PMD + DPS, $t = 270\,\mathrm{s}$
      & (h) Error of Gappy PMD + DPS, $t = 310\,\mathrm{s}$                    \\[6pt]
      \includegraphics[width=0.48\textwidth]{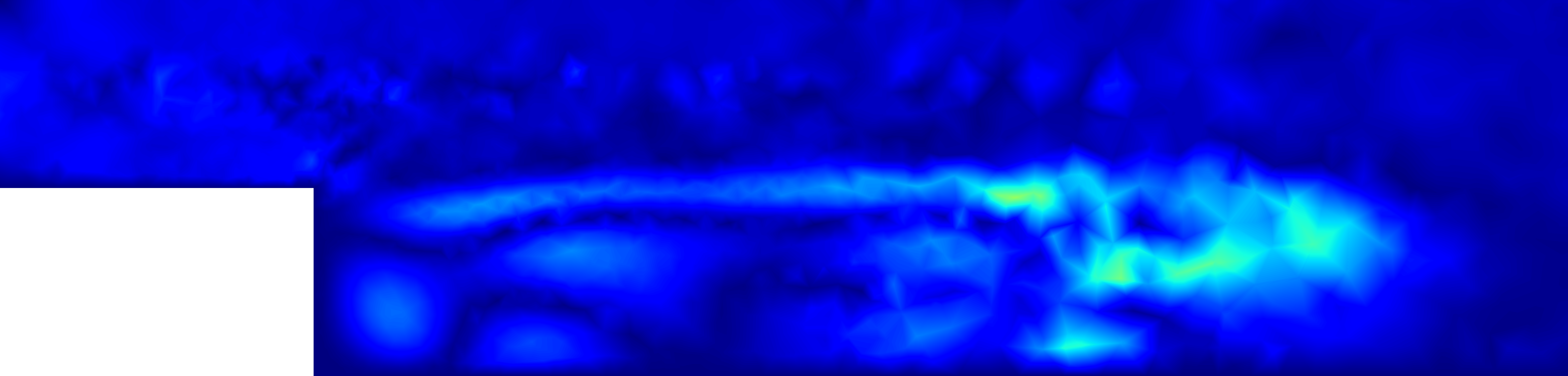}
      & \includegraphics[width=0.48\textwidth]{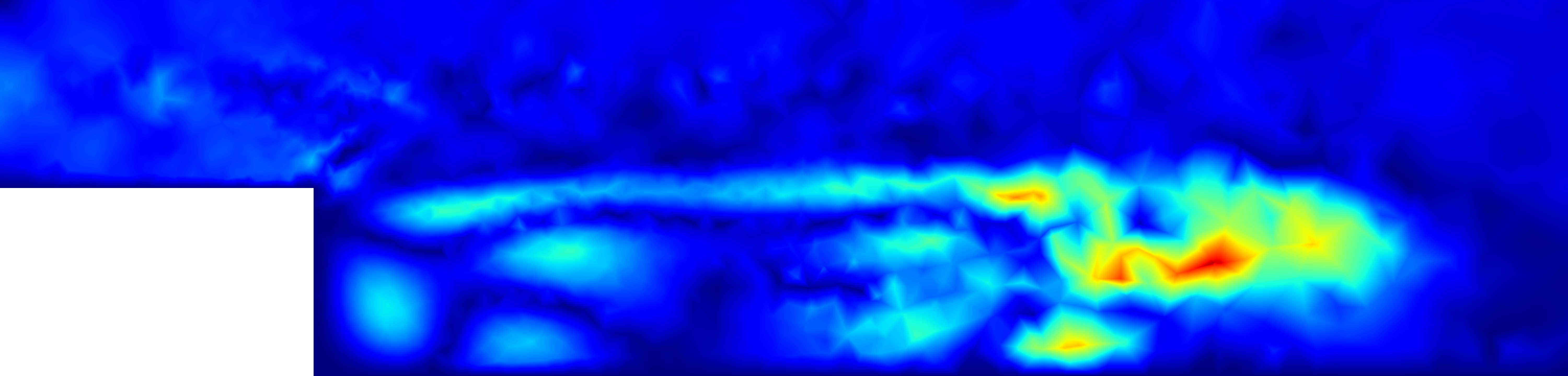} \\
      (i) Error of Gappy PMD + QDEIM, $t = 270\,\mathrm{s}$
      & (j) Error of Gappy PMD + QDEIM, $t = 310\,\mathrm{s}$                   \\[2pt]
      \includegraphics[width=0.48\textwidth]{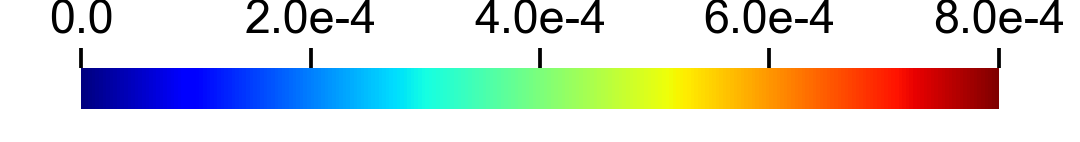}
      & \includegraphics[width=0.48\textwidth]{backward_err_colorbar2.png}       \\
   \end{tabular}
   \caption{Reconstructed streamwise velocity $u_x$ and the corresponding absolute error
   fields for the backward-facing step flow, comparing QDEIM and DPS sampling under Gappy
   PMD.}
   \label{fig:backward_sampling_comparison}
\end{figure}

Figure~\ref{fig:backward_sampling_comparison} compares QDEIM and DPS sampling
within the Gappy PMD reconstruction. The two reconstructions are
indistinguishable, but their error fields differ in amplitude. DPS yields
fewer regions of high error than QDEIM. This reduction is observed at both
instants and is obtained without increasing the number of sampling points.

\begin{figure}[htbp]
   \centering
   \includegraphics[width=0.6\textwidth]{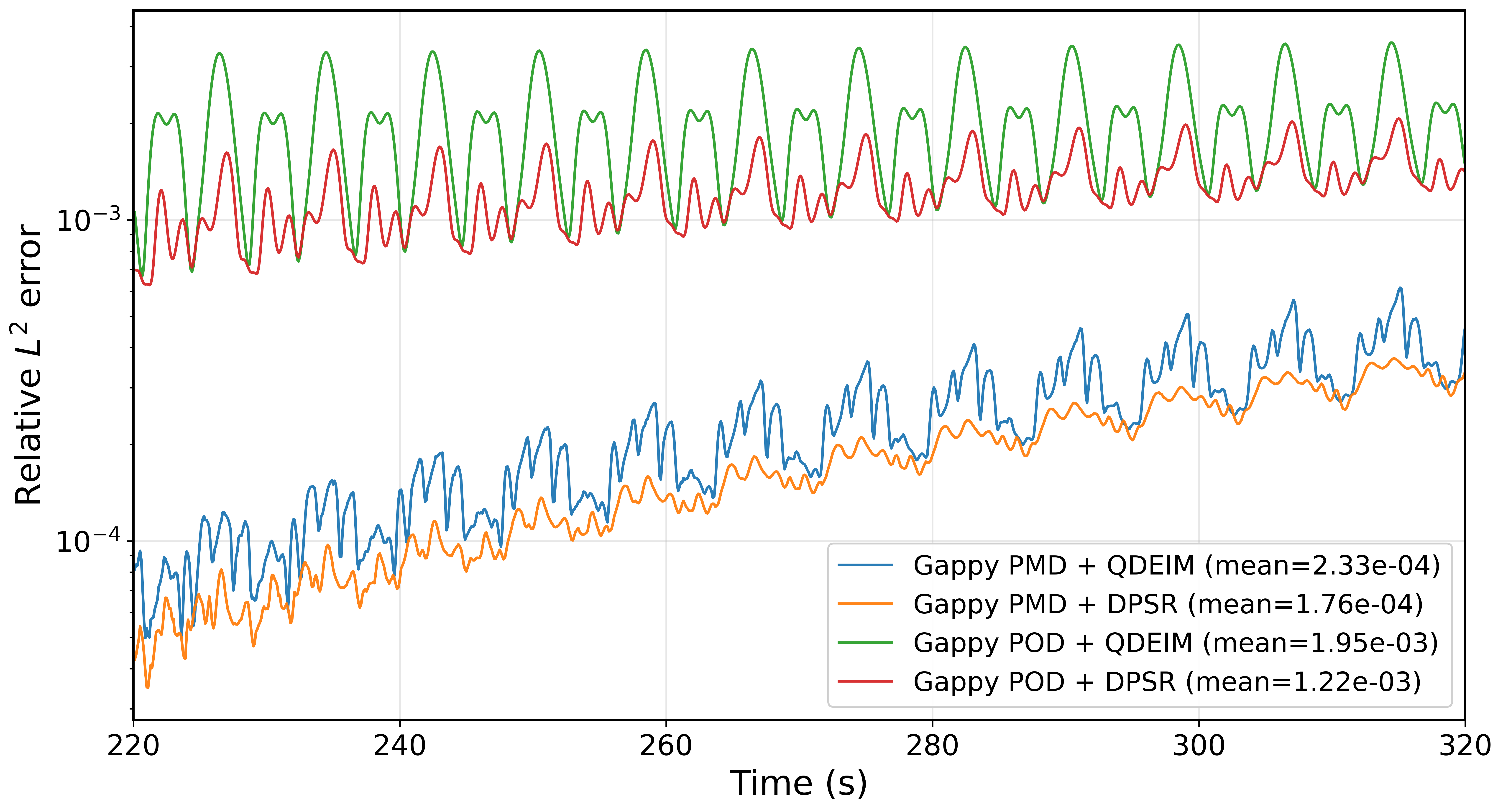}
   \caption{Relative $L^2$ reconstruction error over the test snapshots for the
   backward-facing step flow, for the four combinations of reconstruction method and
   sampling strategy.}
   \label{fig:backward_error_curve}
\end{figure}

Figure~\ref{fig:backward_error_curve} shows the relative $L^2$ error over the
test snapshots for Gappy PMD and Gappy POD, each evaluated with QDEIM and DPS
sampling. The reconstruction method has the stronger influence. For either
sampling strategy, the mean error of Gappy PMD is about one order of magnitude
lower than that of Gappy POD. The choice of sampling points has a smaller but
consistent effect. Replacing QDEIM by DPS reduces the mean error by about one
quarter for Gappy PMD and by slightly more than one third for Gappy POD. Both
comparisons reflect a consistent trend over the test interval, not an
instantaneous effect.

\begin{figure}[htbp]
   \centering
   \includegraphics[width=0.65\textwidth]{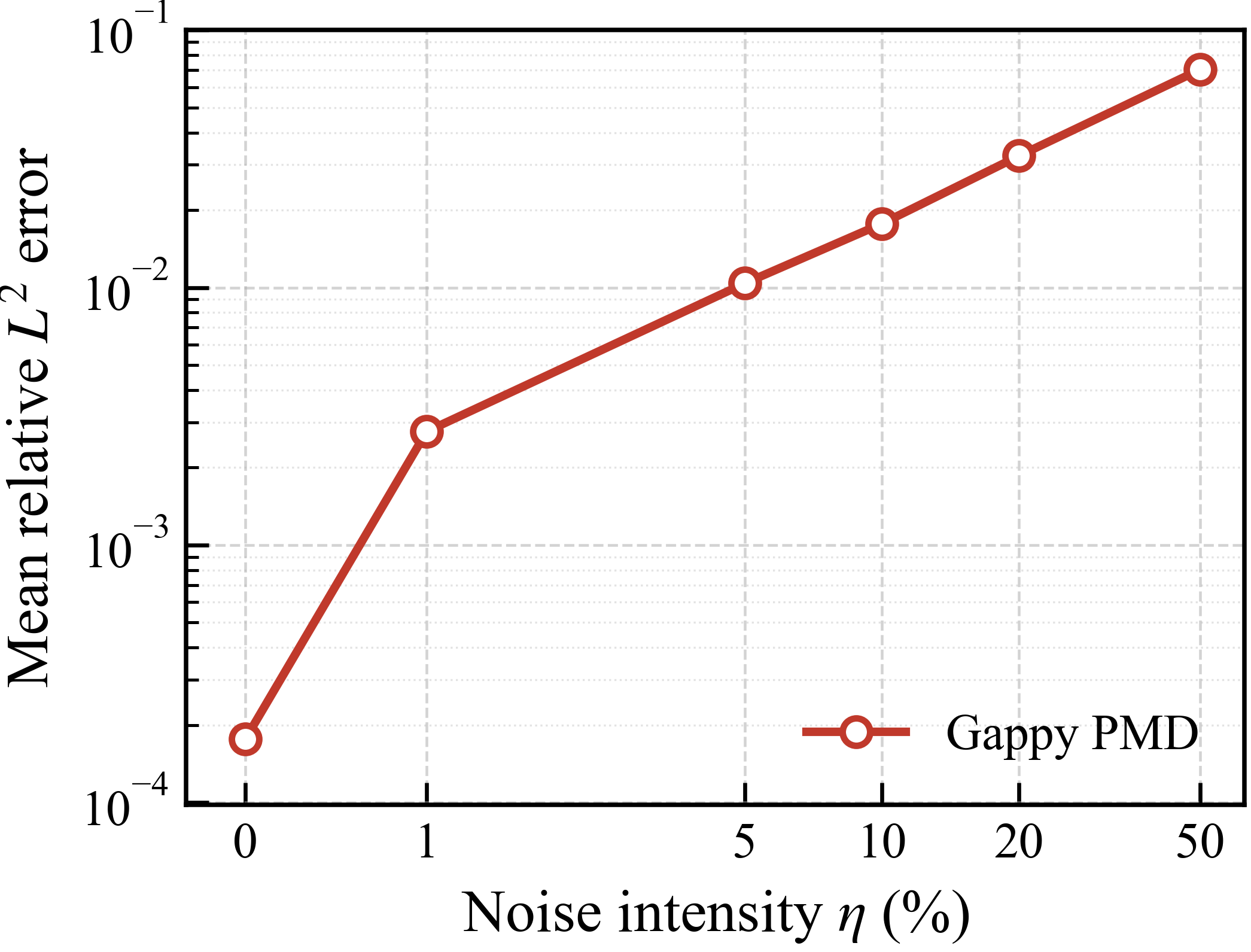}
   \caption{Effect of the noise level $\eta_\%$ on the mean relative $L^2$ error of Gappy PMD
   with DPS sampling for the backward-facing step flow, averaged over the test snapshots
   and over five independent noise realizations.}
   \label{fig:backward_noise_robustness}
\end{figure}

Figure~\ref{fig:backward_noise_robustness} isolates the effect of measurement
noise for the Gappy PMD reconstruction with DPS sampling. As soon as noise is
introduced, the noise contribution dominates the reconstruction error. The mean
error then increases approximately linearly
with the noise level, and this linear trend is maintained up to the largest
noise level considered, showing no accelerated degradation under strong
corruption.

\section{Conclusions}\label{Conclusions}
This paper proposes Gappy PMD, a method for reconstructing high-dimensional fields from
extremely sparse measurements. PMD provides an efficient low-dimensional representation of
complex fields with very few coordinates. Given sparse measurements, Gappy PMD recovers
the full field by solving a nonlinear least-squares problem for the PMD coordinates and
then applying the PMD lift mapping. By reconstructing the field on a nonlinear manifold rather than within a linear subspace, Gappy PMD can recover nonlinear solution structures that are difficult to capture with a low-dimensional linear reconstruction space.

We further propose DPS for selecting the sampling points. Classical criteria such as DEIM
and QDEIM select points from the row structure of a linear basis. DPS instead minimizes
the full-field reconstruction error directly. The discrete indices are relaxed to
continuous coordinates over the physical domain. The state is evaluated at these
coordinates by a differentiable RBF-FD interpolation operator. The reconstruction error
then varies smoothly with the sampling points, and DPS optimizes them by descending this
error directly. The optimized coordinates are finally snapped to distinct mesh nodes.

The theoretical analysis provides an error characterization
for Gappy PMD. In particular, the squared reconstruction error
can be decomposed orthogonally into the approximation error
associated with the reconstruction manifold and an additional
term induced by sparse sampling and measurement noise. Under a
suitable stability condition on the sampling operator, this
characterization further implies that the reconstruction error
vanishes as the PMD approximation error diminishes. The
resulting theory also reduces to the classical Gappy POD error
estimate when the reconstruction relies solely on the linear
basis expansion.

Gappy PMD and DPS were evaluated on three incompressible benchmarks. These are the
two-dimensional flow past a cylinder, the lid-driven cavity flow and the three-dimensional
backward-facing step flow. At equal total reduced dimension and sampling budget, Gappy PMD
attains mean relative $L^2$ errors one to two orders of magnitude smaller than Gappy POD.
Across the three cases, replacing QDEIM by DPS reduces this error by a further 25\% to
90\%. For Gappy POD, where QDEIM already provides a good point selection, DPS still improves
the mean error, but by a smaller margin of 11\% to 37\%. Under noisy measurements the error grows approximately linearly with the noise
level and does not break down even under severe corruption.

Gappy PMD is developed here for a fixed configuration, since the representation and the
sampling points are built offline from data at a single setting. Extending it to the
parametric setting is therefore a natural direction, so that the same representation and
sampling points remain effective as the Reynolds number or the boundary conditions vary
continuously.
\vspace{-6pt}

\section*{Acknowledgments}
\noindent
The authors acknowledge the support of the Fundamental Research Funds for
the Central Universities, the Top Discipline Plan of Shanghai Universities-Class
I, Shanghai Gaofeng Project for University Academic Program Development, National
Key R\&D Program of China(NO. 2022YFE0208000, 2024YFC2816400, and
2024YFC2816401). This work is also supported in part by grants from the
Shanghai Engineering Research Center (No.19DZ2255100) and the Shanghai
Institute of Intelligent Science and Technology, Tongji University. The computations were partially done on the deep learning computers of School of Mathematical Sciences, Tongji University, whose support is gratefully acknowledged.

\clearpage
\bibliographystyle{unsrt}
\bibliography{bibliography}

\end{document}